\documentclass[11pt]{article}
\usepackage{amsmath,amssymb}
\usepackage{graphicx}

\newcommand{\qed}{{\unskip\nobreak\hfil\penalty50\hskip2em\vadjust{}
       \nobreak\hfil$\Box$\parfillskip=0pt\finalhyphendemerits=0\par}}

\newtheorem{theorem}{Theorem}[section]
\newtheorem{lemma}{Lemma}[section]
\newtheorem{definition}{Definition}[section]
\newtheorem{cor}{Corollary}
\newtheorem{prop}{Proposition}
\newtheorem{conjecture}{Conjecture}

\newcommand{{\Z}}{\mathbb Z}
\newcommand{\R}{\mathbb R}
\newcommand{\N}{\mathbb N}

\newcommand{\dist}{\vert \vert}

\renewcommand{\P} {{\mathcal P}}

\newcommand{\cir}{\Gamma_0}
\newcommand{\reg}{{\rm RG}\big( \Gamma_0 \big)}
\newcommand{\reggam}{{\rm RG}\big( \Gamma \big)}
\newcommand{\conv}{{\rm conv}\big( \Gamma_0 \big)}
\newcommand{\delconv}{\partial {\rm conv}\big( \Gamma_0 \big)}
\newcommand{\mfl}{{\rm MFL}\big( \Gamma_0 \big)}
\newcommand{\mlr}{{\rm MLR}\big( \Gamma_0 \big)}
\newcommand{\mlrs}{{\rm MLRF}\big( \Gamma_0 \big)}
\newcommand{\maxseprg}{{\rm MPRG}\big( \Gamma_0 \big)}
\newcommand{\vcir}{V\big( \Gamma_0 \big)}

\newcommand{\acon}{{\rm AREA}_{\bo{0},n^2}}

\newcommand{\aconv}{\big\vert {\rm INT}\big( \Gamma \big) \big\vert \geq n^2}
\newcommand{\intg}{{\rm INT}\big( \Gamma_0 \big)}
\newcommand{\exc}{{\rm EXC}\big( \Gamma_0 \big)}

\newcommand{\globdis}{{\rm GD}\big( \Gamma_0 \big)}
\newcommand{\mc}{\mathcal}
\newcommand{\bo}{\mathbf}
\newcommand{\argu}{{\rm arg}}

\newcommand{\axy}{A_{\bo{x},\bo{y}}}
\newcommand{\axye}{E\big(A_{\bo{x},\bo{y}}\big)}
\newcommand{\axyo}{A_{\bo{x_0},\bo{y_0}}}
\newcommand{\axyoe}{E\big(A_{\bo{x_0},\bo{y_0}}\big)}
\newcommand{\xo}{\bo{x_0}}
\newcommand{\yo}{\bo{y_0}}

\newcommand{\zoz}{\{0,1 \}^{E(\Z^2)}}
\newcommand{\aarg}[2]{A_{\bo{#1},\bo{#2}}}

\newcommand{\ang}{\angle}
\newcommand{\qzero}{q_0}
\newcommand{\clu}[1]{C_{\pi/2 - \qzero}^F \big( #1 \big)}
\newcommand{\clum}[1]{C_{\pi/2 - \qzero}^B \big( #1 \big)} 
\newcommand{\clus}[2]{C_{\pi/2 - \qzero}^F \big( #1 \big) \cup C_{\pi/2 - \qzero}^B \big( #2 \big)}

\newcommand{\lmscb}{{\rm HIGHMFL}}

\newcommand{\mar}{\theta_{\rm RG}^{\rm MAX}\big(\cir\big)}
\newcommand{\margam}{\theta_{\rm RG}^{\rm MAX}\big(\Gamma\big)}

\newcommand{\cco}[1]{W_{#1,c_0}\big( \bo{0} \big)}
\newcommand{\xon}{\bo{x_1}}
\newcommand{\xtw}{\bo{x_2}}

\newcommand{\hvert}{h}
\newcommand{\xmlr}{\bo{x_{\rm MLR}}}
\newcommand{\xmlrm}{\bo{x_{\rm MLR}^-}}
\newcommand{\xmlrp}{\bo{x_{\rm MLR}^+}}

\newcommand{\ccone}{c_1}
\newcommand{\cctwo}{C_1}
\newcommand{\ctilde}{\tilde{C}}
\newcommand{\cprime}{C'}
\newcommand{\ctwo}{C_3}
\newcommand{\cthr}{C_4}
\newcommand{\cfour}{C_5}

\newcommand{\cgac}{2}

\newcommand{\cgen}{c}
\newcommand{\cgeno}{\hat{C}}
\newcommand{\cgentw}{C}
\newcommand{\conka}{C_{\rm rwm}}
\newcommand{\csmall}{c}

\newcommand{\clemgac}{C_{\rm gac}}
\newcommand{\clemkac}{C_2}
\newcommand{\clemkam}{m_0}
\newcommand{\crwm}{C}
\newcommand{\cposen}{c_{\rm be}}
\newcommand{\cpi}{20\pi}
\newcommand{\cgenbig}{C_*}
\newcommand{\perpu}[1]{#1^{\perp}}

\newcommand{\centre}{{\rm cen}}

\newcommand{\area}[1]{{\rm AREA}_{\bo{0},#1}}
\newcommand{\areacon}{\big\vert \intg \big\vert \geq n^2}
\newcommand{\vsat}{{\rm SAT}^+}
\newcommand{\wulff}{\mathcal{W}_\beta}
\newcommand{\outfluc}{{\rm OUTFLUC}}
\newcommand{\aaa}{F}
\newcommand{\bb}{G}
\newcommand{\cluh}[1]{C_{\pi/2 - \qzero/2}^F \big( #1 \big)}
\newcommand{\clumh}[1]{C_{\pi/2 - \qzero/2}^B \big( #1 \big)}
\newcommand{\gacc}{{\rm GAC}}

\def\build#1_#2^#3{\mathrel{ \mathop{\kern 0pt#1}\limits_{#2}^{#3}}}

\def\fff#1{&{{\pageref{#1}}}\cr}
\def\hfff#1{\label{#1}}

\begin{document}
\title{Phase separation in random cluster models I: \\ uniform upper bounds on local deviation}
\author{Alan Hammond\thanks{Department of Statistics, University of Oxford. Supported in part by EPSRC grant EP/I004378/1. This work was undertaken during visits to the Theory Group at Microsoft Research in Redmond, WA, and to Ecole Normale Superieure in Paris. 2010 Mathematics Subject Classification: 60D05, 82B41.}} 
 \maketitle

\begin{abstract}
This is the first in a series of three papers that addresses the behaviour of the droplet that results, in the percolating phase, from conditioning the planar Fortuin-Kasteleyn random cluster model on the presence of an open dual circuit $\cir$ encircling the origin and enclosing an area of at least (or exactly) $n^2$. (By the Fortuin-Kasteleyn representation, the model is a close relative of the droplet formed by conditioning the Potts model on an excess of spins of a given type.)
We consider local deviation of the droplet boundary, measured in a radial sense by the maximum local roughness, $\mlr$, this being the maximum distance from a point in the circuit $\cir$ to the boundary $\delconv$ of the circuit's convex hull; and in a longitudinal sense by what we term {\it maximum facet length}, $\mfl$, namely, the length of the longest line segment of which the polygon $\delconv$ is formed. 
The principal conclusion of the series of papers is the following uniform control on local deviation: that there are constants $0 < c < C < \infty$
such that the conditional probability that the normalized quantity 
$n^{-1/3}\big( \log n \big)^{-2/3} \mlr$ 
lies  in the interval $\big[c,C]$ 
tends to $1$ in the high $n$-limit; and that the same statement holds for  
$n^{-2/3}\big( \log n \big)^{-1/3}  \mfl$. In this way, we confirm the anticipated $n^{1/3}$ scaling of maximum local roughness, and provide a sharp logarithmic power-law correction. 
This local deviation behaviour occurs by means of locally Gaussian effects constrained globally by curvature, and we believe that it arises in many radially defined stochastic interface models, including growth models belonging to the Kardar-Parisi-Zhang universality class. 

The present paper is devoted to proving the upper bounds in these assertions. In fact, we derive bounds valid in the moderate deviations' regime. The second paper \cite{hammondtwo} provides the lower bounds. Crucial to our approach are surgical techniques that renew the conditioned circuit on the scale at which the local deviation manifests itself. 
A successful analysis of the surgeries depends on the circuit possessing a renewal structure, with only local backtracking occurring from its overall progress in a direction specified by the macroscopic Wulff profile. The third paper \cite{hammondthr}
presents the required tool on regeneration structure of the conditioned circuit. 
 
The present paper includes a heuristic presentation of the surgical technique that is also used in \cite{hammondthr}, and a discussion of the significance of local deviation and of problems raised by our approach.
\end{abstract}

\setlength{\baselineskip}{16pt}
\newpage 

\subsection{Glossary of notation}
A certain amount of notation is needed during the proofs in this article.
For the reader's convenience, we begin by listing much of this notation, and provide a summarizing phrase for each item, as well as the page number at which the concept is introduced.

\bigskip
\def\qq{&}

\begin{center}
\halign{
#\quad\hfill&#\quad\hfill&\quad\hfill#\cr
$\Gamma$ \qq a generic circuit \fff{gencir}  
${\rm INT}(\Gamma)$ \qq the region enclosed by $\Gamma$ \fff{intcir}
$\cir$ \qq  the outermost open circuit enclosing $\bo{0}$ \fff{outcir}
$\conv$ \qq the convex hull of $\cir$ \fff{conv}
$\mlr$ \qq  maximum local roughness \fff{mlr}    
$\mfl$ \qq  maximum facet length \fff{mfl}    
$\exc$ \qq the area-excess $\intg - n^2$ \fff{exc}
$T_{\bo{0},\bo{x},\bo{y}}$ \qq the triangle with vertices $\bo{0}$, $\bo{x}$ and $\bo{y}$ \fff{txy}
$\axy$ \qq  the sector with apex $\bo{0}$ bounded by $\bo{x}$ and $\bo{y}$ \fff{axy}   
$W_{\bo{v},c}$ \qq  the cone about $\bo{v}$ with apex $\bo{0}$ and aperture $2c$ \fff{wvc}
$W_{\bo{v},c}^{+/-}$ \qq  an alternative notation for specifying a cone \fff{wvcalt}
$\wulff$ \qq  the unit-area Wulff shape \fff{wulff}   
$\globdis$ \qq  global distortion (from the Wulff shape) \fff{globdis}
$\centre(\Gamma)$ \qq  the centre of a circuit \fff{centre}
$\area{A}$ \qq the event of capture of area $A$ by a circuit centred at the origin \fff{acon}
$C^F_{\pi/2 - \qzero} \big( \bo{v} \big)$ 
\qq the $\qzero$-forward cone with apex $\bo{v}$ \fff{forback}
$C^B_{\pi/2 - \qzero} \big( \bo{v} \big)$ 
\qq the $\qzero$-backward cone  with apex $\bo{v}$ \fff{forback}
 $\reg$ \qq  the set of circuit regeneration sites \fff{reg}
${\rm fluc}_{\bo{x},\bo{y}}\big( \gamma \big)$ \qq the maximal distance of a point in $\gamma$ from $[\bo{x},\bo{y}]$ \fff{fluc} 
$\sigma_{\bo{x},\bo{y}}$ \qq the sector storage-replacement operation (with sector $A_{\bo{x},\bo{y}}$) \fff{ssro}  
$\gamma_{\bo{x},\bo{y}}$ \qq  the outermost open path in $A_{\bo{x},\bo{y}}$ from $\bo{x}$ to $\bo{y}$ \fff{gamoop}   
 $\overline\gamma_{\bo{x},\bo{y}}$ \qq  the open cluster of  $\bo{x}$ and $\bo{y}$ in $A_{\bo{x},\bo{y}}$ \fff{ogamoop}
$I_{\bo{x},\bo{y}} \big( \gamma_{\bo{x},\bo{y}} \big)$ \qq the bounded component of $\axy \setminus \gamma_{\bo{x},\bo{y}}$ \fff{ixy} 
 ${\rm GAC} \big( \bo{x},\bo{y} ,\epsilon\big)$ \qq  configurations in $\axy$ realizing $\epsilon$-good area capture \fff{gac} 
$\phi_{\aaa,\bb,\bo{x}}$ \qq the storage-shift-replacement operation 
 \fff{stsh}
$\maxseprg$ \qq the maximum point-to-regeneration-site distance  \fff{mprg}
$\mlrs$ \qq the length of the facet associated to the site of maximum local roughness
\fff{mlrf}
}\end{center}

\begin{section}{Introduction}
The ferromagnetic Ising model may be considered as describing two populations, the members of each of which prefer the company of their own type. How do the two populations reside if they are forced to live in the same plot of land? Consider the Ising model in a large box, at a supercritical inverse temperature $\beta > \beta_c$ with negative boundary condition. There is a positive magnetization $m(\beta)$, so that   the interior is unfavourable to the populace of positive signs who will typically constitute a minority fraction $\big( 1 - m(\beta) \big)/2$ of the total sites of the box. Conditionally on the plus signs forming a fraction $\lambda > \big( 1 - m(\beta) \big)/2$, an enclave (or droplet) forms in which the plus signs form their own phase, representing there a majority 
fraction $\big( 1 + m(\beta) \big)/2$ of sites, in a surrounding environment in which the negative signs continue to occupy this same fraction of sites. The theory of phase separation is concerned with the study of the random boundary between phases, such as the droplet boundary in this example.

Wulff \cite{wulff} proposed that the profile of such constrained circuits would macroscopically resemble a dilation of an isoperimetrically optimal curve that now bears his name. For the Ising problem,
this claim was first verified rigorously in \cite{dks} at low temperature, 
the derivation being extended up to the critical temperature by \cite{ioffeschonmann}. Fluctuations from this profile have been classified into global or long-wave effects, measured by the deviation (in the Hausdorff metric) of the convex hull of the circuit from an optimally placed dilate of the Wulff shape. Local fluctuations have been measured by the inward deviation of the circuit from the boundary of its convex hull. 

As we will shortly explain, these local fluctuations arise from a competition between the locally Brownian evolution of the interface and a global constraint imposed by the requirement of area capture. The radial and longitudinal magnitudes of these fluctuations exist on the smallest scale at which the area capture requirement competes with the Gaussian path fluctuation. 
This suggests that such local fluctuations should be the object of our attention in seeking to understand the relationship between phase separation and the KPZ universality class \cite{kpz} of random growth models. Indeed, it has been demonstrated physically \cite{majcomt} that the law of the maximum height at equilibrium of certain interfaces subject to roughening and smoothing is given by the Airy distribution, the law of the area trapped by a unit-time one-dimensional Brownian excursion. This same distribution dictates the nature of competition of area capture and Brownian fluctuation in the scale where average and maximum local fluctuation arises in conditioned spin systems. Maximum fluctuations of constrained Brownian motions \cite{schledou} are an instance of the extreme value theory of correlated random variables, a theory which governs, for example, the equilibrium behaviour of low-temperature spin glasses \cite{boumez}.

In this article and its companions \cite{hammondtwo} and \cite{hammondthr}, we study the maximum local fluctuations of  the outermost open circuit enclosing a given point in a subcritical planar random cluster model, when this circuit is conditioned to entrap a large area $n^2$.  
As we will discuss after stating the principal results, the Fortuin-Kasteleyn representation and duality considerations render this problem a close relative to that of the behaviour of a large conditioned droplet in the Ising (or Potts) model arising from conditioning on an excess of signs of a given type. 
We determine sharp logarithmic corrections in the parameter $n$ for maximum local deviations, by means of a flexible new surgical technique. We expect that the scaling exponents here identified hold for a wide range of models for phase separation and in the KPZ class, and that the tools that we introduce will be adaptable for use in other conditioned spin systems and possibly also to kinetic KPZ models.

We now give the definitions required to formulate our results.
\begin{definition}
For $\Lambda \subseteq \Z^2$, let $E(\Lambda)$ denote the set of nearest-neighbour edges whose endpoints lie in $\Lambda$
and write $\partial_{\rm int} \big( \Lambda \big)$ for the interior vertex boundary of $\Lambda$, namely, the subset of $\Lambda$ each of whose elements is an endpoint of some element of 
 $E(\Lambda)^c$. Fix a choice of $\Lambda \subseteq \Z^2$ that is finite. 
The free random cluster model on $\Lambda$ with parameters $p \in [0,1]$ and $q > 0$ on $\Lambda$ 
is the probability space over $\eta: E(\Lambda) \to \{0,1\}$ with measure
$$
\phi_{p,q}^f(\eta) = \frac{1}{Z_{p,q}}   p^{\sum_e \eta(e)} 
\big( 1 - p \big)^{\sum_e (1 - \eta(e))} q^{k(\eta)},
$$
where $k(\eta)$ denotes the number of connected components in the subgraph of 
$\big(\Lambda,E(\Lambda)\big)$ containing all vertices and all edges $e$ such that $\eta(e) = 1$. (The constant $Z_{p,q}$ is a normalization.) The wired random cluster model 
$\phi_{p,q}^w$ is defined similarly, with $k(\eta)$ now denoting the number of such connected components none of whose edges touch $\partial_{\rm int} \big( \Lambda \big)$.

For parameter choices $p \in [0,1]$ and $q \geq 1$, 
either type of random cluster measure $\P$ satisfies the FKG inequality: suppose that 
$f,g:  \{0,1\}^{E(\Lambda)} \to \R$
are increasing functions with respect to the natural partial order on $\{0,1\}^{E(\Lambda)}$. 
Then $\mathbb{E}_\P \big( fg \big) \geq \mathbb{E}_\P \big( f \big) \mathbb{E}_\P \big( g \big)$,
where $\mathbb{E}_\P$ denotes expectation with respect to $\P$.

Consequently, we define the infinite-volume free and wired random cluster measures $\P^f$ and $\P^w$ as limits of the finite-volume counterparts taken along any increasing sequence of finite sets $\Lambda \uparrow \Z^2$. The measures $\P^f$ and $\P^w$ are defined on the space of functions $\eta: E(\Z^2) \to \{0,1\}$ with the product $\sigma$-algebra. In a realization $\eta$, the edges $e \in E(\Z^2)$
 such that $\eta(e) = 1$ are called open; the remainder are called closed.  A subset of $E(\Z^2)$ will be called open (or closed) if all of its elements are open (or closed).
We will record a realization in the form $\omega \in \zoz$, where 
the set of coordinates that are equal to $1$ under $\omega$ is the set of open edges under $\eta$.
Any $\omega \in \zoz$ will be called a configuration. For $\bo{x},\bo{y} \in \Z^d$, we write $\bo{x} \leftrightarrow \bo{y}$ to indicate that $\bo{x}$ and $\bo{y}$ lie in a common connected component of open edges.  
 
Set $\beta \in (0,\infty)$ according to $p = 1 - \exp\{ - 2 \beta \}$. 
In this way, the infinite volume measures are parameterized by $\P_{\beta,q}^w$ and $\P_{\beta,q}^f$ with $\beta > 0$ and $q \geq 1$.
For any $q \geq 1$, $\P^w_{\beta,q} = \P^f_{\beta,q}$ for all but at most countably many values of $\beta$ \cite{grimmett}. 
We may thus define
$$
\beta_c^1 = \inf \big\{  \beta > 0: \P^*_{\beta,q} \big(  0 \leftrightarrow \infty \big) > 0 \big\}
$$
obtaining the same value whether we choose $* = w$ or $* = f$.
\end{definition}
There is a unique random cluster model for each subcritical $\beta < \beta_c^1$ \cite{grimmett}, that we will denote by $\P_{\beta,q}$.
\begin{definition}
Let $\hat{\beta}_c$ denote the supremum over $\beta > 0$ such that the following holds: letting $\Lambda = \big\{ -N,\ldots,N \big\}^d$, there exist constants $C > c > 0$ such that, for any $N$,
$$
\P^w_{\beta,q} \Big( \bo{0} \leftrightarrow \Z^d \setminus  \Lambda_N \Big) \leq C \exp \big\{ - c N  \big\}.
$$
\end{definition}
In the two-dimensional case that is the subject of this article, it has been established that $\beta_c^1 = \hat{\beta}_c$ for $q=1$, $q=2$  and for $q$ sufficiently high, by \cite{alexmix}, and respectively \cite{ab}, \cite{abf} and \cite{lmmrs}. 
A recent advance \cite{beffaraduminilcopin} showed that, on the square lattice, in fact, $\beta_c^1 = \hat{\beta}_c$ holds for all $q \geq 1$. The common value, which is $2^{-1} \big( 1 + \sqrt{q} \big)$, we will denote by $\beta_c$.

The droplet boundary is now defined:
\begin{definition}
A circuit \hfff{gencir} $\Gamma$ is a nearest-neighbour path in $\Z^2$ whose endpoint coincides with its start point but that has no other self-intersections. We set $E(\Gamma)$ equal to the set of nearest-neighbour edges between successive elements of $\Gamma$. 
For notational convenience, when we write $\Gamma$, we refer to the closed subset of $\R^2$ given by the union of the topologically closed intervals corresponding to the elements of $E(\Gamma)$. 
We set $V(\Gamma) = \Gamma \cap \Z^2$.

Let $\omega \in \zoz$. A circuit $\Gamma$ is called open if $E(\Gamma)$ is open. 
For any circuit $\Gamma$, we write  \hfff{intcir} ${\rm INT} \big( \Gamma \big)$ for the bounded component of $\R^2 \setminus \Gamma$, that is, for the set of points enclosed by $\Gamma$.

An open circuit $\Gamma$ is called outermost if any open circuit $\Gamma'$ satisfying 
${\rm INT} \big( \Gamma \big) \subseteq {\rm INT} \big( \Gamma' \big)$ is equal to $\Gamma$. Note that, if
a point $\bo{z} \in \R^2$ is enclosed by a positive but finite number of open circuits in a configuration 
$\omega \in \zoz$ , it
is enclosed by a unique outermost open circuit.

We write  \hfff{outcir} $\cir$ for the outermost open circuit $\Gamma$ for which $\bo{0} \in {\rm INT}(\Gamma)$, taking $\cir = \emptyset$ if no such circuit exists.
\end{definition}
\noindent{\bf Remark.} 
Under any subcritical random cluster measure $P = \P_{\beta,q}$, with $\beta < \beta_c$, there is an exponential decay in distance for the probability that two points lie in the same open cluster. (See Theorem $A$ of \cite{civ}.) As such, $P$-a.s., no point in $\R^2$ is surrounded by infinitely many open circuits, so that $\cir$ exists (and is non-empty) whenever $\bo{0}$ is surrounded by an open circuit. 

Our object of study is the subcritical random cluster model given the event $\areacon$, with $n \in \N$ high.  We now present the radial and longitudinal notions of maximum local deviation in the droplet boundary that will be the principal objects of study. 
\begin{definition}
We write  \hfff{conv} $\conv$ for the convex hull of $V\big( \cir \big)$. 
The maximum local roughness $\mlr$ is defined to be the maximal distance of an element of $V(\cir)$
 to the boundary of the convex hull: that is \hfff{mlr}
$$
 \mlr = \sup \Big\{ d \big( x , \delconv \big): x \in \vcir \Big\},
$$ 
where $d:\R^2 \to [0,\infty)$ denotes the Euclidean distance.
We will denote by \hfff{mfl} {\it maximum facet length}
   $\mfl$ the length of the longest line segment of which the polygon $\delconv$ is comprised.
  \end{definition}
The present article is devoted to proving upper bounds on these two measures of local fluctuation. 

It is proved in  \cite{uzunalex} that $P_p \big( \mlr \geq n^{1/3} (\log n)^{-2/3} \big\vert \areacon \big) \to 1$ as $n \to \infty$, with $P_p$ subcritical bond percolation $p \in (0,1/2)$ on $\Z^2$. The power-law term here was expected to be sharp. 
For a broad class of subcritical models, local roughness was proved in \cite{alexcube} to be bounded above by $O \big(n^{1/3} (\log n)^{2/3} \big)$ 
in an $L^1$-sense, validating the sharpness of the power-law term for an averaged form of local roughness.  
An upper bound on $\mlr$ of $n^{2/3 + o(1)}$  appears in the same work. 

The new upper bounds on local deviation that we obtain here are sharp at both power-law precision and in the logarithmic refinement. They are reached by an application of surgical techniques also used in \cite{hammondthr},
which enable here uniform control on radial and longitudinal deviation. In this paper, we explain these surgeries heuristically before applying them formally, so that the paper both presents the derivations of some of our principal conclusions and acts as a guide to the surgical approach. 

Our central conclusion for radial local deviation is the following upper bound on maximum local roughness.
\begin{theorem}\label{thmmlrbd}
Let $P = \P_{\beta,q}$, with $\beta < \beta_c$ and $q \geq 1$. 
Then there exist $C> c >0$ and $t_0 \geq 1$ such that, for $t \geq t_0$, 
$t = O \big( n^{5/36} (\log n)^{-C}  \big)$,
$$
P \Big( \mlr \geq n^{1/3} \big( \log n \big)^{2/3} t \Big\vert \areacon \Big)
 \leq \exp \Big\{ - c t^{6/5} \log n  \Big\}.
$$
\end{theorem}
Maximum facet length was not an object explicitly considered by K. Alexander and H. Uzun, but it plays a central role in our approach. The principal conclusion  is:
\begin{theorem}\label{thmmflbd}
Let $P = \P_{\beta,q}$, with $\beta < \beta_c$ and $q \geq 1$. 
There exist $0 < c < C < \infty$
such that, for $t \geq C$, 
$t = o \big( n^{1/3} (\log n)^{-5/6}  \big)$,
$$
P \Big( \mfl \geq n^{2/3} \big( \log n \big)^{1/3} t \Big\vert \areacon \Big)
 \leq \exp \Big\{ - c t^{3/2} \log n  \Big\}.
$$
\end{theorem}
The article \cite{hammondtwo} proves lower bounds that complement Theorems \ref{thmmlrbd} and \ref{thmmflbd}.
Taken together, our conclusion is:
\begin{cor}\label{cormlrsum}
Let $P = \P_{\beta,q}$, with $\beta < \beta_c$ and $q \geq 1$. 
Then there exist constants $0 < c < C < \infty$ such that
$$
P \bigg(  c \leq \frac{\mlr}{n^{1/3} \big( \log n \big)^{2/3}} \leq C  \bigg\vert \areacon \bigg)
  \to 1,
\qquad \textrm{as $n \to \infty$,} 
$$
and 
$$
P \bigg(  c \leq \frac{\mfl}{n^{2/3} \big( \log n \big)^{1/3}} \leq C  \bigg\vert \areacon \bigg)
  \to 1, 
\qquad \textrm{as $n \to \infty$.} 
$$
\end{cor}
That is, the techniques of this paper and its counterpart \cite{hammondtwo}  are sufficient to derive the conjectured exponents for the power-laws in radial and longitudinal local deviation, and to identify and prove exponents for the logarithmic correction for these quantities. 

As we will shortly explain, an important ingredient in obtaining these results is an understanding that the conditioned circuit is highly regular, with little backtracking from its overall progress. 
The third article \cite{hammondthr}
in the series presents this result, on the renewal structure of the conditioned circuit.
We will state its main conclusion as Theorem \ref{thmmaxrg}. 
Equipped with this tool, it is straightforward to derive the results under conditioning on a fixed area:
\begin{theorem}\label{thmfixedarea}
Theorems \ref{thmmlrbd}, \ref{thmmflbd} and Corollary  \ref{cormlrsum} are valid (with verbatim statements) under the conditional measure $P \big( \cdot \big\vert \vert \intg \vert = n^2 \big)$.  
\end{theorem}
\begin{subsubsection}{Duality of random cluster models}
Let $\Z^{2,*}$ denote the dual lattice of $\Z^2$, obtained from it by translating by the vector 
$(1/2,1/2)$. We write $E\big( \Z^{2,*} \big)$ for its collection of nearest-neighbour edges.
Evidently, random cluster models may be defined on the dual lattice.

To any configuration $\omega \in \zoz$ is naturally associated a dual configuration 
$\omega^* \in \big\{ 0,1 \big\}^{E(\Z^{2,*})}$, in which a dual edge $e^*$ is declared open, i.e. 
$w^*(e^*)=1$, precisely when the edge $e \in E(\Z^2)$ which which it shares its midpoint is closed, i.e. $\omega(e) = 0$.

The process of dual edges given by a random cluster model is itself a random cluster model: 
with $\Lambda = \big\{ -n,\ldots, n\big\}^2 \subseteq \Z^2$ and 
$\Lambda^* = \big\{ -n-1/2,\ldots, n+1/2\big\}^2 \subseteq \Z^{2,*}$, and recalling that
$p = 1 - \exp \{ - 2 \beta \}$, let $p^*$ solve $\frac{p^*}{1 - p^*} = \frac{q(1-p)}{p}$, and let 
$\beta^*$ be such that $p^* = 1 - \exp \{ -  2 \beta^*  \}$. As explained for example in \cite{chayescs}, 
we have the duality relation
$$
 \P_{\Lambda,\beta,q}^f \big( \omega \big) = \P_{\Lambda^*,\beta^*,q}^w \big( \omega^* \big).
$$
In a high $n$ limit, we obtain 
$\P_{\beta,q}^f \big( \omega \big) = \P_{\beta^*,q}^w \big( A^* \big)$ for all cylinder events 
$A$ and $A^* = \big\{ w^*: w \in A \big\}$.
\end{subsubsection}
\begin{subsubsection}{The appearance of subcritical dual random cluster models via the Fortuin-Kasteleyn representation}
The Ising model (or, indeed, the $q$-state Potts model) 
and the random cluster model may be coupled as the vertex and edge marginals of a single stochastic process that was first introduced by Fortuin and Kasteleyn \cite{forkast} (and then exploited for the purpose of fast simulation by Swendsen and Wang \cite{swendsenwang}, and made more explicit by Edwards and Sokal \cite{edwardssokal}). 
In the coupling, any large Potts droplet has a boundary that is open in the dual, whereas,
in the interior of any open dual circuit, the vertex marginal has the law of the Potts model with free boundary condition, so that a droplet may form there irrespective of the dominant phase in the exterior.
The dual of the random cluster model itself being such a model,
it follows that 
%
conditioning an Ising (or Potts) model on having an excess of spins of a given type (that gather together to form a large droplet), and conditioning a subcritical random cluster model on having a circuit trapping high area, are highly related. The latter is thus a very natural model of phase separation.
\end{subsubsection} 
\begin{subsubsection}{Local roughness and the competition of curvature and fluctuation}
As we have briefly mentioned, among the numerous possible ways of measuring fluctuation of the droplet boundary, the local roughness 
definition is particularly interesting. Random deviation arises from the local Gaussian fluctuation of the subcritical open dual connections of which the circuit is comprised. Globally, however, the circuit is constrained by curvature. Measures of local roughness are interesting because their scaling is determined by the interplay of these fluctuation and curvature effects. It has been expected that the natural scale for such fluctuations is $n^{1/3}$ radially and $n^{2/3}$ longitudinally, and, indeed, this picture has already been validated in an $L^1$-sense: defining average local roughness ${\rm ALR}\big( \cir \big)$ to be the ratio of the area trapped between the circuit and its convex hull and the length of the convex hull, Alexander \cite{alexcube} has proved that ${\rm ALR}\big( \cir \big) = O \big( n^{1/3} (\log n)^{2/3}  \big)$. In this series of articles, we control local roughness in an $L^\infty$-sense up to a constant factor.
The methods of proof by which we derive power-law exponents and sharp logarthmic corrections for $\mlr$ and $\mfl$ illuminate the competition of fluctuation and curvature.  The methods that derive  the upper bounds on $\mlr$ and $\mfl$ in the present article
are a sequence of surgeries, and, for the complementary lower bounds in \cite{hammondtwo}, the analysis of a Markov chain at equilibrium designed to leave invariant the law $P\big( \cdot \big\vert \areacon \big) $ of the conditioned circuit.
We will shortly turn to a heuristic discussion of the method of proof for the upper bounds and of the form of the logarithmic corrections, in which the competition of local fluctuation and global curvature will be apparent. 

Such an interplay of effects is often found in stochastic models of repulsive particles in $1 + 1$ dimensions of space and time, such as the totally asymmetric exclusion process begun with occupation in the negatively indexed sites \cite{johansson}. In an illustrative example of the outcome of this interplay, the limiting law of the maximum of the difference of a Brownian motion and a parabola has been determined \cite{groeneboom} in terms of Airy functions. 
Problems naturally raised by the approach that we introduce and its results 
are to determine the mechanisms of the competition for phase separation boundaries, 
to identify laws for the rescaled circuit under the local-roughness scaling, 
and thereby to investigate the connection of the model to random systems in Kardar-Parisi-Zhang universality class \cite{kpz}. I would like to suggest two directions. 

Firstly, Corollary \ref{cormlrsum} strongly suggests that the following holds.
\begin{conjecture}\label{conjone}
Let $P = \P_{\beta,q}$, with $\beta < \beta_c$ and $q \geq 1$. Maximum facet length and maximum local roughness are asymptotically concentrated, as $n \to \infty$: there exist constants $d_1, d_2 \in (0,\infty)$ such that, for each $\epsilon > 0$, 
$$
P \Bigg( \bigg\vert \frac{\mlr}{n^{1/3}(\log n)^{2/3}} - d_1 \bigg\vert > \epsilon \Bigg\vert \areacon \Bigg) \to 1
$$ 
and 
$$
P \Bigg( \bigg\vert \frac{\mfl}{n^{2/3}(\log n)^{1/3}} - d_2 \bigg\vert > \epsilon \Bigg\vert \areacon \Bigg) \to 1
$$ 
as $n \to \infty$.
\end{conjecture}
Accepting Conjecture \ref{conjone}, it is natural to seek to find the fluctuation scale of $\mlr$ about 
$d_1 n^{1/3}(\log n)^{2/3}$ (or about its mean), and to study the limiting distribution of the difference 
$\mlr - d_1 n^{1/3}(\log n)^{2/3}$ normalized by this scale.

Another question that is very natural to pose is to consider the conditioned circuit, centred in a reasonable way at the origin (see Section \ref{seccircen}), and to examine its behaviour in the region where it cuts the positive $y$-axis. Consider a rescaled coordinate system whose origin is the point of contact of the circuit with the positive $y$-axis and in which the $x$-direction is scaled by a factor of $n^{2/3}$, and the $y$-direction by $n^{1/3}$. Then it would be interesting to find explicitly the law on functions $\R \to \R:t \to y(t)$ such that the rescaled process converges in distribution to the range of the process $t \to \big( t,y(t) \big)$. A priori, backtracking in the circuit may be so extensive to make the existence of such a real-valued function description invalid. However, Theorem \ref{thmmaxrg} of \cite{hammondthr}, and
\cite{hrynivioffe} for a related self-avoiding polygon model, forces circuit regularity down to a logarithmic scale, so that the question is well-posed.  
This question addresses typical behaviour, possibly making it more tractable than fine questions about $L^\infty$-defined objects such as maximum local roughness (although note that the law of maximum roughness may sometimes be expressed in terms of KPZ functionals \cite{majcomt}). 
The analogue of the question for a Brownian bridge $B:[-T,T] \to [0,\infty)$ conditioned to remain above the semi-circle of radius $T$ centred at $(0,0)$ has been found to satisfy a stochastic differential equation with a drift term expressed in terms of the Airy function \cite{ferrarispohn}. The analogue for the discrete polynuclear growth process is closely related to the point-to-point last passage time, which shares these radial and longitudinal scalings, and for which the rescaled last passage time process converges \cite{johdpg} to a translation by a parabola of the Airy process introduced by Pr\"ahofer and Spohn \cite{praspohn}.
 
The model that we discuss is static. It would be interesting to enquire as to whether it may be interpreted as fixed-time marginal of a growth process, with time equal to area. 
In this regard, we mention the true self-repelling motion, constructed rigorously in \cite{tothwerner}, 
a beautiful continuum model of growth and repulsion, in which trapped area plays such a role. 
\end{subsubsection}
\begin{subsubsection}{Further related models and possible extensions}
This work evolved from \cite{hammondperes}, in which the fluctuation from a circular trajectory of a planar Brownian loop conditioned to trap a high area was investigated. The Brownian problem was proposed by Senya Shlosman in part as a model problem, and, as such, it is a little ironic that the present work bounds local fluctuation for a range of more physical models much more sharply than did \cite{hammondperes} in the Brownian case. Some obvious analogues of the upcoming regeneration structure Theorem \ref{thmmaxrg} in the Brownian problem fail, although Vincent Beffara has pointed out that the existence of ``pivoting'' points on the Brownian path \cite{beffara} may provide a useful candidate for an analogous structure. In any case, the discrete analogue of the problem, in which a simple random walk in $\Z^2$ is run for time $n$ and conditioned to trap area $n^2$, could very possibly be tackled with the techniques of the present series of papers. It might be convenient to use a grand canonical variant of the problem, in which the lifetime of the walk is randomized geometrically, to ensure a convenient conditional independence of the walk trajectory under resampling. 

Regarding continuum problems, it would be of much interest to show that the local deviation behaviour proved in this paper is shared by the circuit arising in the subcritical near-critical scaling limit of percolation  on the triangular lattice 
(that has recently been constructed in \cite{gps2b}) under conditioning on the presence of a circuit enclosing the origin and trapping high area. 
\end{subsubsection}
\begin{subsection}{An overview of the surgical procedure}\label{secoverview}
In this section, we wish to communicate some of the main ideas that will be used to establish Theorems \ref{thmmlrbd} and \ref{thmmflbd}. 
We will take the configuration measure $P$ to be an independent subcritical bond percolation model on $\Z^2$, since this case includes the essential features. In the first instance, we will outline the argument that yields the upper bound of $n^{1/3 + o(1)}$ on maximum local roughness for the conditioned measure, and then explain how the argument may be altered to arrive at the sharp logarithmic correction.  

Under the conditioning $\areacon$, the outermost open circuit $\cir$ traces a path that successively visits those vertices contained in the boundary $\delconv$ of its convex hull. The path in $\cir$ that connects the two endpoints $\bo{x}$ and $\bo{y}$ of one of the line segments of which $\delconv$ is comprised is an open subcritical path.
In a subcritical model, conditional on two distant points being connected, the common cluster has 
the macroscopic profile of the interpolating line segment, and a Gaussian fluctuation from this line.
(In Section \ref{secoz}, we will state precise results from \cite{civ} valid for subcritical random cluster models.)

In attempting to rule out the possibility that $\mlr$, under the conditioning $\areacon$, exceeds an order of $n^{1/3}$, it is thus natural to proceed as follows. Under this conditioning, we may seek to rule out the event that $\mfl$ much exceeds $n^{2/3}$. Thus done, we would try to establish that, indeed, the interpolating open connections have an orthogonal fluctuation that has an order given by the square-root of the length of the line segment in question.

In trying to carry out the first of these two steps, set
$$
a(t) = P \Big(    \mfl \geq n^{2/3} t \Big\vert  \big\vert \intg \big\vert \geq n^2  \Big).
$$
We wish to show that $a(t)$ is small, for some $t = t(n) = n^{o(1)}$. 
Such an assertion will be presented formally in Proposition \ref{propmscbexc}. 

Under the event $\lmscb : = \big\{ \areacon \big\} \cap \big\{ \mfl \geq n^{2/3} t \big\}$,
let $\bo{x}$ and $\bo{y}$ denote the endpoints of the longest line segment in $\delconv$; so that, necessarily, $\vert\vert \bo{x} - \bo{y} \vert\vert \geq n^{2/3} t$.
Except with a negligibly small probability, all of $\cir$ lies in the ball $B_{\cctwo n}$, for some large but fixed contant $\cctwo > 0$. This means that there exist two deterministic vertices $\xo,\yo \in B_{\cctwo n}$ such that the conditional probability given $\lmscb$ that $\bo{x}$ is equal to $\xo$ and that $\bo{y}$ is equal to $\yo$ is at least $1/{\rm vol}\big(B_{\cctwo n}\big)^2 \geq \frac{1}{2\pi^2 \cctwo^4} n^{-4}$. (This polynomial decay rate should be considered to be fairly slow for the purposes of this argument.)

We will perform a surgical procedure on the configuration. The procedure will later be formally defined as the {\it sector storage-replacement operation}. Let $\axyo$ denote the angular sector rooted at the origin whose boundary line segments contain $\xo$ and $\yo$.

The procedure will be denoted by $\sigma_{\xo,\yo}$ and will act on a configuration. 
In the first step, the contents of the sector $\axyo$ are lifted out of the plane and recorded separately. The configuration in $\R^2 \setminus \axyo$ is left unaltered.  In the second step, an independent configuration is sampled in the presently empty sector $\axyo$ in the plane, to yield a new configuration in the plane. The outcome of the procedure, then, has two elements:
the original contents of the sector $\axyo$, and a new full-plane configuration where these contents have been independently sampled.

We are taking $P$ to be an independent subcritical bond percolation. As such, it is trivial that the procedure maps the measure $P$ to an output whose two elements are independent bond percolations in the sector and in the plane.

We will consider the action of $\sigma_{\xo,\yo}$ on an input having the law $P$. We will call the input {\it satisfactory} if
$\lmscb  \cap \big\{ \bo{x} = \xo \big\} \cap \big\{ \bo{y} = \yo \big\}$ occurs. (This event, recall, has probability at least $P \big( \areacon \big) a(t) \frac{1}{2\pi^2 \cctwo^4} n^{-4}$.)
The removed contents of the sector $\axyo$ appear necessarily to contain a (subcritical) open path from $\xo$ to $\yo$ under the event in question. (We will return to the imprecision in this statement later in the explanation.) Suppose further that the new configuration in $\axyo$ that is used by $\sigma_{\xo,\yo}$ to define the output full-plane configuration happens to enjoy the following properties (in which case, we will say that the action of $\sigma_{\xo,\yo}$ is {\it successful}):
\begin{itemize}
\item $\xo$ and $\yo$ are connected in $\axyo$ by an open path $\gamma$,
\item and the path $\gamma$ has an orthogonal fluctuation that deviates away from the origin more than it does towards it.
\end{itemize} 
The first event will be denoted by  $\xo \build\leftrightarrow_{}^{\axyo} \yo$.
Obviously, the latter has not been precisely stated. It will be made precise in Lemma \ref{lemgac} with the notion of {\it good area capture}. We roughly mean that $\gamma$ lies in the upper quartile for fluctuating ``outwards'' rather than ``inwards'', away from the origin, rather than towards it. The probability that the action of $\sigma_{\xo,\yo}$ is successful is
$c P \Big( \xo \build\leftrightarrow_{}^{\axyo} \yo \Big)$, where the factor of $c$ corresponds to the occurrence of outward fluctuation.

What can be said about the outcome of $\sigma_{\xo,\yo}$ if its input is satisfactory and its action is successful? Firstly, as we have remarked, we seemingly have that
\begin{eqnarray}
 & (1) &  \textrm{the removed contents of $\axyo$ have an open path from $\xo$ to $\yo$.}
\end{eqnarray}
Secondly, we appear to have that, for some small but fixed constant $c > 0$, 
\begin{eqnarray}
 & (2) &  \textrm{the new full-plane configuration contains an open circuit that traps an area} \nonumber \\
  & & \qquad \textrm{of at least $n^2 + c n t^{3/2}$.} \nonumber
\end{eqnarray}
\begin{figure}\label{figfirstop}
\begin{center}
\includegraphics[width=0.8\textwidth]{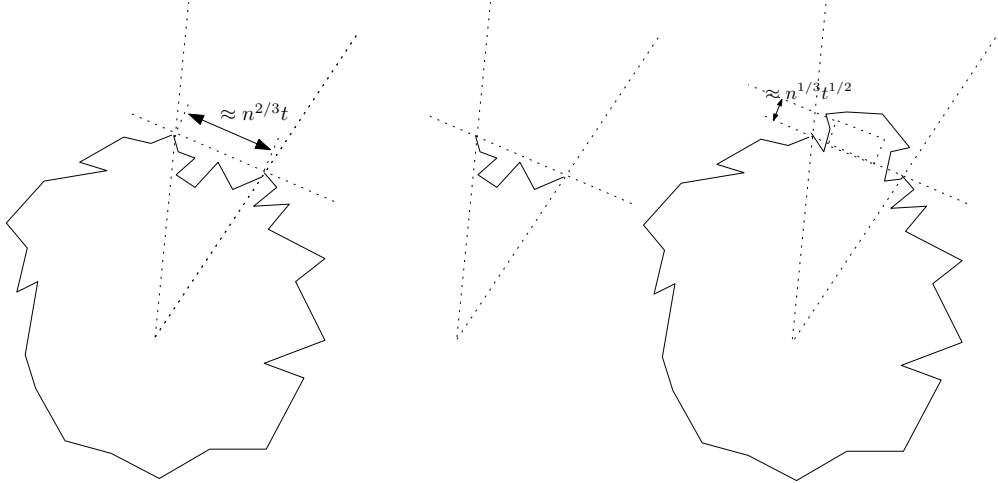} \\
\end{center}
\caption{A satisfactory input, and the two output components after successful action.}
\end{figure}
Figure $1$ 
illustrates the two output properties.
The reason for the second property is that we might expect 
the output full-plane configuration to contain an open circuit, that
coincides with the original circuit $\cir$ in $\R^2 \setminus \axyo$, and
with $\gamma$ in $\axyo$. The path $\gamma$ traps an area on the side of
the line segment $[\xo,\yo]$ opposite to the origin of the order of the length $\vert\vert \xo -
\yo \vert\vert$ of the segment multiplied by an orthogonal fluctuation of
the order of  $\sqrt{\vert\vert \xo - \yo \vert\vert}$. The original
circuit $\cir$ never crossed this line segment, however, because its
endpoints are members of $\delconv$. Thus, the change in the circuit gives
rise to an area-gain of at least the order of $\vert\vert \xo - \yo
\vert\vert^{3/2}$. Noting that $\vert\vert \xo - \yo \vert\vert =
\vert\vert \bo{x} - \bo{y} \vert\vert = \mfl \geq n^{2/3} t$ yields the
heuristic for why the second event should occur.

To summarize, in using the procedure on an input with law $P$, we make a ``probabilistic payment'' of
\begin{equation}\label{probinp}
P \Big( \areacon \Big) a(t) \frac{1}{2\pi^2 \cctwo^4 n^4}  c P \Big( \xo \build\leftrightarrow_{}^{\axyo} \yo \Big)  
\end{equation} 
to buy satisfactory input and successful action, and, in this way, to manufacture the circumstances under which the procedure yields an outcome for which properties 
$(1)$ and $(2)$ hold. However, we know that the law of the output is simply independent bond percolation on $\axyo$ and on $\Z^2$.  Thus, the outcome entails an event of probability at most
\begin{equation}\label{proboutp}
  P \Big( \xo \build\leftrightarrow_{}^{\axyo} \yo \Big) 
  P \Big( \big\vert \intg \big\vert \geq n^2 + c n t^{3/2} \Big).
\end{equation}
We learn that the quantity in (\ref{probinp}) is at most that in (\ref{proboutp}). Rearranging, we find that the conditional probability $a(t) = P \big( \mfl \geq n^{2/3} t \big\vert \areacon \big)$ that is the object of study satisfies
\begin{equation}\label{abd}
 a(t) \leq 2 \pi^2 \cctwo^4 c^{-1} n^4 
  P \Big( \big\vert \intg \big\vert \geq n^2 + c n t^{3/2} \Big\vert \areacon \Big).
\end{equation}
That is, a long line segment in the convex boundary can be surgically altered to force a certain excess in the area trapped over that mandated by the conditioning $\areacon$.
Note the massive cancellation of terms that led from the inequality $(\ref{probinp}) \leq (\ref{proboutp})$ to the conclusion (\ref{abd}).  With the purchase of satisfactory input and successful action, we receive output configurations satisfying properties (1) and (2). The large-scale cost of purchase and benefit of outcome have been designed to cancel out.
What is left over after that cancellation - in this case, the bound (\ref{abd}) - is a meaningful discussion of the fluctuation of the conditioned circuit, despite the probability of this conditioning being exponentially small in $n$.  

The conclusion (\ref{abd}) 
prompts us to seek to control the law of the area-excess \hfff{exc} $\exc : = \intg - n^2$ under the event $\areacon$, and we present the required result in Proposition \ref{propexc}. 
We will not at this moment indicate the method of proof. However, note that the exponential decay in $n \in \N$ of $P \big( \areacon \big)$ is consistent with a typical value of $\exc$ under $\areacon$ of order $n$, and with
$$
P \Big( \exc \geq nt \Big\vert \areacon\Big) \leq \exp \big\{ -ct \big\},
$$
at least for choices of $t = O(n)$ that are not extremely high. 

Allied with (\ref{abd}), we would obtain 
\begin{equation}\label{twothtwoth}
P \Big( \mfl \geq n^{2/3}t \Big\vert \areacon \Big) \leq \exp \big\{ - c t^{3/2} \big\},
\end{equation}
for $t \geq C \big( \log n \big)^{2/3}$. This is a conclusion of the form we sought.

\begin{figure}\label{figthreebox}
\begin{center}
\includegraphics[width=0.8\textwidth]{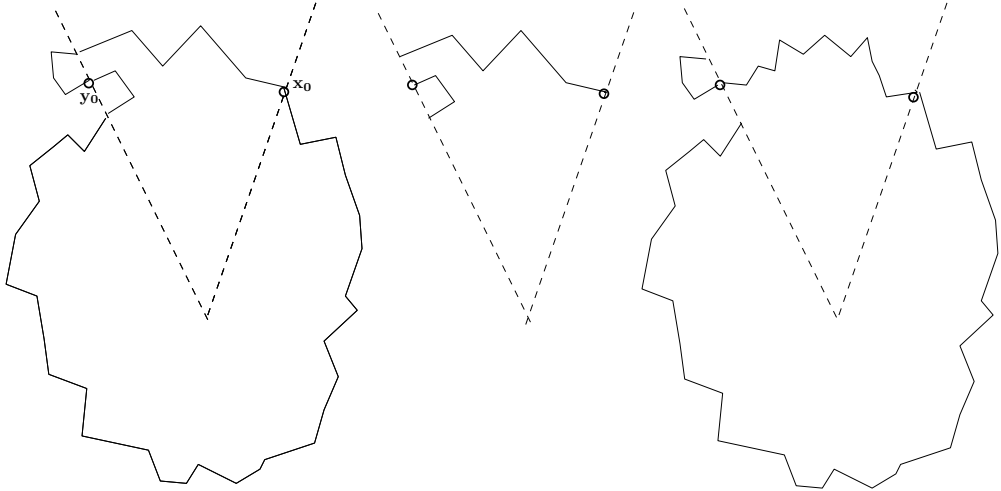} \\
\end{center}
\caption{An input for $\sigma_{\xo,\yo}$, and the sector and full-plane parts of the output.}
\end{figure}
However, the argument contains a flaw. Consider the action of
$\sigma_{\xo,\yo}$ in the case illustrated in Figure $2$.
The input is satisfactory and the action is successful, but neither property (1) nor (2) is satisfied by the output. The circuit in the input doubles back over the line segment from $\bo{0}$ through $\yo$, in such a way that the occurrence of $\xo \build\leftrightarrow_{}^{\axyo} \yo$ is insufficient to glue together the ends of the original circuit in $\axyo^c$. For similar reasons, the stored sector contains no open path from $\xo$ to $\yo$.

This problem motivates the following:
\begin{definition}
A site $\bo{v}$ in a circuit $\Gamma$ for which $\bo{0} \in {\rm INT}(\Gamma)$ is called a cutpoint of $\Gamma$ if the line segment $\big\{ t\bo{v}: t \geq 0 \big\}$ intersects $\Gamma$ only at $\bo{v}$.
\end{definition}
Indeed, as the reader may readily verify, we have that:
\begin{lemma}\label{lemcutp}
Let $\bo{v_1}, \bo{v_2} \in \Gamma$ be two cutpoints of $\Gamma$. 
Then
$$
{\rm INT}(\Gamma)    =  
{\rm INT} \Big( 
\big( \Gamma \cap \aarg{v_1}{v_2} \big) \cup \big[\bo{0},\bo{v_1}  \big] \cup 
 \big[\bo{0},\bo{v_2}  \big] \Big)  
\cup 
 {\rm INT} \Big(
\big( \Gamma \cap \aarg{v_1}{v_2}^c \big)
\cup \big[\bo{0},\bo{v_1}  \big] \cup 
 \big[\bo{0},\bo{v_2}  \big]
  \Big).
$$
\end{lemma}
In fact, we will actually make use of a notion that is slightly stronger than that of a cutpoint vertex, namely,
that of a regeneration site of $\cir$, that will be introduced in Section \ref{secrrg}.

The strategy for repairing the proof of the upper bound on 
$P \big( \mfl \geq n^{2/3}t \big\vert \intg \geq n^2 \big)$
requires that regeneration sites are, in fact, typically found in all parts of the circuit $\cir$ under the event $\areacon$. We will shortly state Theorem \ref{thmmaxrg} that is sufficient for this purpose. Its proof appears in \cite{hammondthr}. 
  We then let $\bo{x'}$ and $\bo{y'}$ denote the nearest regeneration sites to the endpoints $\bo{x}$ and $\bo{y}$, and apply the earlier argument, choosing instead $\xo$ and $\yo$ to be such that $\bo{x'}$ and $\bo{y'}$ have a probability of order $n^{-4}$ to be equal to $\xo$ and $\yo$. We then analyse the action of the procedure in the event that 
$\big\{ \bo{x'} = \xo \big\} \cap \big\{ \bo{y'} = \yo \big\}$ occurs. We are now resampling $\axyo$ in the case that $\xo$ and $\yo$ are regeneration sites, so that Lemma \ref{lemcutp} may be applied to bound the area enclosed by the (now unbroken) open circuit in the full-plane configuration resulting from a successful surgical procedure. Theorem \ref{thmmaxrg} tells us that $\bo{x'}$ is close to $\bo{x}$, and $\bo{y'}$ to $\bo{y}$, so that the earlier argument on area is valid with an insignificant additional error term.

Logically, then, our first step is to establish Theorem \ref{thmmaxrg}. In this sense, the paper \cite{hammondthr} is a precursor to this one. However, the result is technical, and relies on  
variants of the sector storage-replacement operation introduced in this article. The proofs of that paper have been structured similarly to those found here. To the reader who is primarily  interested in understanding the proof of Theorem \ref{thmmaxrg}, we suggest to begin by reading the proof of Proposition \ref{propmscbexc} in the present article (which is the formal argument corresponding to the preceding discussion). 
\end{subsection}
\begin{subsection}{An explanation of the form of the logarithmic power corrections}
\begin{subsubsection}{Maximum facet length}
With Theorem \ref{thmmaxrg} available, the argument that we have just sketched  leads to the conclusion (\ref{twothtwoth}), which implies that
$$
P \Big(  \mfl \leq C n^{2/3}  \big( \log n \big)^{2/3} \Big\vert  \areacon \Big) \to 1.
$$
As we turn to make the argument rigorous in the following sections, we obtain an improvement, namely,
that  $\mfl = O \big( n^{2/3} \big( \log n \big)^{1/3} \big)$. We here summarize the reason for this. The existing argument used a definition of satisfactory input that was prepared to make an order $n^{-4}$ probabilistic payment to ensure that the ``right'' parameters $\xo,\yo$ were used in the action of $\sigma_{\xo,\yo}$.  
Having made one polynomially-sized payment, we may as well make another: by changing the definition of successful action of the operation to require that the outward orthogonal
fluctuation of the  open  connection in the updated configuration in $\axyo$ lies in a polynomially small (in $n$) upper tail, we should be able to trap a little more area by the open circuit in the full-plane output. 

The open path $\gamma$ connecting $\xo$ to $\yo$ in $\axyo$ on whose existence in the updated configuation in $\axyo$ we insist as part of the definition of successful action of $\sigma_{\xo,\yo}$ has an asymptotically Gaussian fluctuation. It deviates to maximal distance of $\vert\vert \xo - \yo \vert\vert^{1/2} \big( \log \vert\vert \xo - \yo \vert\vert  \big)^{\alpha}$
with probability of the order of the Gaussian tail,
$$
\exp \Bigg\{ - \frac{\Big( \vert\vert \xo - \yo  \vert\vert^{1/2} (\log n)^{\alpha} \Big)^2}{\vert\vert \xo - \yo  \vert\vert} \Bigg\} = \exp \Big\{ - \big( \log n \big)^{2\alpha} \Big\},
$$ 
since $\log \vert\vert \xo - \yo  \vert\vert \approx \log n$. For $\alpha \leq 1/2$, this is an at most polynomial decay. 

However, if the orthogonal fluctuation behaves as $\vert\vert \xo - \yo \vert\vert^{1/2} \big( \log n \big)^{1/2}$, then the open circuit of property $(2)$ satisfied by the output would trap an area of at least $n^2 + c n t^{3/2} \big( \log n \big)^{1/2}$,
so that an additional factor of $\big( \log n \big)^{1/2}$ would appear in the area-excess term of (\ref{abd}), and (\ref{twothtwoth}) becomes
$$
P \Big(  \mfl \geq n^{2/3}t   \Big\vert   \acon \Big) \leq \exp \Big\{ - c t^{3/2} \big( \log n \big)^{1/2} \Big\},
$$
if $t = t(n)$ is chosen so that the right-hand-side decays at least at a fast polynomial rate (so that, indeed, the polynomial in $n$ factors are insignificant). It is the inequality $t \geq C \big( \log n \big)^{1/3}$
that ensures this. This explains how we obtain the upper bound of $C n^{2/3} \big( \log n \big)^{1/3}$ for $\mfl$.
\end{subsubsection}
\begin{subsubsection}{Maximum local roughness}\label{secmlrlog}
Similar reasoning yields the conclusion that $\mlr$ is at most a constant multiple of $n^{1/3} \big( \log n \big)^{2/3}$, typically. Indeed, suppose that,
under $\areacon$,
maximum local roughness were, with a significant probability, rather larger
than this, larger than 
$n^{1/3} \big( \log n \big)^{2/3 + \epsilon}$ for some $\epsilon > 0$, say. By now, we
know that we may assume that $\mfl \leq C n^{2/3} \big( \log n
\big)^{1/3}$. 
We may apply the sector storage-replacement operation $\sigma_{\xo,\yo}$
for a choice of $(\xo,\yo)$ that has a probability of order $n^{-4}$  
of coinciding with the endpoint-pair of the line segment in $\delconv$
corresponding to the point in $\vcir$ that attains $\mlr$. Buying an input
configuration satisfying $\areacon$ and an open path from $\xo$ to $\yo$ in
$\axyo$ in the updated sector
configuration produced by $\sigma_{\xo,\yo}$, and further insisting on an
upper-quartile outward orthogonal fluctuation for this connection, we
achieve an open circuit in the full-plane output with an area capture of at
least the area enclosed by $\cir$ in the input, that is, at least
$n^2$. This mirrors the cost of the event $\areacon$ for the input. 
However, we also obtain
a diagnozable abnormality in the sector component of the output. This sector
contains an open path from  $\xo$ to $\yo$ in $\axyo$, 
which is no more unusual than the
price we chose to pay; but, in addition, this path has an orthogonal
fluctuation of at least $\mlr \geq 
n^{1/3} \big( \log n \big)^{2/3 + \epsilon}$. The path interpolates two
points at distance at most $\mfl \leq  C n^{2/3} \big( \log n
\big)^{1/3}$. An open path having an order of orthogonal fluctuation given by a
Gaussian, such a fluctuation has a probability of the order of at most
$$
\exp \Bigg\{ - c \frac{ \Big( n^{1/3}  \big( \log n \big)^{2/3 + \epsilon}
  \Big)^2}{C n^{2/3} \big( \log n \big)^{1/3}} \Bigg\}  = \exp \Big\{ - c \big(
\log n \big)^{1 + 2\epsilon} \Big\}.
$$ 
An oddity, then, of superpolynomial rarity. It cannot, in fact, occur,
because only a polynomial price has been paid to achieve it. So the scenario that   
$\mlr \geq C n^{1/3} \big( \log n \big)^{2/3 + \epsilon}$ under $\areacon$ 
has been excluded;
or, more accurately, has been shown to arise with a super-polynomially
decaying conditional probability. 
\end{subsubsection}
\end{subsection}
\begin{subsection}{The structure of the paper}
We describe how the rest of the paper is structured.
The formal counterpart of (\ref{abd}), Proposition \ref{propmscbexc}, is presented and proved in Section \ref{secunu}.
The tail of the law of the area-excess is then bounded above in Proposition \ref{propexc}, appearing in Section \ref{secarea}, the result being proved by a different surgical procedure that is introduced there. In combination, these two propositions yield Theorem \ref{thmmflbd}. The upper bound on maximum local roughness (Theorem \ref{thmmlrbd}) uses Theorem \ref{thmmflbd}. The proof of Theorem \ref{thmmlrbd}, which is the rigorous version of the argument presented in Subsection \ref{secmlrlog}, is supplied in Section \ref{secubdmlr}.

Various preliminaries are required before we give these proofs. 
As we have seen, we require that the conditioned circuit has a high degree of regularity, with regeneration points pervading the circuit. 
In the next section, Section \ref{secnott}, after fixing some notation, we explain our convention for fixing the centre of the circuit 
 and then
state in the radial definition of regeneration site and Theorem \ref{thmmaxrg} on the profusion of such sites in the conditioned circuit. (The proof of Theorem \ref{thmmaxrg} appears in \cite{hammondthr}.)
The surgeries that we use rely on information about the fluctuation and renewal structure of point-to-point conditioned connections in subcritical random cluster models, and we conclude Section \ref{secnott} by presenting the results, developed by \cite{civ}, that we require from this theory. \\
\noindent{\bf Acknowledgments.} I would like to thank Kenneth Alexander, Rapha\"{e}l Cerf, Dima Ioffe, James Martin, Yuval Peres and Senya Shlosman for helpful discussions.
\end{subsection}
\end{section}
\begin{section}{Notation and tools}\label{secnott}
\begin{subsection}{Notation}\label{secnot}
\begin{definition}\label{defpathedge}
Elements of $\R^2$ will be denoted by boldface symbols. 
By a discrete path, we mean a list of elements of $\Z^2$, each being a nearest-neighbour of the preceding one, and without repetitions. 
In referring to a path, we mean a subset of $\R^2$ given by the union of the topologically closed edges formed from the set of consecutive pairs of vertices of some discrete path. 
(As such, a path is defined to be self-avoiding, including at its vertices.)
In a similar vein, any subset of $\R^2$ that is introduced as a connected set is understood to be a union of closed intervals $\big[ \bo{u},\bo{v} \big]$ corresponding to nearest-neighbour edges $(\bo{u},\bo{v})$.
For such a set $A$, we write $V(A) = A \cap \Z^2$ and $E(A)$ for the set of edges of which $A$
is comprised. 

For a general subset $A \subseteq \R^2$, we write $E(A)$ for the set of nearest-neighbour edges 
$(\bo{u},\bo{v}) \in E(\Z^2)$ such that $\big[ \bo{u},\bo{v} \big] \subseteq A$. (This is of course consistent with the preceding definition.) We write 
$E^*(A)$ for the set of nearest-neighbour edges $(\bo{u},\bo{v}) \in E(\Z^2)$ such that $\big[ \bo{u},\bo{v} \big] \cap A \not= \emptyset$. 
\end{definition}
\begin{definition}\label{deftriandothers}
For  $\bo{x},\bo{y} \in \Z^2$, $\bo{y} \not= \bo{x}$,
we write $\ell_{\bo{x},\bo{y}}$ for the planar line containing $\bo{x}$ and $\bo{y}$, and $\ell^+_{\bo{x},\bo{y}}$ for the semi-infinite line segment that contains $\bo{y}$ and has endpoint $\bo{x}$. We write $\big[ \bo{x},\bo{y} \big]$ for the line segment whose endpoints are $\bo{x}$ and $\bo{y}$. 
We write \hfff{txy} $T_{\bo{0},\bo{x},\bo{y}}$ for the
closed triangle with vertices $\bo{0}$, $\bo{x}$ and $\bo{y}$.
For $\bo{x},\bo{y} \in \R^2 \setminus \{ \bo{0} \}$, 
we write $\ang\big(\bo{x},\bo{y} \big) \in [0,\pi]$ for the angle between these two vectors.
Borrowing complex notation, we write $\argu\big(\bo{x}\big)$ for the argument of $\bo{x}$.
In many derivations, the cones, line segments and points in question all lie in a cone, rooted at the origin, whose aperture has angle strictly less than $2\pi$. As such, it is understood that $\argu$
denotes a continuous branch of the argument that is defined throughout the region under consideration. 
\end{definition}
Sometimes we wish to specify a cone by a pair of boundary points, and sometimes by the argument-values of its boundary lines:
\begin{definition}
For $\bo{x},\bo{y} \in \Z^2$, $\argu(\bo{x}) < \argu(\bo{y})$, write  \hfff{axy}
$$
 A_{\bo{x},\bo{y}} = \Big\{ \bo{z} \in \R^2 \setminus \{ \bo{0} \} : \argu\big( \bo{x} \big) \leq  \argu\big( \bo{z} \big) \leq  \argu\big( \bo{y} \big)  \Big\} \cup \big\{ \bo{0} \big\}.
$$ 
To specify a cone by the argument-values of its boundary lines, take $\bo{v} \in \Z^2$ and $c \in [0,\pi)$, and let \hfff{wvc}
\begin{equation}\label{wnot}
W_{\bo{v},c} = \Big\{ \bo{z} \in \R^2 \setminus \{ \bo{0} \}: \argu(\bo{v}) - c \leq  \argu(\bo{z}) \leq  \argu(\bo{v}) + c \Big\}
\cup \big\{ \bo{0} \big\}
\end{equation}
denote the cone of points whose angular displacement from $\bo{v}$ is at most $c$. 
Extending
this notation, for any $\bo{x} \in \Z^2$ and $c \in [0,\pi)$, we write 
$W_{\bo{v},c}\big( \bo{x} \big) = \bo{x} + W_{\bo{v},c}$.
We also write, for $\bo{x} \in \R^2$ and $c \in (0,2\pi)$, \hfff{wvcalt}
$$
W_{\bo{x},c}^+  = \Big\{ \bo{z} \in \R^2 \setminus \{ \bo{0} \}: \argu(\bo{x})  \leq  \argu(\bo{z}) \leq  \argu(\bo{x}) + c \Big\}
\cup \big\{ \bo{0} \big\}
$$
and 
$$
W_{\bo{x},c}^-  = \Big\{ \bo{z} \in \R^2 \setminus \{ \bo{0} \}: \argu(\bo{x}) - c \leq  \argu(\bo{z}) \leq  \argu(\bo{x})  
\Big\}
\cup \big\{ \bo{0} \big\}.
$$
\end{definition}
\begin{definition}
For $\bo{v} \in \R^2$, let $\bo{v}^{\perp} \in S^1 $ denote the direction vector 
obtained by a counterclockwise rotation of $\pi/2$ radians from the direction of $\bo{v}$. 
\end{definition}
\begin{definition}\label{defmarg}
For $P$ a probability measure on $\zoz$ and for $\omega' \in \{ 0,1 \}^A$ for some $A \subseteq E(\Z^2)$, we write $P_{\omega'}$ for the conditional law of $P$ given $\omega\big\vert_A = \omega'$. We will also write $P \big( \cdot \big\vert \omega' \big)$ for $P_{\omega'}$.
\end{definition}
\begin{definition}
Given a subset $A \subseteq \R^2$, two elements $\bo{x},\bo{y} \in \Z^2 \cap A$, 
and a configuration  $\omega \in \{ 0,1 \}^{E(A)}$,  we write 
$\bo{x} \build\leftrightarrow_{}^A \bo{y}$
 for the event that there exists an $\omega$-open path from $\bo{x}$ to $\bo{y}$ all of
whose edges lie in $E(A)$. By the ($\omega$-)open component of $\bo{x}$ in $A$, 
we mean the connected subset of
$A$ whose members lie in an edge belonging to an ($\omega$-)open path in $E(A)$ 
that begins at $\bo{x}$.
\end{definition}
Throughout, the notation $\dist \cdot \dist$ and $d \big( \cdot , \cdot \big)$ refers to the Euclidean metric on $\R^2$. For $\gamma \subseteq \R^2$, we set ${\rm diam}(\gamma) = \sup \big\{ d\big(\bo{x},\bo{y}\big): \bo{x},\bo{y} \in \gamma \big\}$. For $K > 0$, we set $B_K = \big\{ \bo{x} \in \R^2: \dist \bo{x} \dist \leq K \big\}$.
\begin{subsubsection}{Convention regarding constants}
Some constants are fixed in all arguments: in particular, $c_0$ and $\qzero$ in the upcoming Definition \ref{defrg} of 
$\cir$-regeneration site. A few constants will be fixed in upcoming lemmas, in which case, they carry  letter subscripts that are acronyms, evident from the defining context. 
Throughout, constants with an upper case $C$ indicate large constants, and, with a lower case, small constants.
 The constants $C$ and $c$ may change from line to line and are used to absorb and simplify expressions involving other constants. 
\end{subsubsection}
\end{subsection}
\begin{subsection}{Decorrelation of well-separated regions: ratio-weak-mixing}\label{secrwm}
The following spatial decorrelation property is well-suited to analysing the conditioned circuit.
\begin{definition}\label{defrwm}
A probability measure $P$ on $\{ 0,1  \}^{E(\Z^2)}$ is said to satisfy the
ratio-weak-mixing property if, for some $\crwm,\lambda > 0$, and for all sets 
$\mc{D}, \mc{F} \subseteq E \big( \Z^2 \big)$,
$$
\sup \Big\{ \Big\vert \frac{P \big( D \cap F \big)}{P\big( D \big) P \big(
  F \big)} - 1 \Big\vert: D \in \sigma_{\mc{D}}, F \in \sigma_{\mc{F}}, 
  P \big( D \big) P \big( F \big) > 0   \Big\} 
\leq \crwm \sum_{x \in V(\mc{D}), y \in V(\mc{F})} e^{- \lambda \vert x - y \vert},
$$
whenever the right-hand side of this expression is less than one. Here, for $A \subseteq E(\Z^2)$, $\sigma_A$ denotes the set of configuration events measurable with respect to the variables $\big\{ \omega(e): e \in A \big\}$.
\end{definition}
The ratio-weak-mixing property is satisfied by any  $P = \P_{\beta,q}$, with $\beta < \beta_c$, that is, by any random cluster model with exponential decay of connectivity (Theorem 3.4, \cite{alexon}).
Any such measure trivially satisfies the following condition:
\begin{definition}\label{defbden}
A probability measure $P$ on $\{ 0,1  \}^{E(\Z^2)}$ satisfies 
the bounded energy property if there exists a constant $\cposen > 0$ such that,  for any $\omega' \in \zoz$ and an edge $e \in E(\Z^2)$, the conditional probability that $\omega(e) = 1$ given the marginal $\omega' \big\vert_{E(\Z^2) \setminus \{ e \}}$ is bounded between $\cposen$ and $1 - \cposen$.
\end{definition}
The following tool will establish near-independence of separated regions:
\begin{lemma}\label{lemkapab}
Given $A,B \subseteq E(\Z^2)$ and $m > 0$, let 
$$
\kappa_m \big( A,B \big) = \sum_{\bo{x} \in V(A),\bo{y} \in V(B), \vert\vert
  \bo{x} - \bo{y} \vert\vert \geq m} \exp \big\{ - \lambda \vert\vert
\bo{x} - \bo{y} \vert\vert \big\},
$$
where 
we write $V(A)$ for the set of vertices incident to one of the edges comprising $A$ (and similarly, of course, for $V(B)$). The constant $\lambda > 0$ is fixed; subsequent constants may depend on it implicitly.
Set $\phi_m \big( A,B \big) = \big\vert \big\{ \bo{a} \in A, \bo{b} \in B:
\vert\vert \bo{a} - \bo{b}
\vert\vert \leq m \big\} \big\vert$.
We say that $A$ and $B$ are $(m,C_0)$-well separated if $A \cap B = \emptyset$,
$\kappa_m \big( A,B \big) \leq 1/(2\crwm)$ and $\phi_m \big( A,B \big) \leq
C_0$. Here, $\crwm$ denotes the constant appearing in the definition (\ref{defrwm}) of ratio-weak-mixing. 
Let $P$ be a probability measure on $\{ 0,1  \}^{E(\Z^2)}$  satisfying the ratio-weak-mixing and bounded energy properties.
For $\clemkam \in \N$, $\clemkac \in \N$, there exists $\conka =
C\big(\clemkam,\clemkac\big)$ such that, if $A,B \subseteq E(\Z^2)$ are $(\clemkam,\clemkac)$-well
separated, then, for any $G \in \sigma_B$,
\begin{equation}\label{omubd}
\sup_{\omega \in \{ 0,1 \}^A} P_\omega \big( G \big) \leq \conka P(G)
\end{equation}
and
\begin{equation}\label{omlbd}
\inf_{\omega \in \{ 0,1 \}^A} P_\omega \big( G \big) \geq \conka^{-1} P(G).
\end{equation}
\end{lemma}
\noindent{\bf Proof.} For $n \in \N$, set $B_{\clemkam} = \big\{ b \in B: d \big( b, A \big) > \clemkam \big\}$.
Write 
$$
G_{\clemkam} = \bigg\{  \omega \in \{ 0,1 \}^{B_{\clemkam}}:  P_\omega(G) > 0 \bigg\}.
$$
Clearly, the event that the configuration in $B_{\clemkam}$
 satisfies $G_{\clemkam}$ is at least as probable under the measure $P$ as that of the configuration in $B$ satisfying $G$, or, indeed, under the measure $P_{\omega'}$ for arbitrary $\omega' \in \{ 0,1 \}^A$. For $\omega' \in A$, then, $P_{\omega'}(G) \leq P_{\omega'}(G_m)$. Now, if $A$ is finite,
 $$
    \frac{P_{\omega'}(G_{\clemkam})}{P(G_{\clemkam})} =   \frac{P \Big( G_{\clemkam} \cap \Big\{ \omega\big\vert_A = \omega' \Big\} \Big)}{P \big(G_{\clemkam} \big) P \big( \omega\big\vert_A = \omega' \big) }.
 $$ 
(The case of infinite $A$ is treated by taking a limit.)
 By the definition (\ref{defrwm}) of ratio weak mixing, the right-hand-side differs from one by at most $C \sum_{\bo{x} \in V(B_{\clemkam}), \bo{y} \in V(A)} \exp \big\{ - \lambda \vert\vert \bo{x} - \bo{y} \vert\vert \big\}$. However, this sum is at most $\kappa_{\clemkam}(A,B)$. Since $A$ and $B$ are $(\clemkam,\clemkac)$-well separated, we learn that
 \begin{equation}\label{pomgm}
  P_{\omega'}\big( G \big) \leq P_{\omega'} \big( G_{\clemkam} \big) \leq \frac{3}{2}  P \big( G_{\clemkam} \big).
 \end{equation}
 Note that, for any $\omega' \in G_{\clemkam}$, there exists $\tilde\omega \in \{ 0,1 \}^{B \setminus B_{\clemkam}}$ such that $(\omega',\tilde\omega) \in \{ 0,1 \}^B$ realizes $G$. However,
 $$
 \inf_{\tilde\omega \in \{ 0,1 \}^{B \setminus B_{\clemkam}}} P_{\omega'} \Big( \omega\big\vert_{B \setminus B_{\clemkam}} = \tilde{\omega} \Big) \geq \cposen^{\vert B \setminus B_{\clemkam} \vert},
 $$   
 by the bounded energy property satisfied by the measure $P$. Thus,
 \begin{equation}\label{pbbm}
  P \big( G \big) \geq P \big( G_{\clemkam} \big) \cposen^{\vert B \setminus B_{\clemkam} \vert}.
 \end{equation}
 Note that $\vert B \setminus B_{\clemkam} \vert \leq \phi_{\clemkam} \big( A , B \big)$. 
  Since $A$ and $B$ are $(\clemkam,\clemkac)$-well separated, we find that
   $P \big( G \big) \geq P \big( G_{\clemkam} \big) \cposen^{\clemkac}$. Returning to (\ref{pomgm}), we find that 
(\ref{omubd}) is valid with the choice $\conka = (3/2) \cposen^{-\clemkac}$. It follows similarly to (\ref{pbbm}) that, for any $\omega' \in \{ 0,1 \}^A$, $P_\omega \big( G \big) \geq P_\omega \big( G_m \big) \cposen^{\vert B \setminus B_{\clemkam} \vert}$. By the first inequality of (\ref{pomgm}) and as previously, we obtain (\ref{omlbd}) with the choice $\conka = \cposen^{- \clemkac}$. \qed
\end{subsection}
\begin{subsection}{The Wulff shape and circuit centering}\label{seccircen}
The macroscopic profile of the conditioned circuit is given by the boundary of the Wulff shape.
\begin{definition}
 We define the {\it inverse correlation length}: for $\bo{x} \in \R^2$,
$$
\xi(\bo{x}) = - \lim_{k \to \infty} k^{-1} \log P \big( \bo{0} \leftrightarrow \lfloor k\bo{x} \rfloor  \big),
$$
where $\lfloor \bo{y} \rfloor \in \Z^2$ is the component-wise integer part of $\bo{y} \in \R^2$. 
\end{definition}
\begin{definition}
The unit-area Wulff shape  \hfff{wulff} $\wulff$ is the compact set given by
$$
\wulff 
= \lambda \bigcap_{\bo{u} \in S^1} \Big\{ \bo{t} \in \R^2: \big( \bo{t},\bo{u} \big) \leq \xi\big(\bo{u}\big)  \Big\},
$$
with the dilation factor $\lambda > 0$ chosen to ensure that $\big\vert \wulff \big\vert = 1$.
\end{definition}
The following appears in Theorem B of \cite{civ}:
\begin{lemma}\label{lemozstr}
Let $P = \P_{\beta,q}$ with $\beta < \beta_c$ and $q \geq 1$. 
Then $\wulff$ has a locally analytic, strictly convex boundary.
\end{lemma}
Global deviations of the conditioned circuit from the Wulff shape may be measured in the following way.
\begin{definition}
Let $\Gamma \subseteq \R^2$ denote a circuit.  Define its
global distortion ${\rm GD} \big( \Gamma \big)$ (from an factor $n$ dilate of the Wulff shape boundary) by means of \hfff{globdis}
\begin{equation}\label{eqdefgd}
{\rm GD} \big( \Gamma \big) = \inf_{\bo{z} \in \Z^2} d_H \Big(  n \partial \wulff + \bo{z}, \cir  \Big),
\end{equation}
where $d_H$ denotes the Hausdorff distance on sets in $\R^2$. 
\end{definition}
(In a general context, this would be a peculiar definition. However, we will work with this quantity only in the case of circuits that are conditioned to trap an area of at least, or exactly, $n^2$.)

In undertaking surgery, it is convenient to work with circuits that are centred at the origin in a sense we now define.
\begin{definition}\label{defncentr}
Let $\Gamma \subseteq \R^2$ denote a circuit. 
The lattice point $\bo{z}$ attaining the minimum in (\ref{eqdefgd}) will be called the centre  \hfff{centre} $\centre(\Gamma)$ of $\Gamma$. In the case that the minimum is not uniquely attained, we take $\bo{z}$ to be the lexicographically minimal among those points in $\Z^2$ that attain the minimum. 
\end{definition}
\begin{definition}\label{defarea}
Let $A \in \N$. We write \hfff{acon} $\area{A}$ for the event $\big\{ \big\vert {\rm INT}(\Gamma_0) \big\vert \geq A \big\} \cap \big\{ \centre(\cir) = \bo{0} \big\}$.
\end{definition}
\end{subsection}
\begin{subsection}{Radial regeneration structure}\label{secrrg}
We now formulate the assertion on circuit regularity that is a vital element in our approach.
\begin{definition}\label{defrg}
For $\bo{u} \in S^1$,
let $w_{\bo{u}}$ denote the counterclockwise-oriented unit tangent vector
to $\partial \wulff$ at $\partial \wulff \cap \big\{ t\bo{u}: t \geq 0  \big\}$.
Let $\qzero > 0$ satisfy
\begin{equation}\label{supang}
\sup_{\bo{z} \in S^1}  \ang \big( w_{\bo{z}} , \bo{z}^{\perp} \big) \leq \pi/2 - 4\qzero,
\end{equation}
a choice made possible by the compactness and convexity of $\wulff$.

\hfff{forback} 
The forward cone $C^F_{\pi/2 - \qzero} \big( \bo{v} \big)$ denotes the set of vectors $\bo{w} \in \R^2$ for which $\ang\big( \bo{w} - \bo{v}, \bo{v}^{\perp} \big) \leq \pi/2 - \qzero$. The backward cone $C^B_{\pi/2 - \qzero} \big( \bo{v} \big)$
 denotes the set of vectors $\bo{w} \in \R^2$ for which $\ang\big( \bo{w} - \bo{v}, - \bo{v}^{\perp} \big) \leq \pi/2 - \qzero$.

Let $c_0 \in (0,q_0/2)$ be chosen in such a way that,  whenever $\bo{x},\bo{y} \in \partial
\wulff$,
$\argu(\bo{x}) < \argu(\bo{y})$ and 
$\ang \big( \bo{x} , \bo{y} \big) \leq 2 c_0$, then
\begin{equation}\label{czercond}
 \ang \big( \bo{x} - \bo{y}, - \bo{y}^{\perp}  \big) \leq
\pi/2 - 3\qzero.  
\end{equation}
This is possible by (\ref{supang}) and Lemma \ref{lemozstr}.

We henceforth fix a choice of constants $c_0$ and $\qzero$ satisfying these conditions, (as well as another one, that will be stated in Subsection \ref{secldgd}).  A site $\bo{v} \in \vcir$ in a circuit $\cir$ for which $\bo{0} \in \intg$ is called a $\cir$-regeneration site if
\begin{equation}\label{eqdefreg}
\cir \cap 
 W_{\bo{v},c_0} \subseteq C^F_{\pi/2 - \qzero} \big(  \bo{v} \big) \cup C^B_{\pi/2 - \qzero} \big( \bo{v} \big).
\end{equation}
See Figure $3$.
We write \hfff{reg} $\reg$ for the set of $\cir$-regeneration sites.
\end{definition} 
\begin{figure}\label{figregdef}
\begin{center}
\includegraphics[width=0.3\textwidth]{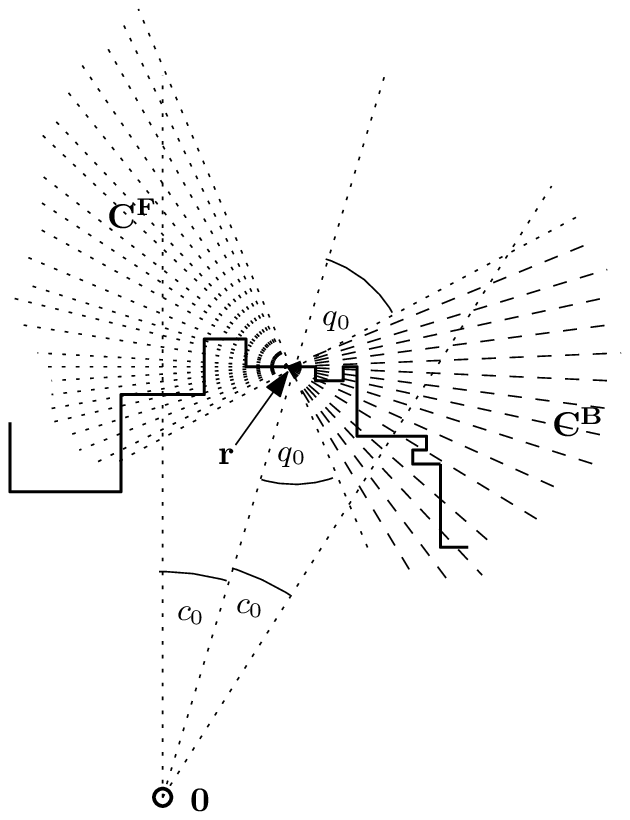} \\
\end{center}
\caption{A $\cir$-regeneration site $\bo{r}$ and the nearby circuit.}
\end{figure}
\begin{definition}\label{defmar}
Let $\Gamma$ denote a circuit for which $0 \in {\rm INT}\big(\Gamma \big)$.
We write $\margam \in [0,2\pi]$ for the angle of the largest angular sector rooted at the origin that contains no $\Gamma$-regeneration sites. That is,
\begin{equation}\label{eqdefmar}
\margam = \sup \Big\{ r \in [0,2\pi]: \exists \, \bo{a} \in S^1, 
 W_{\bo{a},r/2} \big( \bo{0} \big) \cap \reggam = \emptyset  \Big\}.
\end{equation}
\end{definition}
Our result on the regularity of the conditioned circuit is now stated.
\begin{theorem}\label{thmmaxrg}
Let $P = \P_{\beta,q}$ with $\beta < \beta_c$ and $q \geq 1$. 
There exist $c > 0$ and $C > 0$ such that
$$
P \Big(  \mar > u/n   \Big\vert \acon \Big)
 \leq  \exp \Big\{ - c u \Big\} 
$$
for $C \log n   \leq u \leq c n$.
\end{theorem}
We also record a simple lemma regarding control of $\cir$ near regeneration sites:
\begin{lemma}\label{lemdistang}
If $\bo{x},\bo{y} \in \R^2$ satisfy $\ang \big( \bo{x}, \bo{y} \big) \leq c_0$   and $\bo{y} \in \clu{\bo{x}} \cup \clum{\bo{x}}$, then $\vert\vert \bo{y} - \bo{x}  \vert\vert \leq \csc \big( \qzero/2 \big) \vert\vert \bo{x} \vert\vert \ang \big( \bo{x} , \bo{y} \big)$.
\end{lemma}
\noindent{\bf Proof.}
Without loss of generality, $\argu(\bo{y}) < \argu(\bo{x})$.  Let $\bo{q}$ denote the point on $\ell_{\bo{0},\bo{y}}$ closest to $\bo{x}$. Let $\theta$ denote the angle at $\bo{x}$ of the right-angled triangle with vertices $\bo{x}$, $\bo{q}$ and $\bo{y}$. Then  $d \big( \bo{x}, \bo{y} \big) = d \big( \bo{x}, \bo{q} \big) \sec(\theta)$. Now, 
\begin{equation}\label{eqnangt}
\theta = \ang \big( \bo{y} - \bo{x} , \bo{q} - \bo{x} \big) \leq  \ang \big( \bo{y} - \bo{x} , - \perpu{\bo{x}} \big) + 
\ang \big( - \perpu{\bo{x}} , \bo{q} - \bo{x}  \big).  
\end{equation}
We have that $\ang \big( \bo{y} - \bo{x} , - \perpu{\bo{x}} \big) \leq \pi/2 - \qzero$.
The vectors $\bo{q} - \bo{x}$ and $- \perpu{\bo{y}}$ being parallel, 
we have that 
$\ang \big( - \perpu{\bo{x}} , \bo{q} - \bo{x}  \big) = \ang \big( \bo{x},\bo{y} \big) \leq c_0 \leq \qzero/2$. Revisiting (\ref{eqnangt}),
we find that $\theta \leq \pi/2 - \qzero/2$. Thus, $d\big(\bo{x},\bo{y}\big) \leq d\big(\bo{x},\bo{q}\big) \sec \big( \pi/2 - \qzero/2 \big)$. We have that
$d\big(\bo{x},\bo{q}\big) 
= \vert\vert \bo{x} \vert\vert 
\sin \ang \big( \bo{x} , \bo{y} \big)$. Hence, the statement of the lemma. \qed
\end{subsection}
\begin{subsection}{Ornstein-Zernike results for point-to-point connections}\label{secoz}
We recount the statements that we require from the theory \cite{civ} of point-to-point conditioned connections in a subcritical random cluster model. 

We record Theorem $A$ of \cite{civ} in the two-dimensional case:
 \begin{lemma}\label{lemozciv}
Let $P = \P_{\beta,q}$ with $\beta < \beta_c$ and $q \geq 1$. Then
$$
P \Big(  \bo{0} \leftrightarrow \bo{x}  \Big) = 
\big\vert\big\vert  \bo{x}  \big\vert\big\vert^{- \frac{1}{2}}
 \Psi \big( \bo{n_x} \big) 
 \exp \Big\{   - \xi \big( \bo{n_x} \big)   \vert\vert \bo{n_x} \vert\vert \Big\}
   \Big( 1 + o(1) \Big),
$$
uniformly as $\bo{x} \to \infty$. The functions $\Psi$ and $\xi$ are positive, locally analytic functions on $S^1$, and $\bo{n_x} = \frac{\bo{x}}{\vert\vert \bo{x} \vert\vert}$. 
\end{lemma}
The next lemma follows directly from (1.6) of \cite{civ}.
\begin{lemma}\label{lemoznor}
Let $P = \P_{\beta,q}$ with $\beta < \beta_c$ and $q \geq 1$. 
For all $\delta > 0$, there exists $K = K(\delta) \in \N$ and $c = c(\delta) > 0$ such that,
for all $\bo{x},\bo{y} \in \Z^2$,
$$
P \bigg(  
   \gamma'_{\bo{x},\bo{y}} \subseteq
   \Big( W_{\bo{y}-\bo{x},\delta}\big( \bo{x} \big) \cap 
  W_{-(\bo{y}-\bo{x}),\delta}\big( \bo{x} \big)  \Big) \cup B_K\big(\bo{x}\big)
 \cup B_K\big(\bo{y}\big)
\bigg\vert
 \bo{x} \leftrightarrow \bo{y}  \bigg) \geq c,
$$
where $\gamma'_{\bo{x},\bo{y}}$ denotes the common open component of $\bo{x}$ and $\bo{y}$.
\end{lemma}
\begin{definition}\label{deffluc}
Let $\gamma$ denote a connected set containing $\bo{x},\bo{y} \in \Z^2$.
We write \hfff{fluc}
$$
 {\rm fluc}_{\bo{x},\bo{y}}\big( \gamma \big) = \sup \Big\{ d\big( \bo{z} , \big[ \bo{x},\bo{y} \big] \big): \bo{z} \in \gamma \Big\}.
$$
 \end{definition}
We require bounds on moderate fluctuations of conditioned connections:
\begin{lemma}\label{lemmdf}
Let $P = \P_{\beta,q}$ with $\beta < \beta_c$ and $q \geq 1$. 
There exists a constant $c > 0$ such that, for all $\xo,\yo \in \Z^2$ and $0 < t < c \vert\vert \xo - \yo
\vert\vert^{1/2}$,
$$
P \Big(   {\rm fluc}_{\xo,\yo}\big( \gamma'_{\xo,\yo} \big) 
 \geq  \vert\vert \xo - \yo
\vert\vert^{1/2} t   \Big\vert \xo \leftrightarrow \yo \Big)
 \leq \exp \big\{ - c t^2 \big\},
$$
where $\gamma'_{\xo,\yo}$ again denotes the common open component of $\xo$ and $\yo$. 
\end{lemma}
\noindent{\bf Proof.} 
Define a 
$(\bo{x},\bo{y},\delta,K)$-connection regeneration site of $\gamma_{\xo,\yo}$
to be a vertex $\bo{v} \in V\big(\overline\gamma_{\xo,\yo}\big)$ such that
$$
\gamma'_{\xo,\yo}  \setminus B_K (\bo{v}) \subseteq W_{-(\bo{y}-\bo{x}),\delta} \big( \bo{v} \big) \cup 
 W_{\bo{y}-\bo{x},\delta} \big( \bo{v} \big).
$$
Then it follows directly from (1.8) of \cite{civ} that, for any $\delta> 0$, there exists $K = K(\delta)$ such that, for $s > C \log \dist \bo{x} - \bo{y} \dist$, except with conditional probability $\exp \big\{  - c s \big\}$, every $\bo{w} \in \gamma_{\xo,\yo}$ is within distance $s$ of such a connection regeneration site.
If ${\rm fluc}_{\xo,\yo}\big( \gamma'_{\xo,\yo} \big) 
 \geq  \vert\vert \xo - \yo
\vert\vert^{1/2} t$, we may thus suppose that
some such connection regeneration site $\bo{r}$ is at distance at least  $(t/2)\vert\vert \bo{x} - \bo{y} \vert\vert^{1/2}$ from $\big[ \bo{x_0},\bo{y_0} \big]$. 
By Lemma \ref{lemkapab}, we have that the $P$-probability that $\xo,\yo$ and $\bo{r}$ lie in a common open cluster of which $\bo{r}$ is a connection regeneration site is $\Theta \big( p_{\xo,\bo{r}}p_{\bo{r},\yo} \big)$, where $p_{\bo{x},\bo{y}}= P \big( \bo{x} \leftrightarrow \bo{y} \big)$. 
We apply Lemma \ref{lemozciv} to bound $p_{\xo,\bo{r}}$ and $p_{\bo{r},\yo}$, 
and use the strict convexity of $\xi$ (Lemma \ref{lemozstr}), to  find that $p_{\xo,\bo{r}} p_{\bo{r},\yo}/p_{\xo,\yo} \leq \exp \big\{ - c t^2  \big\}$. \qed 
\end{subsection}
\begin{subsubsection}{Large deviations of global distortion}\label{secldgd}
A large deviations' estimate on the macroscopic profile of the conditioned circuit will be valuable.
\begin{prop}\label{propglobdis}
Let $P = \P_{\beta,q}$ with $\beta < \beta_c$ and $q \geq 1$. 
There exists $c > 0$
such that, for any  $\epsilon \in \big( 0, c \big)$, and for all $n \in \N$,
\begin{equation}\label{eqgd}
P \Big( \globdis > \epsilon n \Big\vert \big\vert {\rm INT} \big( \cir \big) \big\vert \geq n^2  \Big)  \leq \exp \big\{ - c \epsilon n \big\}.
\end{equation}
Under this measure,
 $\centre\big( \cir \big) \in {\rm INT} \big( \cir \big)$
except with exponentially decaying probability in $n$. Moreover, (\ref{eqgd}) holds under the conditional measure
$P \big( \cdot \big\vert \acon \big)$. 
\end{prop}
\noindent{\bf Proof.} There exists $c > 0$ such that, for $\delta \in (0,c)$ and $n$ sufficiently high, under $P \big( \cdot \big\vert \vert {\rm INT} \big( \cir \big) \vert \geq n^2 \big)$, 
the convex boundary $\partial {\rm conv} \big( \cir \big)$ has Hausdorff distance from the dilate $n \partial \wulff$ of the Wulff curve of at most $\delta n$, except with conditional probability at most $\exp \big\{ - c \delta n \big\}$. This follows readily from Theorem 5.7 of \cite{alexcube}. As such, by Lemma \ref{lemozstr}, we may suppose that $\mfl$ is at most $\delta n$. From this, we may argue that, for $C > 0$ a sufficiently high constant, the conditional probability that $\mlr \geq C \delta n$ is at most $\exp \big\{ - \delta n \big\}$.
The author's argument for this assertion takes several paragraphs, but, its being not particularly instructive, we have chosen to omit it. Note that $\globdis \leq {\rm GD}\big( \delconv \big) + \mlr$.
We may thus choose $\delta = \epsilon/(2C)$ to obtain (\ref{eqgd}). The two other assertions follow directly. \qed
\noindent{\bf Remark.} 
We mention that an analogue of Proposition \ref{propglobdis} in dimensions $d \geq 3$ is much more subtle. Proofs of such analogues have been undertaken by \cite{cerf} and \cite{cerfpisztora}. 

The following is an immediate consequence of Proposition \ref{propglobdis}. 
\begin{lemma}\label{lemmac}
There exists $\epsilon > 0$, $\ccone > 0$ and $\cctwo > 0$ such that, for each $n \in \N$,
$$
P \Big(   \cir \subseteq B_{\cctwo n} \setminus B_{\ccone n}    \Big\vert \acon \Big)
   \geq 1 - \exp \big\{ - \epsilon n \big\} 
$$
and
\begin{equation}\label{eqclaim}
P \Big(   \cir \subseteq B_{\cctwo n}     \Big\vert \areacon \Big)
   \geq 1 - \exp \big\{ - \epsilon n \big\}. 
\end{equation}
\end{lemma}
In fixing the value of $\qzero > 0$ in Definition \ref{defrg}, we impose the condition $\qzero \leq \ccone/(2\cctwo)$.
\end{subsubsection}
\begin{subsubsection}{Translating to a centred circuit}
In all of the proofs, we will prove upper bounds on events defined in terms of $\cir$ under the measure 
$P \big( \cdot \big\vert \acon \big)$. The corresponding statement under 
$P \big( \cdot \big\vert \areacon \big)$ will follow from the next result.
\begin{lemma}\label{lemcendisp}
There exist constants $C > c > 0$ such that the following holds. 
Let $\mathcal{M}$ denote a set of circuits that is invariant under translation by any $\bo{v} \in \Z^2$.
Then 
$$
P \Big( \cir \in \mathcal{M} \Big\vert \areacon \Big)
\leq
 C n^2 P \Big( \cir \in \mathcal{M} \Big\vert \acon \Big) + \exp \big\{ - c n \big\}.
$$ 
\end{lemma}
\noindent{\bf Proof.}
It is a simple consequence of (\ref{eqclaim}) and the definition of $\centre(\cir)$
that there exist $C > c > 0$ such that, for each $n \in \N$,
$P \big( \centre(\cir) \not\in B_{Cn} \big\vert \areacon \big) \leq \exp \big\{ - cn \big\}$. Alongside 
 Proposition \ref{propglobdis}, we find that such constants exist that
$$
P \Big( \centre(\cir) \not\in B_{Cn} \cap \intg \Big\vert \areacon \Big) \leq \exp \big\{ - cn \big\}.
$$
By the definition of conditional probability and an application of this bound, we learn that
\begin{eqnarray}
 & & P \Big( \cir \in \mathcal{M} \Big\vert \areacon  \Big) \label{eqclone} \\
 & \leq & \frac{P \Big( \cir \in \mathcal{M} , \areacon ,  \centre(\cir) \in \intg \cap B_{Cn} \Big)}{
 P \big( \areacon \big)} + \exp \big\{ - cn \big\}. \nonumber
\end{eqnarray}
 For a configuration $\omega \in \zoz$ and a lattice point $\bo{v} \in \Z^2$, let $\omega_{\bo{v}} : = \omega \big( \cdot + \bo{v} \big)$ denote the translation of $\omega$ by $- \bo{v}$. Note that, if $\centre\big(\cir\big)(\omega) = \bo{v} \in {\rm INT} \big( \cir(\omega) \big)$, then
$\cir \big( \omega_{\bo{v}} \big) = \cir \big( \omega \big) - \bo{v}$. 
In light of this fact, we find that, for any $\bo{v} \in \Z^2$,
\begin{eqnarray}
 & & P \Big( \cir \in \mathcal{M},  \areacon ,  \centre(\cir) \in \intg, \centre(\cir) = \bo{v}  \Big)
   \nonumber \\
 & \leq &  
  P \Big( \cir \in \mathcal{M},  \areacon ,  \centre(\cir) = \bo{0}  \Big). \label{eqcltwo}
\end{eqnarray}
Note then that there exists $\bo{v} \in Z^2$ such that the following holds.
\begin{eqnarray}
 & & P \Big( \cir \in \mathcal{M},  \areacon ,  \centre(\cir) \in \intg \cap B_{Cn}  \Big) \nonumber \\
 & \leq &  C n^2 P \Big( \cir \in \mathcal{M},  \areacon ,  \centre(\cir) \in \intg, \centre(\cir) = \bo{v}  \Big)
            \nonumber \\
 & \leq &  C n^2 P \Big( \cir \in \mathcal{M},  \areacon ,  \centre(\cir) = \bo{0}  \Big). \nonumber
\end{eqnarray}
We apply this to bound above the numerator of the fraction on the right-hand side of (\ref{eqclone}). Applying the lower bound $P \big( \areacon \big) \geq P \big( \big\{ \areacon \big\} \cap \big\{ \centre(\cir) = \bo{0} \big\} \big)$ to the denominator of this fraction then leads to the statement of the lemma. \qed
\end{subsubsection}
\end{section}
\begin{subsection}{Some comments on the required hypotheses}
Most of the arguments in this paper and its companions \cite{hammondtwo} and \cite{hammondthr} use hypotheses that are a little weaker than insisting that $P = \P_{\beta,q}$, (with $\beta < \beta_c$ and $q \geq 1$), be a subcritical random cluster measure. We now list the basic hypotheses of which we have need and the specific instances at which more is required:
\begin{itemize}\label{assump}
\item $P$ satisfies the ratio-weak-mixing property (\ref{defrwm}),
\item $P$ has exponential decay of connectivity; that is, 
there exists $c > 0$ such that $P_{\omega} \big( \bo{0} \to \partial B_n \big) \leq \exp \big\{ -cn \big\}$ for all $n \in \N$ and $\omega \in \{0,1\}^{E(Z^2) \setminus E(B_n)}$. 
\item $P$ satisfies the bounded energy property (\ref{defbden}),
\item $P$ is translation-invariant,
\item $P$ is invariant under reflection in the coordinate axes.
\end{itemize}
As we have noted, the ratio-weak-mixing and bounded energy properties are satisfied by any $\P_{\beta,q}$, $\beta < \beta_c$. Exponential decay of connectivity follows from $\beta < \beta_c$ and the ratio-weak-mixing property, and the other listed hypotheses are trivially satisfied.
 
Of the listed hypotheses, the last is used only in an inessential way: it is used in Proposition \ref{propexc} to simplify the argument. Note that the FKG inequality is not on the list.

Beyond these hypotheses, there are two instances in the present paper and \cite{hammondtwo} and \cite{hammondthr} at which we require something further. The first is in making use of the theory of conditioned point-to-point subcritical connections, i.e., Lemmas \ref{lemozstr}, \ref{lemozciv}, \ref{lemmdf}, as well as Lemma 2.5 in \cite{hammondthr}. These lemmas rely on the theory of \cite{civ}, in which, at one key moment (the final displayed equation of page 1311), a specific decoupling property of the random cluster measures is used. It would be of interest to attempt to dispense with this property. However, this appears to require a significant alteration of method. Indeed, it is easy to construct examples of models satisfying the above listed assumptions in which the deconstruction of long conditioned connections into irreducible pieces does not result in these pieces satisfying the decay condition (1.8) of \cite{civ}  uniformly under conditioning on all the others. Rather, (1.8) may fail on an exponentially small part of the conditioned space. This absence of uniformity disables the use of Ruelle's theory of spectral operators on which the local limit analysis of \cite{civ} relies. We mention, however, that we do not require the full force of Lemma \ref{lemozciv}  to implement our approach: it suffices to know this statement (and that of Lemma \ref{lemoznor}) up to a factor of polynomial growth in $n$.

The second instance of a use of hypotheses beyond the listed assumptions is in the 
 macroscopic distortion bound Proposition \ref{propglobdis}. This invokes 
Theorem 5.7 of \cite{alexcube}, which, beyond the listed hypotheses, also requires
that each of the conditional distributions $P_\omega$ for $\omega \in \{0,1\}^A$, $A \subseteq E(\Z^2)$, satisfy the FKG inequality (as well as invariance of $P$ under rotation by an angle of $\pi/2$). 
\end{subsection}
\begin{section}{From a long facet to a high area-excess: \\ Proposition \ref{propmscbexc}}\label{secunu}
In this section, we exploit the regeneration structure Theorem \ref{thmmaxrg} to obtain:
\begin{prop}\label{propmscbexc}
On the event $\areacon$, we define the area-excess $\exc$ to be the quantity $\big\vert \intg \big\vert - n^2$.
There exist $c,C > 0$ such that, for $1 \leq t = o\big( n^{4/3} \big)$,
 \begin{eqnarray}
& &P \Big(  \mfl \geq n^{2/3}t  \Big\vert \acon \Big) \nonumber \\
& \leq & n^C 
P \Big(  \exc \geq  c n t^{3/2} \big( \log n \big)^{1/2}  \Big\vert \acon \Big)
 + \exp \Big\{ - c n^{1/3} t^{1/2} \Big\}. \nonumber
\end{eqnarray}
\end{prop}
\noindent{\bf Remark.} Note that, if $t \geq Cn^{1/3}$, a bound that decays exponentially in $n$ holds trivially.  
\begin{subsection}{The sector storage-replacement operation}
We now formally define the surgical procedure that is the main tool in the proof of Proposition \ref{propmscbexc}.
\begin{definition}\label{defsrro}
Let $\bo{x},\bo{y} \in \Z^2$, $\argu(\bo{x}) < \argu(\bo{y})$, be given.
Let $P$ be a given measure on configurations $\zoz$. 
The sector storage-replacement operation \hfff{ssro}
$\sigma_{\bo{x},\bo{y}}$ is a random map 
$$
\sigma_{\bo{x},\bo{y}}: \zoz \to \zoz \times \big\{ 0,1 \big\}^{\axye},
$$
whose law is specified in terms of $\bo{x},\bo{y}$ and $P$.
The output $\sigma_{\bo{x},\bo{y}}(\omega) = \big(\omega_1,\omega_2  \big)$
is given as follows. 

We set $\omega_2 = \omega \big\vert_{E(A_{\bo{x},\bo{y}})}$. We then define
$\omega' \in  \{0,1 \}^{E(A_{\bo{x},\bo{y}})}$ to be a random variable whose law 
is the marginal in $E(A_{\bo{x},\bo{y}})$ of  $P \big(\cdot \big\vert \omega\vert_{E(\Z^2) \setminus
    E(A_{\bo{x},\bo{y}})} \big)$ (specified in Definition \ref{defmarg}).
 We set 
\begin{equation*}
 \omega_1 =
\begin{cases}
 \omega' & \text{on $E(A_{\bo{x},\bo{y}})$,}  \\
 \omega  & \text{on $E(\Z^2) \setminus E(A_{\bo{x},\bo{y}})$.}
\nonumber
\end{cases}
\end{equation*}
That is, in acting on $\omega$, we begin by removing the contents of  $E\big(A_{\bo{x},\bo{y}}\big)$ and storing this information as
$\omega_2$.
We then resample these bonds subject to the untouched
information in the complement  $E(\Z^2) \setminus E(A_{\bo{x},\bo{y}})$. The new configuration,
that coincides with the original one in  $E(\Z^2) \setminus
E(A_{\bo{x},\bo{y}})$, is recorded as $\omega_1$.

We will call $\omega_1 \in \zoz$ the full-plane output, 
$\omega_2 \in \big\{ 0,1 \big\}^{\axye}$ the sector output, and
$\omega_1 \big\vert_{\axye}$ the updated configuration.
\end{definition}
In all its applications, the sector storage-replacement operation will act in the following way.
\begin{definition}\label{defregact}
Let $\bo{x},\bo{y} \in \Z^2$, $\argu(\bo{x}) < \argu(\bo{y})$, be given. 
The sector storage-replacement operation 
$\sigma_{\bo{x},\bo{y}}$ will be said to act regularly if 
\begin{itemize}
 \item the input configuration has the distribution $P$, and
 \item 
the randomness of the action is chosen such that, 
given the input $\omega \big\vert_{E(\Z^2) \setminus \axyoe}$, the updated configuration
$\omega_1 \big\vert_{\axyoe}$ is conditionally independent of the stored sector configuration $\omega_2 \big\vert_{\axye} = \omega \big\vert_{\axye}$.
\end{itemize}
\end{definition}
\end{subsection}
\begin{subsection}{The good area capture lemma}
Before proceeding to the proof of Proposition \ref{propmscbexc}, we make precise one of the objects required, the notion of {\it good area capture}:

\begin{definition}\label{defgac}
Let $\bo{x},\bo{y} \in \Z^2$.
For a configuration $\omega \in \{ 0,1 \}^{\axye}$ such that 
$\bo{x}  \build\leftrightarrow_{}^{\axy} \bo{y}$ under $\omega$,
write \hfff{ogamoop} $\overline{\gamma}_{\bo{x},\bo{y}}$
for the common $\omega$-open cluster of $\bo{x}$ and $\bo{y}$ in $\axy$, and
 \hfff{gamoop}  $\gamma_{\bo{x},\bo{y}}$ for the outermost open path in $\axy$
from $\bo{x}$ to $\bo{y}$ (i.e., that open path $\gamma$ in $\axy$
from $\bo{x}$ to $\bo{y}$ for which $I_{\bo{x},\bo{y}}\big( \gamma \big)$ is maximal under containment). 

Let $\epsilon > 0$.
\hfff{gac}
Let the set ${\rm GAC}\big( \bo{x},\bo{y},\epsilon \big)$ of $\epsilon$-good area capture configurations in $A_{\bo{x},\bo{y}}$ denote the subset of 
$\omega \in \{ 0,1 \}^{\axye}$
such that the following conditions apply:
\begin{itemize}
 \item $\bo{x}  \build\leftrightarrow_{}^{\axy} \bo{y}$ under $\omega$,
 \item 
$\overline\gamma_{\bo{x},\bo{y}} \subseteq \cluh{\bo{x}} \cap \clumh{\bo{y}}$,
 \item  ${\rm diam} \big( \gamma_{\bo{x},\bo{y}} \big) \leq \cgac \vert\vert
\bo{x} - \bo{y} \vert\vert$, and
 \item writing  \hfff{ixy} $I_{\bo{x},\bo{y}}\big( \gamma_{\bo{x},\bo{y}} \big) \subseteq \R^2$ for the bounded component of $\axy  \setminus \gamma_{\bo{x},\bo{y}}$, 
$$
 \Big\vert I_{\bo{x},\bo{y}} \big( \gamma_{\bo{x},\bo{y}} \big) \Big\vert
  \geq \big\vert  T_{\bo{0},\bo{x},\bo{y}} \big\vert + \epsilon
  \vert\vert \bo{x} - \bo{y} \vert\vert^{3/2} \big( \log \vert\vert \bo{x} -
  \bo{y} \vert\vert \big)^{1/2},
$$
where $T_{\bo{0},\bo{x},\bo{y}}$ is specified in Definition \ref{deftriandothers}.
\end{itemize}
We set  ${\rm GAC}\big( \bo{x},\bo{y}\big) =  {\rm GAC}\big( \bo{x},\bo{y},1/10 \big)$.
\end{definition}
\begin{lemma}\label{lemgac}
Let $P = \P_{\beta,q}$, with $\beta < \beta_c$  and $q \geq 1$.  There exists $\clemgac > 0$ and $n_0:(0,\infty) \to (0,\infty)$ such that the following holds. Let $\epsilon > 0$ and let $n \in \N$ 
satisfy $n \geq n_0(\epsilon)$. 
Let $\bo{x},\bo{y} \in \Z^2$ 
satisfy $\arg(\bo{x}) < \arg(\bo{y})$, $\vert\vert \bo{x} \vert\vert, \vert\vert \bo{y} \vert\vert \leq \clemgac n$, $\vert\vert \bo{y} - \bo{x} \vert\vert \geq  \clemgac \log n$, 
\begin{equation}\label{xyincone}
\bo{y} \in  C_{\pi/2 - \qzero}^F \big( \bo{x} \big)
\end{equation}
and
\begin{equation}\label{xyinctwo}
\bo{x} \in  C_{\pi/2 - \qzero}^B \big( \bo{y} \big).
\end{equation}

Let $\omega \in \{0,1\}^{E(\Z^2)
  \setminus \axye}$ be arbitrary. 
Then
$$
P_{\omega} \Big( {\rm GAC} \big( \bo{x}, \bo{y},\epsilon \big)    \Big) \geq n^{-\clemgac \epsilon^2}  P \Big( \bo{x}  \build\leftrightarrow_{}^{\axy} \bo{y} \Big).
$$
\end{lemma}
\noindent{\bf Proof.} The argument is a straightforward application of Lemma \ref{lemozciv}. We provide an explicit event $\outfluc \in \sigma \big\{ \axy \big\}$
that ensures   ${\rm GAC}\big( \bo{x},\bo{y},\epsilon \big)$
and for which the required bound is satisfied.

We will make use of a coordinate frame for $\R^2$ in which $\ell_{\bo{x},\bo{y}}$ is horizontal, with origin equal to $\bo{y}$ and with $\bo{x}$ having positive $x$-coordinate,  and with $\bo{0}$ in the usual coordinates lying in the lower half-plane. 
Set $h = \vert\vert \bo{y} - \bo{x} \vert\vert$. (We omit integer-rounding from our notation, and assume that $h$ and related quantities are integers.) 
Using the new coordinate system, we write $\bo{x_1} = \big( h/4, 10\epsilon \sqrt{h}  ( \log h )^{1/2} \big)$ and 
$\bo{x_2} = \big( 3h/4, 10\epsilon \sqrt{h}  ( \log h )^{1/2} \big)$. Further set $\bo{x_0} = \bo{y}$ and $\bo{x_3} = \bo{x}$. 
\begin{figure}\label{figguide}
\begin{center}
\includegraphics[width=0.45\textwidth]{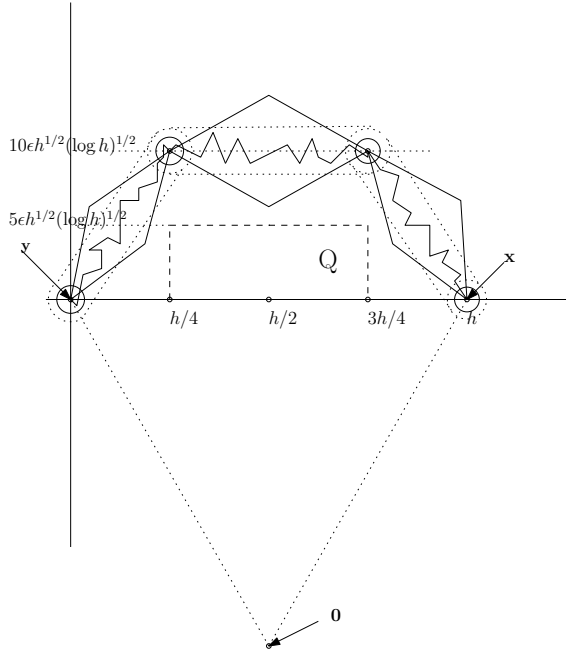} \\
\end{center}
\caption{Illustrating the proof of Lemma \ref{lemgac}. The four circles are the boundaries of the  radius-$K$ balls about the successive points $\bo{y} = \bo{x_0}$ on the left up to $\bo{x} = \bo{x_3}$ on the right. The three diamond-shaped regions each enclosing a line segment $[\bo{x_i},\bo{x_{i+1}}]$  are the $R_i$ for $i=0,1,2,3$ from left to right. 
The region $N$ is the union of the regions bounded by the three dotted curves each of which surrounds one of the line segments  $[\bo{x_i},\bo{x_{i+1}}]$. An instance of $\overline{\gamma}_{\bo{x},\bo{y}}$ that realizes the event $\outfluc$ is depicted.}
\end{figure} 

For $i \in \big\{ 0, 1,2 \big\}$, 
let $R_i =  W_{\bo{x_{i + 1}} - \bo{x_i},\qzero/4} \big( \bo{x_i}  \big) \cap  W_{\bo{x_i} - \bo{x_{i+1}},\qzero/4} \big( \bo{x_{i+1}}\big)$. Note that we may find $K \in \N$ such that there exists an infinite simple (lattice) path from $\bo{0}$ in $B_K \cup W$, for every aperture-$\qzero/2$ cone $W$ with apex at $\bo{0}$. We fix $K \in \N$ at such a value, independently of the value of $h$.
Set $B_0 = B_K(\bo{x_0}) \cap A_{\bo{x},\bo{y}}$, $B_1 =  B_K(\bo{x_1})$,  $B_2 =  B_K(\bo{x_2})$ and  $B_3 = B_K(\bo{x_3}) \cap A_{\bo{x},\bo{y}}$. For such $i$, we denote by $\overline{\gamma}_i$
the common connected component of $\bo{x_i}$ and $\bo{x_{i+1}}$ in  $B_i \cup R_i \cup B_{i+1}$ (if such a component exists), and  
let
$H_i$ denote the event that $\overline{\gamma}_i$ does indeed exist, 
with $\overline{\gamma}_i$  
intersecting  
  $\partial \big( B_i \cup R_i \cup B_{i+1} \big)$ only in $\partial \big( B_i \cup B_{i+1}\big)$ and satisfying
$$
\sup \Big\{ d \big( \bo{v}, \big[  \bo{x_i}, \bo{x_{i+1}} \big] \big): \bo{v} \in \overline{\gamma}_i \Big\}
\leq 10 \vert\vert  \bo{x_{i+1}} - \bo{x_i} \vert\vert^{1/2}.
$$
Let $J$ denote the event that $\overline\gamma_{\bo{x},\bo{y}} \cap B_K(\bo{x}) \subseteq 
C^F_{\pi/2 - \qzero/2}(\bo{x})$ and
 $\overline\gamma_{\bo{x},\bo{y}} \cap B_K(\bo{y}) \subseteq 
C^B_{\pi/2 - \qzero/2}(\bo{x})$. 
 
Set $\outfluc = H_0 \cap H_1 \cap H_2  \cap J$. 

Set $L$ equal to the union of the line segments  
$\big[\bo{x_i},\bo{x_{i+1}} \big]$ for $0 \leq i \leq 2$.
Let $N$ denote the $10 h^{1/2}$-neigbourhood of $L$. Note that $\outfluc$ implies that
 $\bo{x} \build\leftrightarrow_{}^{\axy} \bo{y}$, with $\bo{x_i} \in \overline{\gamma}_{\bo{x},\bo{y}}$ for $i \in \{ 1,2 \}$. Moreover, $\overline{\gamma}_{\bo{x},\bo{y}} \subseteq N$. Hence,
$$
 I_{\bo{x},\bo{y}} \big( L \big) 
\subseteq I_{\bo{x},\bo{y}}\big(\gamma_{\bo{x},\bo{y}}\big) \cup N.
$$
Let $Q$ denote the rectangle that, in the chosen coordinates, has the form 
$Q = \big[ h/4,3h/4 \big] \times \big[ 0, 5 \epsilon h^{1/2} \big( \log h \big)^{1/2} \big]$.
We have that the sets $T_{\bo{0},\bo{x},\bo{y}}$ and $Q$ are disjoint, each of them being a subset of 
$I_{\bo{x},\bo{y}}\big( L \big)$. Noting that $\vert N \vert \leq 20 h^{3/2}$, and that
$\vert Q \vert = (5/2) \epsilon h^{3/2} \big( \log h  \big)^{1/2}$, it follows that, for $n$ high enough,
$$
  \Big\vert I_{\bo{x},\bo{y}}\big( \gamma_{\bo{x},\bo{y}} \big) \Big\vert
 \geq  \big\vert T_{\bo{0},\bo{x},\bo{y}} \big\vert + \epsilon  h^{3/2} \big( \log h  \big)^{1/2},
$$
since $h = \dist \bo{x} - \bo{y} \dist \to \infty$
as $n \to \infty$. To confirm
$\outfluc \subseteq \gacc \big( \bo{x},\bo{y},\epsilon \big)$, it remains to verify the second condition listed in Definition \ref{defgac}. 
To this end, set $\overline{B} = \cup_{i=0}^2 \big( R_i \cup B_i \big) \cup R_3$.
Note that $\outfluc$ entails that $\overline{\gamma}_{\bo{x},\bo{y}} \subseteq \overline{B}$, and hence that
 $\gamma_{\bo{x},\bo{y}} \subseteq \overline{B}$. From (\ref{xyincone}) and (\ref{xyinctwo}) and the definitions of the constituent sets of $\overline{B}$, any point $\bo{u} \in \overline{B} \setminus \big( B_0 \cup B_3 \big)$ 
 is such that each of the angles 
$\ang \big( \bo{u} - \bo{x}, \bo{x}^{\perp} \big)$
 and $\ang \big( \bo{u} - \bo{x}, - \bo{y}^{\perp} \big)$ is at most $\pi/2 - \qzero/2$.
Hence, $\overline\gamma_{\bo{x},\bo{y}}$ satisfies the containment in this second condition except possibly in $B_0 \cup B_3$; but, in this region, the containment is assured by  the occurrence of $J$. Hence, indeed, we have  the second condition listed in Definition \ref{defgac}.

Note that there exist $\clemkac > 0$ and $\clemkam \in \N$ such that  
each of the four following edge-set pairs are  $(\clemkac,\clemkam)$-well separated:
$\big( E(R_i),E(R_{i+1}) \big)$ for $i=0$ and $i=1$, $\big(E(\Z^2) \setminus E(\axy),E(R_0) \big)$, 
and  $\big(E(\Z^2) \setminus E(\axy),E(R_2) \big)$.
It follows readily from Lemma \ref{lemkapab}, and a use of the bounded energy property of $P$ to treat the configuration in the bounded regions $B_i$, that, for a small constant $c > 0$, 
\begin{equation}\label{pombd}
 P_\omega \Big( \outfluc \Big) \geq c P \big( H_0 \big)  P \big( H_1 \big) P \big( H_2 \big).
\end{equation}
We note that, for $i=0,1,2$,
\begin{equation}\label{probhibd}
P \big( H_i \big) \geq \frac{c}{2} P \Big( \bo{x_i} \leftrightarrow \bo{x_{i+1}} \Big).
\end{equation}
Indeed, recalling Definition \ref{deffluc}, we have that 
$$
 \Big\{ \bo{x_i} \leftrightarrow \bo{x_{i+1}} \Big\} 
   \cap \Big\{ {\rm fluc}_{\bo{x_i},\bo{x_{i+1}}} \big( \overline{\gamma}_i \big)  \leq 10 \vert\vert \bo{x_{i+1}} - \bo{x_i} \vert\vert^{1/2} \Big\} \cap \Big\{ \overline{\gamma}_i \in R_i \cup B_K \big( \bo{x_i} \big) \cup      
    B_K \big( \bo{x_{i+1}} \big) \Big\} \subseteq H_i.
$$
Thus, (\ref{probhibd}) follows from Lemmas \ref{lemoznor} and \ref{lemmdf}.

We now bound the three terms $P \big( \bo{x_i} \leftrightarrow \bo{x_{i+1}} \big)$. Set $\bo{e} = \frac{\bo{y} - \bo{x}}{h} \in S^1$. 
Note that $\bo{x_2} - \bo{x_1}
=  \frac{\hvert}{2} \bo{e}$. 
We have that 
$$
\bo{x_1} - \bo{x_0} =  \Big( \bo{e} + \bo{f} \hvert^{-1/2}  \Big)
\sqrt{ (h/4)^2 +  100 \epsilon^2 h \log h},
$$ 
where $\vert\vert \bo{f} \vert\vert \leq C \epsilon \big( \log \hvert \big)^{1/2}$ is such that 
$\bo{e} + \bo{f} \hvert^{-1/2},\bo{e} - \bo{f} \hvert^{-1/2} \in S^1$. 
We also have that
$$
\bo{x_3} - \bo{x_2} =  \Big( \bo{e} - \bo{f} \hvert^{-1/2}  \Big)
\sqrt{ (h/4)^2 +  100 \epsilon^2 h \log h}.
$$  
We find that
$\prod_{i=0}^{2} P \Big( \bo{x_i} \leftrightarrow \bo{x_{i+1}} \Big)$ is at least
\begin{eqnarray}
 & &  \big( \hvert/2 \big)^{-1/2} (h/4)^{-1}  c^3
      \exp \bigg\{ - \bigg(  \xi \Big( e + f  \hvert^{-1/2}  \Big) +
        \xi \Big( e - f  \hvert^{-1/2}  \Big)   \bigg)  \nonumber \\
         & & \qquad \qquad \qquad 
\sqrt{\big( \hvert/4 \big)^2 +   100 \epsilon^2  \hvert \log \hvert} 
 \bigg\}  \
     \exp \Big\{ - \xi \big( \bo{e} \big) \frac{\hvert}{2}  \Big\} \nonumber \\
 & \geq &
      4 \sqrt{2} \hvert^{-3/2}  c^3  
 \exp \Big\{ - C \epsilon^2 \log \hvert \Big\}
      \exp \Big\{ - \xi\big(\bo{e} \big)  \hvert \Big\} \nonumber \\
       & \geq &  4 \sqrt{2} \hvert^{-C \epsilon^2} h^{ -1} c^3 P \Big(  \bo{x} \leftrightarrow \bo{y}  \Big), \label{pprodbd}
\end{eqnarray}
where, to bound the product, we used Lemma \ref{lemozciv} and $\psi \geq c > 0$ uniformly; in the first displayed inequality, we used that $\xi$ is locally analytic; 
and, in the second, Lemma \ref{lemozciv} again. 
By (\ref{pombd}), (\ref{probhibd}), (\ref{pprodbd}), we have that, for $\cgen > 0$ small, 
\begin{equation}\label{probfbd}
P_\omega \Big( \outfluc \Big) \geq \cgen  h^{-C \epsilon} h^{ -1} 
 P \Big(  \bo{x} \leftrightarrow \bo{y}  \Big).
\end{equation}
This conclusion is not sufficient for our purpose, due to the presence of the factor of $h^{-1}$ on the right-hand-side. This term arises as a product of two terms of order $h^{-1/2}$ from the use of Lemma \ref{lemozciv}: in the definition of $\outfluc$, we prescribed the $y$-coordinate (in the present frame) of the points $\bo{x_1}$ and $\bo{x_2}$. We may define a collection of events by varying the $y$-coordinates of $\bo{x_1}$ and $\bo{x_2}$ over intervals of length $h^{1/2}$ centred at the value $10 \epsilon h^{1/2} (\log h)^{1/2}$ used in the definition of $F$. By insisting that successive $y$-coordinates to be used differ by $2K+1$, we ensure the disjointness of the resulting collection. Reprising the argument shows that each event realizes ${\rm GAC} \big( \bo{x},\bo{y},\epsilon \big)$ and satisfies the bound 
(\ref{probfbd}). There being $\big( \cgen \sqrt{h} \big)^2$ such events, we obtain the statement of the lemma. \qed
\end{subsection}
\begin{subsection}{Proof of Proposition \ref{propmscbexc}}
We will use a format for describing the action of the sector storage-replacement 
operation that will be a template for later proofs. \\
\noindent{\bf Definition of satisfactory input and of operation parameters.} 
Let 
$\omega \in \zoz$  denote a configuration. If $\cir = \emptyset$, then the input is not satisfactory. If $\cir \not= \emptyset$, we let $\bo{x},\bo{y} \in \vcir$,
$\arg(\bo{x}) < \arg(\bo{y})$, denote the endpoints
of the longest line segment of which $\delconv$ is comprised. 
An arbitrary deterministic rule should be applied to find the pair 
$\big(\bo{x},\bo{y}\big)$  if the longest line segment is not unique. 
Let $\bo{x'},\bo{y'} \in \reg$, $\arg(\bo{x'}) < \arg(\bo{y'})$,
denote the pair of regeneration sites of $\cir$ for which 
$\axy \subseteq A_{\bo{x'},\bo{y'}}$, and for which $A_{\bo{x'},\bo{y'}}$ is minimal subject to this constraint. That is,
$\bo{x'}$ is the first regeneration site encountered in a clockwise
direction from $\bo{x}$, and $\bo{y'}$ the first such located
counterclockwise from $\bo{y}$. (If, due to $\reg = \emptyset$, this
procedure is not well-defined, then the input is not satisfactory.)

Set ${\rm SAT}_1$ to be the event that
\begin{eqnarray}
 & & \acon
 \cap \Big\{ \mfl \geq n^{2/3} t \Big\} \nonumber \\
 & \cap &  \Big\{ \max \big\{ \vert\vert \bo{x'} - \bo{x} \vert\vert,
  \vert\vert \bo{y'} - \bo{y} \vert\vert  \big\} 
  \leq n^{1/3} t^{1/2}   \Big\}  
\cap \Big\{  \cir \subseteq B_{\cctwo  n}  \setminus B_{\ccone n}
  \Big\}. \label{varevm}
\end{eqnarray}
Note that ${\rm SAT}_1  \subseteq \big\{ \cir \subseteq B_{\cctwo n} \big\}$ permits us to find two points $\bo{x_0},\bo{y_0} \in B_{\cctwo n}$ for which
\begin{equation}\label{xpxoineq}
 P \Big(   \bo{x'} = \bo{x_0}, \bo{y'} = \bo{y_0}  \Big\vert {\rm SAT}_1 \Big) \geq \frac{1}{2 \pi^2 \cctwo^4 n^4}.
\end{equation}
Writing ${\rm SAT}_2 = \big\{ \bo{x'} = \bo{x_0}, \bo{y'} = \bo{y_0}
\big\}$, we declare
the input 
$\omega \in \zoz$ 
to be {\it satisfactory} if it realizes the event ${\rm SAT} : = {\rm
  SAT}_1 \cap {\rm SAT}_2$. \\
\noindent{\bf Specifying operation randomness and the definition of successful action.}
We will apply
$\sigma_{\bo{x_0},\bo{y_0}}$ so that it acts regularly.
The operation will be said to {\it act successfully} 
if the updated configuration $\omega_1 \big\vert_{\axyoe}$ produced by
$\sigma_{\bo{x_0},\bo{y_0}}$ realizes the event {\rm GAC}$(\xo,\yo)$. \\
\noindent{\bf Properties enjoyed by the output.} 
We now claim that, if the input configuration $\omega$ is satisfactory  and the operation $\sigma_{\bo{x_0},\bo{y_0}}$ acts successfully, then the output 
$\sigma_{\bo{x_0},\bo{y_0}}(\omega) = \big( \omega_1,\omega_2 \big)$ has the following properties.
\begin{itemize}
\item {\bf Full-plane circuit property:} the full-plane configuration $\omega_1 \in \zoz$ contains an open
  circuit $\Gamma$ for which $\Gamma \subseteq B_{5 \cctwo n}$, and  
 $$
  \Big\vert  {\rm INT}  \big( \Gamma \big) \Big\vert \geq n^2 + (2/3)^{1/2}
  \frac{1}{40} n t^{3/2} \big( \log n \big)^{1/2};
 $$
\item {\bf Sector open-path property:} and the sector configuration $\omega_2 \in \big\{ 0,1 \big\}^{\axyoe}$
  realizes the event, to be denoted by ${\rm SOPP}$,
that there exists an open path $\gamma$ 
connecting $\xo$ to 
$\yo$ such that $\{ \xo \} \cup \{ \yo \} \in V(\gamma)$,  
$\gamma 
\subseteq \axyo \cap \big( B_{\cctwo n} \setminus B_{\ccone n} \big)$,
$\gamma  \cap W_{\xo,Cn^{-1}\log n}^+ \subseteq \clu{\xo}$ and  
$\gamma  \cap W_{\yo,Cn^{-1}\log n}^- \subseteq \clum{\yo}$. 
\end{itemize}
\noindent{\it Remark.} Note that ${\rm SOPP} \in \sigma \big\{ E(R_{\xo,\yo}) \cap
E \big( B_{\cctwo n} \setminus B_{\ccone n} \big)  \big\}$,
where 
\begin{eqnarray}
R_{\xo,\yo} & = &  \Big( W_{\xo,Cn^{-1}\log n}^+ \cap \clu{\xo} \Big) \cup
   \Big( W_{\yo,Cn^{-1}\log n}^- \cap \clum{\yo} \Big) \nonumber \\
 & & \qquad  \cup \,
 A_{\argu(\xo) + C n^{-1} \log n,\argu(\yo) - C n^{-1} \log n } \, ; \nonumber 
\end{eqnarray}
here we abused the $A$-notation for cones by writing, for $u < v$, $A_{u,v}$ for the cone of points with argument value between $u$ and $v$. See Figure $5$ for an illustration.
Note that there exist constants $\clemkac > 0$ and $\clemkam > 0$ such that, for all $n \in \mathbb{N}$, the two
edge-sets $E \big( R_{\xo,\yo} \big) \cap
E \big( B_{\cctwo n} \setminus B_{\ccone n} \big)$ and $E^* \big( \R^2 \setminus \axyo \big)$
are $(\clemkam,\clemkac)$-well separated, in the sense of Lemma \ref{lemkapab}. (Recall from 
Definition \ref{defpathedge} the notation $E^*$.) We
will later use this lemma to bound the conditional probability of the event
${\rm SOPP}$, given information regarding the configuration on edges that touch
$\R^2 \setminus \axyo$.\\ 
\begin{figure}\label{figsopp}
\begin{center}
\includegraphics[width=0.3\textwidth]{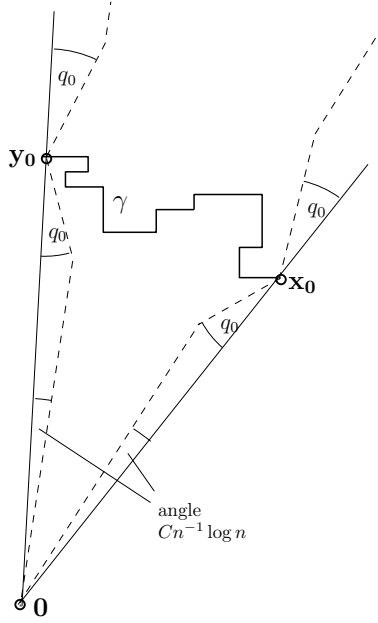} \\
\end{center}
\caption{An instance of the open path $\gamma$ in $\axyo$ when {\rm SOPP} occurs. The path $\gamma $ is contained  in the region $R_{\xo,\yo}$ delimited by the dashed line segments.} 
\end{figure} 
\noindent{\bf Proof of the sector open-path property.}
Note that $\omega_2$ is
equal to the input configuration $\omega\big\vert_{\axyoe}$ in $\axyoe$. 
From $\bo{0} \in \intg$, we find that $\bo{0}$ is separated from $\infty$
in $\axyo \setminus \Gamma_0$. 
Note that the occurrence of ${\rm SAT}_2$ ensures that $\xo,\yo \in \reg$,
and that this, alongside the previous sentence, forces $\cir \cap \axyo$ to
take the form of a path from $\xo$ to $\yo$. Define $\gamma = \cir \cap \axyo$. 
Note that 
$\xo,\yo \in \reg$ 
imply the last two inclusions referred to in the sector open-path property. \\
\noindent{\bf Proof of the full-plane circuit property.}
Let $\gamma_{\bo{x_0},\bo{y_0}} \subseteq \axyo$, 
$\{ \xo \} \cup \{ \yo \} \subseteq \gamma_{\bo{x_0}, \bo{y_0}}$,
denote the $\omega'$-open path from $\xo$ to $\yo$ whose existence is ensured by 
$\omega' : = \omega_1 \big\vert_{\axyoe}$
satisfying ${\rm GAC}(\xo,\yo)$. We set 
$\Gamma = \big( \cir \cap \axyo^c \big) \cup \gamma_{\bo{x_0},\bo{y_0}}$. That is,
$\Gamma$ is the path formed from $\cir$ by replacing $\cir \cap \axyo$
(which touches $\partial \axyo$ just at $\xo$ and $\yo$) with $\gamma_{\xo,\yo}$
(which has the same intersection with $\partial \axyo$). 

We now show that $\Gamma \subseteq B_{5 \cctwo n}$. Note that $\Gamma
\subseteq \cir \cup \gamma_{\bo{x_0},\bo{y_0}}$, that $\cir \subseteq B_{\cctwo n}$ by
assumption, and that $\gamma_{\bo{x_0},\bo{y_0}} \subseteq \xo + B_{4  \cctwo n}$, 
since 
${\rm diam}
\big( \gamma_{\bo{x},\bo{y}} \big) \leq 
2
\vert\vert \xo - \yo \vert\vert \leq 4
\cctwo n$
by means of $\omega'$ realizing ${\rm GAC}(\xo,\yo)$, and by $\xo,\yo \in B_{\cctwo n}$. From $\xo \in  \cir \subseteq B_{\cctwo n}$, we indeed 
obtain $\Gamma \subseteq B_{5 \cctwo n}$.

We will establish the following further properties of $\cir$ and $\Gamma$:
$$
 \Big\vert {\rm INT} (\Gamma) \Big\vert = 
 \Big\vert {\rm INT} (\Gamma) \cap \aarg{x_0}{y_0} \Big\vert
 +   \Big\vert {\rm INT} (\Gamma) \cap \Big( \R^2
 \setminus \aarg{x_0}{y_0} \Big) \Big\vert,
$$
\begin{equation}\label{intinc}
{\rm INT} \big(\Gamma\big) \cap \Big( \R^2
 \setminus \aarg{x_0}{y_0} \Big) = 
 {\rm INT} (\Gamma_0) \cap \Big( \R^2
 \setminus \aarg{x_0}{y_0} \Big), 
\end{equation}
\begin{equation}\label{intaxyo}
 \Big\vert {\rm INT} (\Gamma) \cap \aarg{x_0}{y_0} \Big\vert
 \geq   \big\vert T_{\bo{0},\bo{x_0},\bo{y_0}} \big\vert  +
 \frac{\vert\vert \bo{x_0} - \bo{y_0} \vert\vert^{3/2}}{10} \big( \log
 \vert\vert \xo - \yo \vert\vert \big)^{1/2},
\end{equation}
\begin{equation}\label{propn}
 \Big\vert {\rm INT} (\cir) \cap \aarg{x_0}{y_0} \Big\vert
 \leq   \big\vert T_{\bo{0},\bo{x_0},\bo{y_0}} \big\vert  + 
 \Big( \vert\vert \bo{x} - \bo{y} \vert\vert + 
   \big( 2 \cctwo \ccone^{-1} + 1 \big) \big( 
 \vert\vert \bo{x_0} - \bo{x} \vert\vert +
 \vert\vert \bo{y_0} - \bo{y} \vert\vert \big)  \Big)
 \max \Big\{  
 \vert\vert \bo{x_0} - \bo{x} \vert\vert , 
\vert\vert \bo{y_0} - \bo{y} \vert\vert  \Big\},
\end{equation}
\begin{equation}\label{propnpl}
 \max \Big\{  
 \vert\vert \bo{x_0} - \bo{x} \vert\vert , 
\vert\vert \bo{y_0} - \bo{y} \vert\vert  \Big\}
 \leq 
 \vert\vert \bo{x} - \bo{y} \vert\vert^{1/2}.
\end{equation}
To obtain (\ref{intinc}), we will argue that 
\begin{equation}\label{intlzx}
  {\rm INT} \big( \cir \big) \cap 
 \big( \ell_{\bo{0},\xo}^+ \cup  \ell_{\bo{0},\yo}^+  \big) 
= 
\big[ \bo{0},\xo
  \big] \cup
\big[ \bo{0},\yo
  \big],  
\end{equation}
and that the same statement holds for $\Gamma$.
Indeed, note that, since $\xo,\yo \in \reg$,
\begin{equation}\label{eqcir}
   \cir  \cap 
 \big( \ell_{\bo{0},\xo}^+ \cup  \ell_{\bo{0},\yo}^+  \big) 
= 
\big\{ \xo
  \big\} \cup
\big\{ \yo
  \big\}.  
\end{equation}

From $\bo{0} \in \intg$ and (\ref{eqcir}), it follows that $\cir$ is formed of two sub-paths, one of which crosses $A_{\xo,\yo}$ from $\xo$ to $\yo$, and the other of which crosses $\axyo^c$ from $\yo$ to $\xo$, with each path intersecting $\partial \axyo$ only at $\xo$ and $\yo$. We see then that 
$\Gamma = \big( \cir \cap \axyo^c \big) \cup \gamma_{\xo,\yo}$ also has such a form. The two circuits having such a form justifies both (\ref{intlzx}) and its counterpart for $\Gamma$. 
%

From (\ref{intlzx}), we see that 
$$
 {\rm INT} \big( \Gamma_0 \big) \cap \Big( \R^2 \setminus \axyo \Big)
 = {\rm INT} \bigg(  \Big( \cir \cap \big( \R^2 \setminus \axyo \big) \Big) \cup \big[ \bo{0}, \xo \big]  \cup
 \big[ \bo{0}, \yo \big] \bigg).
$$
We obtain this statement for $\Gamma$ from 
the assertion (\ref{intlzx}) for $\Gamma$. However, 
$ \Gamma \cap \big( \R^2 \setminus \axyo \big) =  \cir \cap \big( \R^2
\setminus \axyo \big)$ by definition, so that (\ref{intinc}) is obtained.

To see (\ref{intaxyo}),
note that
\begin{eqnarray}
  {\rm INT} \big( \Gamma \big) \cap \axyo
   & = & {\rm INT} \Big(  \big( \Gamma \cap \axyo \big) \cup  
 \big[ \bo{0}, \xo \big]  \cup
 \big[ \bo{0}, \yo \big] \Big) \nonumber \\
  & = & {\rm INT} \Big(  \gamma_{\xo,\yo} \cup  
 \big[ \bo{0}, \xo \big]  \cup
 \big[ \bo{0}, \yo \big] \Big), \nonumber 
\end{eqnarray} 
the first equality  by (\ref{intlzx}) for $\Gamma$. Note, however, that 
the final expression coincides with the region 
$I_{\xo,\yo} \big( \gamma_{\xo,\yo} \big)$ appearing in Definition \ref{defgac}.
Thus, (\ref{intaxyo}) follows from $\gamma_{\xo,\yo}$ 
having good area capture.

Turning to (\ref{propn}), let $T$ denote the closed triangle
bounded by the lines $l_{\bo{0}\bo{x_0}}$,  $l_{\bo{0}\bo{y_0}}$ and  $l_{\bo{x}\bo{y}}$.
Note that ${\rm INT}(\cir) \cap \axyo \subseteq T$, because the line
segment from $\bo{x}$ to $\bo{y}$ lies in $\delconv$. Note further that
$\xo,\yo \in \vcir$, $\bo{0} \in \intg$ and 
$\big[ \bo{x}, \bo{y} \big] \subseteq \delconv$ imply that $\xo,\yo \in T$.
Since $T$ is convex, we find that $T_{\bo{0},\xo,\yo} \subseteq T$. We have then
that
$\big\vert {\rm INT}(\cir) \cap \axyo \big\vert \leq \vert
T_{\bo{0},\xo,\yo} \vert + \vert T \setminus T_{\bo{0},\xo,\yo} \vert$. We
will obtain (\ref{propn}) by showing the bound
\begin{equation}\label{tdtx}
  \big\vert T \setminus T_{\bo{0},\xo,\yo} \big\vert  \leq
 \Big( \vert\vert \bo{x} - \bo{y} \vert\vert + 
 \big(2C \ccone^{-1} + 1 \big) 
   \big( \vert\vert \bo{x_0} - \bo{x} \vert\vert +
  \vert\vert \bo{y_0} - \bo{y} \vert\vert \big) \Big)
 \max \Big\{  
 \vert\vert \bo{x_0} - \bo{x} \vert\vert , 
\vert\vert \bo{y_0} - \bo{y} \vert\vert  \Big\}.
\end{equation} 

Consider the quadrilateral $\mathcal{Q} =  T \setminus T_{\bo{0},\xo,\yo}$
whose vertex set $\big\{ \bo{A},\bo{B},\bo{C},\bo{D} \big\}$ is given by
$\bo{A} = l_{\bo{x}\bo{y}} \cap  l_{\bo{0}\yo}$, 
$\bo{B} = l_{\bo{x}\bo{y}} \cap  l_{\bo{0}\xo}$, $\bo{C} = \xo$ and $\bo{D} = \yo$. 
We claim that
\begin{equation}\label{quadare}
\big\vert \mathcal{Q} \big\vert \leq d \big(  \bo{A} , \bo{B} \big)
 \max \Big\{ d \big( \bo{C}, \bo{A}\bo{B} \big)  , d \big( \bo{D}, \bo{A}\bo{B} \big) \Big\}.
\end{equation}
To see this, note that the two sides $\big[\bo{A},\bo{D}\big]$ and $\big[\bo{B},\bo{C}\big]$ of $\mathcal{Q}$,
if continued, would meet at the point $\bo{0}$, on the same side of the line
through $\bo{A}$ and $\bo{B}$ as the whole of $\mathcal{Q}$. This means that, if we
orient $\mathcal{Q}$ so that the side $\big[\bo{A},\bo{B}\big]$ is its base, (namely, so that
$\big[\bo{A},\bo{B}\big]$ lies on the $x$-axis and $\mathcal{Q}$ lies in the upper half-plane),
then the longest intersection that $\mathcal{Q}$ has with a horizontal line
occurs along the side $\big[\bo{A},\bo{B}\big]$. The bound (\ref{quadare}) then arises
because the height of $\mathcal{Q}$ in this coordinate frame is equal to  
$\max \big\{ d \big( \bo{C} , \bo{A}\bo{B} \big)  , d \big( \bo{D} , \bo{A}\bo{B} \big) \big\}$.

Let $\theta = \ang\big( \bo{B} \to \bo{A}, \bo{B} \to \bo{0} \big) = \ang \big( \bo{y} - \bo{x}, - \xo \big)$.
Note that 
\begin{equation}\label{thtangdif}
\theta = \ang \big( \bo{x} \to \bo{y} , -
\bo{x} \big) - \ang \big( \xo, \bo{x} \big).
\end{equation} 
We claim that 
\begin{equation}\label{angtanot}
  \pi - \ccone/\cctwo \geq \ang \big( \bo{x} \to \bo{y}, - \bo{x} \big) \geq 
\ccone/\cctwo.
\end{equation}
To see this, consider the right-angled triangle $\tau$ with vertices
$\bo{0}$ and $\bo{x}$, one of whose sides lies in
$\ell_{\bo{x},\bo{y}}$. Denoting the angle of $\tau$ at the vertex $\bo{x}$ by $\theta'$, note that  
$\theta'  = \min \big\{ \ang \big( \bo{x} \to \bo{y}, - \bo{x} \big) , \pi -  \ang \big( \bo{x} \to \bo{y}, - \bo{x} \big) \big\}$. Let $\bo{z}$
denote the remaining vertex of $\tau$, whose angle is a right angle. We
claim that $\vert\vert \bo{z} \vert\vert \geq \ccone n$. Indeed, $\big[ \bo{x},
\bo{y} \big] \subseteq \delconv$ forces $\cir$ to lie on one side of
$\ell_{\bo{x},\bo{y}}$. From $\bo{0} \in \intg$
and $\cir \cap B_{\ccone n} = \emptyset$, we infer that
$\ell_{\bo{x},\bo{y}} \cap B_{\ccone n} = \emptyset$: for otherwise, a path
might be drawn from $\bo{0}$ to that part of $\partial B_{\ccone n}$ lying
in the half-plane component of $\R^2 \setminus \ell_{\bo{x},\bo{y}}$ that is disjoint
from $\cir$, contradicting $\bo{0} \in \intg$. From $\bo{z} \in
\ell_{\bo{x},\bo{y}}$, we indeed learn that $\vert\vert \bo{z} \vert\vert
\geq \ccone n$.
We now obtain (\ref{angtanot}):
$$
 \theta' \geq \sin \theta' =
 \frac{\vert\vert \bo{z}
  \vert\vert}{\vert\vert \bo{x} \vert\vert} \geq \ccone/\cctwo,
$$
because $\bo{x} \in \vcir$ and the occurrence of ${\rm SAT}_1$ imply that
$\bo{x} \in B_{\cctwo n}$. 

Note also that 
\begin{equation}\label{eqangtwo}
\ang \big( \xo , \bo{x} \big) \leq \frac{\pi}{2} \frac{\vert\vert \xo - \bo{x}
  \vert\vert}{\min \big\{ \vert\vert \xo \vert\vert,\vert\vert \bo{x}
  \vert\vert \big\} } \leq  \frac{\pi}{2 \ccone} n^{-2/3} t^{1/2},
\end{equation}
since ${\rm SAT}_1 \cap {\rm SAT}_2$ implies that
$\vert\vert \xo - \bo{x} \vert\vert \leq n^{1/3} t^{1/2}$
and $\bo{x},\xo  \in \vcir$, $\cir \cap B_{\ccone n} = \emptyset$
imply that $\vert\vert \bo{x} \vert\vert, \vert\vert \xo \vert\vert \geq
\ccone n$.

\begin{figure}\label{figthreetwosix}
\begin{center}
\includegraphics[width=0.6\textwidth]{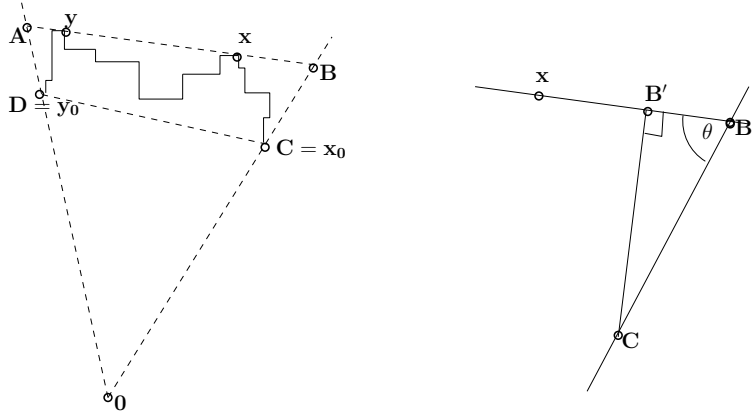} \\
\end{center}
\caption{Deriving (\ref{propn}). The left-hand sketch shows some of the denoted points in the sector $\axyo$.
The right-hand sketch illustrates an instance of case $1$ ($\theta < \pi/2$) considered in the proof of (\ref{tdtx}).} 
\end{figure}

We distinguish the cases that $\theta < \pi/2$ (case $1$) and $\theta \geq \pi/2$ (case $2$). 
Let
$\bo{B'}$
denote the point on  $\ell_{\bo{x}\bo{y}}$ closest to $\xo$. 
Case $1$ holds precisely when $\argu\big(\bo{B'}\big) > \argu\big(\bo{B}\big)$.

We now argue that, in case $1$,
\begin{equation}\label{eqcaseone}
d \big( \bo{x} , \bo{B} \big) \leq \big( 1 + 2 \cctwo/\ccone \big) d \big( \bo{x} , \xo
 \big).  
\end{equation}

To verify this, note that
\begin{equation}\label{xbineq}
 d \big( \bo{x}, \bo{B} \big) \leq d \big( \bo{x}, \bo{B'} \big) + d \big( \bo{B'} , \bo{B} \big). 
\end{equation}  
We have that $d \big( \bo{x}, \bo{B'} \big) \leq d \big( \bo{x} , \xo \big)$
by considering the right-angled triangle with vertex set $\big\{ \bo{x},
\bo{B'}, \xo \big\}$. We have that 
$d \big( \bo{B'} , \bo{B} \big) =  d \big( \xo  , \bo{B'} \big) \cot \theta$.
By (\ref{thtangdif}), (\ref{angtanot}) and (\ref{eqangtwo}), we find that, 
for $n$ high enough, 
$\theta \geq \ccone/(2\cctwo)$. Alongside $\theta < \pi/2$ (we are in case $1$), this implies that
$\cot \theta  \leq \cot \big( \ccone/(2\cctwo) \big) \leq 2\cctwo/\ccone$.
Applying these facts to  (\ref{xbineq}) yields (\ref{eqcaseone}). 

In case $2$, we will show that
\begin{equation}\label{eqcasetwo}
d \big( \bo{x} , \bo{B} \big) \leq d \big( \bo{x} , \xo  \big).  
\end{equation}
Indeed, in this case, the collinear points $\bo{B'}$, $\bo{B}$ and $\bo{x}$ satisfy 
$\argu(\bo{B'}) < \argu(\bo{B}) < \argu(\bo{x})$. Hence, $d(\bo{x},\bo{B})$ is at most $d(\bo{x},\bo{B'})$. However, by considering the right-angled triangle with vertices $\bo{x}$, $\bo{B'}$ and $\bo{x_0}$, 
we see that  $d(\bo{x},\bo{B'}) \leq d(\bo{x},\bo{x_0})$. Hence, we obtain (\ref{eqcasetwo}).

Of course, the weaker inequality (\ref{eqcaseone}) holds in either case; and the same arguments lead to
\begin{equation}\label{eqothercase}
d \big( \bo{y} , \bo{A} \big) \leq \big( 1 + 2 \cctwo/\ccone \big) d \big( \bo{y} , \yo
 \big).  
\end{equation}

Using (\ref{eqcaseone}) and (\ref{eqothercase}),
we find then that
$$
d \big( \bo{A} , \bo{B} \big) \leq d \big( \bo{A} , \bo{y} \big) + d \big( \bo{y} , \bo{x}
\big) + d \big( \bo{x} , \bo{B} \big)
 \leq  d \big( \bo{x} , \bo{y} \big) + \big( 1 + 2\cctwo/\ccone \big) 
 \Big( d \big( \bo{x} , \xo \big) +  d \big( \bo{y} , \yo \big)   \Big).
$$
Noting that 
$d \big( \bo{C} , \bo{A}\bo{B} \big) \leq  d \big( \bo{x} , \xo \big)$
and that 
$d \big( \bo{D} , \bo{A}\bo{B} \big) \leq  d \big( \bo{y} , \yo \big)$, we use 
(\ref{quadare}) to obtain (\ref{tdtx}). 
This completes the
derivation of (\ref{propn}).

To verify (\ref{propnpl}), note that $\bo{x'} = \xo$, $\bo{y'} = \yo$
and ${\rm SAT}_1$ imply that  $\max \big\{  
 \vert\vert \bo{x_0} - \bo{x} \vert\vert , 
\vert\vert \bo{y_0} - \bo{y} \vert\vert  \big\}
 \leq  n^{1/3} t^{1/2}$, while
$\vert\vert \bo{x} - \bo{y} \vert\vert = \mfl \geq n^{2/3}t$ by the
occurrence of ${\rm SAT}_1$. 

We may now bound $\big\vert {\rm INT} (\Gamma)  \big\vert$ from below:
\begin{eqnarray}
 \big\vert {\rm INT} (\Gamma)  \big\vert
 & = & \Big\vert {\rm INT} \big(\Gamma\big)  \cap \Big( \R^2
 \setminus \aarg{x_0}{y_0} \Big) \Big\vert +
 \Big\vert {\rm INT} \big(\Gamma \big)  
    \cap \aarg{x_0}{y_0} \Big\vert \label{intlbd} \\
 & \geq & \Big\vert {\rm INT} \big( \cir \big)  \cap \Big( \R^2
 \setminus \aarg{x_0}{y_0} \Big) \Big\vert 
   +   \big\vert T_{\bo{0},\bo{x_0},\bo{y_0}} \big\vert +
  \frac{\vert\vert \bo{x_0} - \bo{y_0} \vert\vert^{3/2}}{10}  \big( \log
 \vert\vert \xo - \yo \vert\vert \big)^{1/2} \nonumber \\
 & \geq & \Big\vert {\rm INT} \big(\cir \big)  \cap \Big( \R^2
 \setminus \aarg{x_0}{y_0} \Big) \Big\vert
   +  \Big\vert {\rm INT} \big(\cir \big)  \cap \aarg{x_0}{y_0} \Big\vert
   \nonumber \\
 & & \qquad
 -  2 \vert\vert \bo{x} - \bo{y}
   \vert\vert \cdot \vert\vert \bo{x} - \bo{y}
   \vert\vert^{1/2} + \frac{1}{20} \vert\vert \bo{x} - \bo{y}
   \vert\vert^{3/2} \big( \log \vert\vert \bo{x} - \bo{y}
   \vert\vert \big)^{1/2} \nonumber \\
 & \geq & \big\vert \intg \big\vert +  \frac{\vert\vert \bo{x} - \bo{y}
   \vert\vert^{3/2}}{40} \big( \log \vert\vert \bo{x} - \bo{y}
   \vert\vert \big)^{1/2} \nonumber \\
 & \geq & 
\big\vert \intg \big\vert + (2/3)^{1/2} (1/40) n t^{3/2} \big( \log n
\big)^{1/2}, \nonumber
\end{eqnarray}
the first inequality by (\ref{intinc}) and (\ref{intaxyo}),
the second for $n$ high by (\ref{propn}), (\ref{propnpl}) and $\vert\vert \bo{x} - \bo{y} \vert\vert \to \infty$ as $n \to \infty$, and the fourth by 
$\vert\vert \bo{x} - \bo{y} \vert\vert = \mfl \geq n^{2/3} t$ and $t \geq
1$.
This completes the verification of the full-plane circuit property. \\
\noindent{\bf The lower bound on the probability of satisfactory input and successful action.}
%
The input is satisfactory and the operation acts successfully with probability at least
\begin{equation}\label{bdsatcn}
P \Big( {\rm SAT}_1 \Big) \times \frac{1}{2 \pi^2 \cctwo^4 n^4} \times \inf_{\tilde\omega \in \{ 0,1 \}^{E(\Z^2) \setminus \axyoe}}
 P_{\tilde\omega} \Big( {\rm GAC} \big( \xo,\yo \big) \Big).
\end{equation}
The middle term is present due to (\ref{xpxoineq}). 
The third term is present because, given that the input $\omega \in \zoz$
satisfies ${\rm SAT}$, the updated configuration $\omega' \in
\{0,1\}^{\axyoe}$ has law equal to the marginal on $\axyoe$ of 
$\int_{\omega \in \zoz} P \big( \cdot \big\vert \omega\vert_{E(\Z^2) \setminus \axyoe} \big) d P'(\omega)$,
where $P' = P \big( \cdot \big\vert {\rm SAT} \big)$.

We bound the
third term in (\ref{bdsatcn}) 
by
using the good area capture  Lemma \ref{lemgac}. We now check that its hypotheses are satisfied by
$\xo,\yo$. We have that $\dist \xo \dist , \dist \yo \dist \leq \cctwo n \leq \clemgac n$ if we choose $\clemgac \geq \cctwo$. The hypothesis that $d \big( \xo , \yo \big) \geq 4 \clemgac \log n$ is a consequence of 
$\vert\vert \bo{x} - \bo{y}
\vert\vert \geq n^{2/3}t$, $\max \big\{ \vert\vert \xo - \bo{x}
\vert\vert,\vert\vert \yo - \bo{y}
\vert\vert  \big\} \leq n^{1/3} t^{1/2}$ and $t \geq 1$.

The conclusions $\yo \in \clu{\xo}$ and $\xo \in \clum{\yo}$ would, in the
case that $\ang(\xo,\yo) \leq c_0$, follow from $\xo,\yo \in \reg$. However,
we need an alternative argument to handle the case of larger angle. The
two norm inequalities in the preceding paragraph imply that
 \begin{equation}\label{angrone}
\ang \big( \xo \to \yo, \bo{x} \to \bo{y} \big) = O(n^{-1/3}). 
\end{equation}
Note that (\ref{angtanot}) 
implies that     
 \begin{equation}\label{angrtwo}
\ang \big( \bo{x} \to \bo{y},  \bo{x}^{\perp} \big) \leq \pi/2 - \ccone/\cctwo.
\end{equation}
We also have that
\begin{equation}\label{angrthr}
\ang \big( \bo{x} ,  \xo \big) = o \big( 1 \big)
\end{equation}  
in high $n$,
due to $\vert\vert \bo{x} \vert\vert \geq \ccone n$, $\vert\vert \xo -
\bo{x} \vert\vert = O\big( n^{1/3} t^{1/2} \big)$ and $t = o \big( n^{4/3} \big)$.
From (\ref{angrone}), (\ref{angrtwo}) and (\ref{angrthr}), we learn that
$\ang\big( \xo \to \yo,  \bo{x_0^\perp} \big) \leq \pi/2 - \ccone/(2\cctwo)$, whence $\yo \in
\clu{\xo}$ by means of the condition $\qzero \leq \ccone/(2\cctwo)$ that we imposed in setting the value of $\qzero$. Similarly,  $\xo \in \clum{\yo}$. 
The hypotheses of Lemma \ref{lemgac} verified, we apply the lemma 
to bound the third term in 
 (\ref{bdsatcn}). In this way, we find that  (\ref{bdsatcn})  
is at least
\begin{equation}\label{prsatbd}
 P \Big( {\rm SAT}_1 \Big) \times \frac{1}{2 \pi^2 \cctwo^4 n^4} n^{- 10^{-2} \clemgac}  P \Big(
 \bo{x_0} \build\leftrightarrow_{}^{\axyo}  \bo{y_0} \Big).
 \end{equation}
\noindent{\bf The upper bound on the probability of the two output properties.}
 We now find an upper bound on the probability that the procedure has an output with both the full-plane circuit and sector open-path properties. 
Note that, since the input configuration has law $P$, the full-plane configuration $\omega_1 \in \zoz$ in the output also has this law. Thus, the full-plane circuit property is satisfied by the output with probability equal to 
\begin{equation}\label{fpcpub}
P \Big( \exists \, \textrm{open circuit} \, \, \Gamma: \Gamma \subseteq B_{5 \cctwo n}, 
  \big\vert {\rm INT} \big( \Gamma \big) \big\vert \geq n^2 + (2/3)^{1/2}
  (1/40)  n t^{3/2} \big( \log n \big)^{1/2} \Big).
\end{equation}
We wish to argue that this probability is bounded above a similar expression involving the $\bo{0}$-centred circuit $\cir$. 
The circuit $\Gamma$ arising under the full-plane circuit property may not be centred at $\bo{0}$. 
This minor difficulty will arise at the corresponding moment in later arguments. We need to translate and centre such circuits as $\Gamma$. Our device for doing so uses the event $\area{A}$ introduced in Definition~\ref{defarea}.
\begin{lemma}\label{lemcirprop}
Let $A \in \N$. 
For any constant $\cgenbig > 0$, 
and for each $n,A \in \N$ with $A \geq n^2$, we have that
$$
P \Big( \exists \,  \textrm{an open circuit} \, \, \Gamma:  \Gamma \subseteq B_{\cgenbig n},  
   \big\vert {\rm INT}(\Gamma) \big\vert \geq A  \Big) \leq
\cpi  \cgenbig^2 n^2 P \big(  \area{A} \big).
$$
\end{lemma}
\noindent{\bf Proof.} It follows from the exponential decay of connectivity and the ratio-weak-mixing property of $P$ that
\begin{equation}\label{eqonel}
P \Big( \exists \,  \textrm{an open circuit} \, \, \Gamma:  \Gamma \subseteq B_{\cgenbig n},  
   \big\vert  {\rm INT}(\Gamma)  \big\vert \geq A  \Big) 
\leq 2 P \big( M_{2 \cgenbig n} \big),
\end{equation}
where $M_R$ denotes the event that there exists 
an outermost open circuit $\Gamma$ such that  $\Gamma \subseteq B_R$, 
  and $\big\vert {\rm INT}(\Gamma)   \big\vert \geq A$.
Under $M_{2 \cgenbig n}$, select arbitrarily a circuit $\Gamma$ realizing this event. 
It is easy to see that  $\centre(\Gamma) \in B_{4\cgenbig n}$ under $M_{2 \cgenbig n}$. 
Let $\bo{v} \in  B_{4\cgenbig n}$.
On the event 
$M_{2 \cgenbig n} \cap \big\{ \centre(\Gamma) = \bo{v} \big\} \cap \big\{ \bo{v} \in {\rm INT}(\Gamma) \big\}$, 
the shifted configuration $\omega_{\bo{v}} : = \omega \big( \bo{v} + \cdot \big)$ 
realizes $\centre(\cir) = \bo{0}$ and thus $\area{A}$. On the other hand, the event 
$M_{2 \cgenbig n} 
\cap 
\big\{ 
\centre(\Gamma) 
\not\in {\rm INT}(\Gamma) \big\}$
entails the presence in $B_{2 \cgenbig n}$ of an open circuit that encloses a region whose area is least $A$ but which does not contain its own centre. This eventuality is much less probable than the occurrence of $\area{A}$, as we see from Proposition \ref{propglobdis} and the bound $A \geq n^2$.   
Hence, $P\big( M_{2\cgenbig n} \big) \leq \pi \big( 4 \cgenbig n \big)^2 P \big( \area{A} \big) \big( 1 + o(1) \big)$, as required. \qed
By Lemma \ref{lemcirprop} and (\ref{fpcpub}), we find that that the full-plane circuit property is satisfied by the output with probability at most 
\begin{equation}\label{proponebd}
\cpi \big( 5 \cctwo \big)^2 n^2    P \Big( 
    \big\vert \intg \big\vert \geq n^2 +  (2/3)^{1/2}
 (1/40)  n t^{3/2} \big( \log n \big)^{1/2}   \Big).
\end{equation}
Recall that we apply $\sigma_{\xo,\yo}$ so that it acts regularly (as specified in Definition \ref{defregact}).  
Hence,
given the occurrence of the full-plane circuit property, the conditional distribution of the sector output 
 $\omega_2$ is given by the marginal in $\axyoe$ of
 $\int_{\tilde{\omega} \in \{0,1\}^{E(\Z^2) \setminus \axyoe}} P_{\tilde{\omega}} d \mu (\tilde{\omega})$, where the measure $\mu$ on the set of configurations
 $\{0,1\}^{E(\Z^2) \setminus \axyoe}$ is given by the marginal on
$E(\Z^2) \setminus \axyoe$ of the measure $P$ conditioned on  the full-plane circuit property. 
Hence, the conditional probability of the sector open-path property, given the full-plane circuit property, is given by $\int_{\tilde{\omega} \in \{0,1\}^{E(\Z^2) \setminus \axyoe}} P_{\tilde{\omega}} \big( {\rm SOPP} \big) d \mu (\tilde{\omega})$, and thus is at most
\begin{equation}\label{ptomf}
 \sup \Big\{ 
  P_{\tilde\omega} \big( {\rm SOPP} \big) : 
\tilde\omega \in \{ 0,1 \}^{E(\Z^2) \setminus \axyoe} \Big\}.
\end{equation}
Noting that $E(\R^2) \setminus \axyoe = E^* \big( \R^2 \setminus \axyo
\big)$,
and recalling the remark following the statement of the sector open-path property, we may
apply Lemma \ref{lemkapab} with the choice of edge-sets $A = E^* \big( \R^2 \setminus \axyo \big)$
and $B = E \big( R_{\xo,\yo} \big)  \cap E \big( B_{\cctwo n} \setminus B_{\ccone n} \big)$
to conclude that the last displayed quantity is at most $\conka P \big( {\rm SOPP} \big)$. This in
turn is at most $\conka P \big( \bo{x_0} \build\leftrightarrow_{}^{\axyo}
\bo{y_0} \big)$, due to ${\rm SOPP} \subseteq \big\{ \bo{x_0} \build\leftrightarrow_{}^{\axyo}
\bo{y_0} \big\}$. 

In summary, in acting on an input with law $P$, $\sigma_{\bo{x_0},\bo{y_0}}$ will return an output having the full-plane circuit property and the sector open-path property with probability at most 
\begin{equation}\label{prsigubd}
  \cpi \big( 5 \cctwo \big)^2 \conka  n^2 P \Big(  \big\vert \intg \big\vert \geq   n^2 +  (2/3)^{1/2}
  (1/40)  n t^{3/2} \big( \log n \big)^{1/2}  \Big)
  P \Big( \bo{x_0} \build\leftrightarrow_{}^{\axyo}  \bo{y_0} \Big).
 \end{equation}
\noindent{\bf Conclusion by comparison of the obtained bounds.} 
However, we have seen that circumstances arise in which such an output will definitely be produced whose probability is at least the quantity given in (\ref{prsatbd}). Thus, the quantity in (\ref{prsatbd}) is at most that in (\ref{prsigubd}). That is, for high $n$,
\begin{equation}\label{satbd}
 P \Big( {\rm SAT}_1 \Big) \leq n^{10^{-2} \clemgac + 7}
  P \Big(  \big\vert \intg \big\vert \geq   n^2 + (2/3)^{1/2}
  (1/40)  n t^{3/2} \big( \log n \big)^{1/2} \Big).
 \end{equation}
Set $\vsat$ equal to
\begin{eqnarray}
& & \acon   \cap  \Big\{ \mfl \geq n^{2/3} t  \Big\} \nonumber \\
 & & \cap
\Big\{ \mar \leq  \sin \big( \qzero/2\big) \cctwo^{-1} n^{-2/3} t^{1/2}   \Big\} \cap  \Big\{
\cir \subseteq B_{\cctwo n}  \setminus  B_{\ccone n} \Big\}.
\end{eqnarray}
We have that $\vsat \subseteq {\rm SAT}_1$. Indeed,
$\mar \leq  \sin \big( \qzero/2\big)  \cctwo^{-1} n^{-2/3} t^{1/2}$ implies that the vertices 
$\bo{x'}$ and
$\bo{y'}$ are well-defined
and satisfy $\argu(\bo{x}) -  \sin \big( \qzero/2\big)  \cctwo^{-1} n^{-2/3} t^{1/2} \leq
\argu\big(\bo{x'}\big) \leq \argu(\bo{x})$ and 
$\argu(\bo{y}) \leq \argu(\bo{y'}) \leq \argu(\bo{y}) + \sin \big( \qzero/2\big)  \cctwo^{-1} n^{-2/3}
  t^{1/2}$. 
Note that $\bo{x} \in
 \clu{\bo{x'}}$, because $\bo{x'} \in \reg$, $\bo{x} \in \vcir$ and  
$\argu(\bo{x'}) \leq \argu(\bo{x}) \leq \argu(\bo{x'}) +   \sin \big( \qzero/2\big)  \cctwo^{-1} n^{-2/3}
  t^{1/2} \leq \argu(\bo{x'}) + c_0$, the last inequality due to $t
= o\big( n^{4/3} \big)$. From $\argu(\bo{x}) \leq \argu(\bo{x'}) + o(1)$
and $\bo{x} \in \clu{\bo{x'}}$, Lemma \ref{lemdistang} implies that $\vert\vert \bo{x'}
- \bo{x}\vert\vert \leq  \csc \big( \qzero/2 \big)  \vert\vert \bo{x'} \vert\vert \big( \argu(\bo{x})
- \argu(\bo{x'})\big)$. From $\bo{x'} \in \vcir$ and the occurrence of $\vsat$,
we have that $\vert\vert \bo{x'} \vert\vert  \leq \cctwo n$, so that  
$\vert\vert \bo{x'} - \bo{x} \vert\vert \leq  
n^{1/3} t^{1/2}$. The same inequality holds for 
$\vert\vert \bo{y'} - \bo{y} \vert\vert$. We see that, indeed,  $\vsat
\subseteq {\rm SAT}_1$.

Note then that
\begin{eqnarray}
P \big( \vsat \big)  & \geq & 
 P \Big(  \acon   \cap  \Big\{ \mfl \geq n^{2/3}t  \Big\} 
 \Big) \nonumber \\
 & & \qquad - \, P \Big( \Big\{ \acon  \Big\} \cap 
\Big\{ \mar > \sin \big( \qzero/2\big) \cctwo^{-1} n^{-2/3} t^{1/2}   \Big\} \Big)
 \nonumber \\
 & & \qquad  - \,  P \Big( \acon   \cap 
   \Big\{
\cir \not\subseteq B_{\cctwo n} \setminus B_{\ccone}  \Big\} \Big)  \nonumber \\
 & \geq & 
 P \Big(  \acon   \cap  \Big\{ \mfl \geq n^{2/3}t  \Big\} 
 \Big) \nonumber \\
 & & \quad - \, \exp \Big\{ - c \sin \big( \qzero/2\big) \cctwo^{-1}  n^{1/3} t^{1/2} \Big\} 
 P \Big( \acon  \Big)  - \exp \big\{ - c n \big\} 
 P \Big( \acon  \Big), \nonumber 
 \end{eqnarray}
 the second inequality by means of Theorem \ref{thmmaxrg} and 
Lemma \ref{lemmac}.
The second
 inequality requires that $t \leq c^2 \cctwo^2  \csc^2\big( \qzero/2 \big) n^{4/3}$ for the application
 of Theorem \ref{thmmaxrg}.
 
 For $t \geq 1$, $ t = o \big( n^{4/3} \big)$, then,
 $$
P \big( \vsat \big)   \geq  
 P \Big(  \acon  \ \cap  \Big\{ \mfl \geq n^{2/3}t  \Big\} 
 \Big) -  \exp \Big\{ - c n^{1/3} t^{1/2} \Big\} P \Big( \acon \Big). 
 $$
We divide this inequality by $P \big( \acon \big)$; we then apply 
$P \big( \acon \big) \geq c n^{-2} P \big( \intg \big)$, which is a consequence of (\ref{eqclaim}) in Lemma \ref{lemmac}.
Allied with (\ref{satbd}) and $\vsat \subseteq {\rm SAT}_1$, we find that, for
 such $t$,
 \begin{eqnarray}
& &P \Big(  \mfl \geq n^{2/3}t  \Big\vert \acon \Big) \nonumber \\
& \leq & C n^{10^{-2} \clemgac + 9} 
P \Big(  \exc \geq  (2/3)^{1/2} (1/40) n t^{3/2} \big( \log n \big)^{1/2}  \Big\vert \acon \Big)
 + \exp \Big\{ - c n^{1/3} t^{1/2} \Big\}, \nonumber
\end{eqnarray}
as required for the statement of the proposition. \qed
\end{subsection}
\end{section}
\begin{section}{The tail of the area-excess}\label{secarea}
Here, we prove:
\begin{prop}\label{propexc}
There exist constants $\csmall,\cgentw > 0$ 
and 
such that, for $\cgentw \log n  \leq t \leq \csmall n$,
$$
P \Big( \exc \geq tn  \Big\vert \acon  \Big) \leq
 \exp \big\{ - \csmall t \big\}.
$$
\end{prop}
Note that Theorem \ref{thmmflbd} follows from Proposition \ref{propmscbexc}, Proposition \ref{propexc}
and Lemma \ref{lemcendisp}. \\
\noindent{\bf Remark.} For a model of self-avoiding polygons in the first quadrant, geometrically penalized according to length, \cite{hrynivioffe} provides an expression for the probability that the polygon (and the coordinate axes) enclose a high area that is asymptotically sharp up to $1 + o(1)$ terms. The conclusion of Proposition \ref{propexc} would be strengthened, were these techniques adapted to our setting, 
suggesting that a combination of the techniques of the present study, and of those of \cite{hrynivioffe}, might be valuable for approaching such questions as Conjecture \ref{conjone} and its extensions.\\ 
\subsection{A sketch of the proof of Proposition \ref{propexc}}
Why is an area-excess much exceeding $n$ unlikely? Given $\acon$, the circuit $\cir$ is highly likely to have a diameter of order $n$. If the area-excess is $tn$, then we might say that the region $\intg$ is ``too wide'' by roughly $t$, in the sense that the condition $\acon$ would still be satisfied even if we cut $t$ units from the width of the circuit. 
We would like to define an operation making precise this notion of cutting. Here is a first attempt. Consider an input having the law  of subcritical percolation for the sake of this argument. Take the strip of width $t$ whose left-hand border is the $y$-axis, rip it out from the plane and store it, and push the right-hand portion remaining in the plane to the left by $t$ units. The output consists of a configuration in the removed strip, and a full-plane configuration, each of which has the law of the given subcritical percolation, independently. 
Suppose that we define the input $\omega \in \zoz$ of this operation
to be satisfactory if $\omega$ realizes $\acon \cap \big\{ \intg \geq n^2 + nt \big\}$, and if the circuit $\cir$ cuts through each of $\big\{ ( 0, s ): s \geq 0 \big\}$ and $\big\{ ( t, s ): s \geq 0 \big\}$
just once, at points with a common $y$-coordinate (and with a similar condition holding in the lower half-plane). If the input is satisfactory, 
then the full-plane output 
contains a circuit trapping an area of at least $n^2$, because the portions of the original circuit in the left-hand and right-hand portions were successfully reattached after the strip was removed, and the area-loss from the strip-removal is at most $\Theta(nt)$. (It is easy to take care of the implied constant.) The strip output is crossed by two open subcritical paths, above and below the $x$-axis, that are fragments of the input circuit $\cir$, an 
event of probability $\exp \big\{ -ct \big\}$. This extraction of independent randomness in the 
strip output points to the desired conclusion. However, to implement this proposal, we would need to bound the probability of satisfactory input, including the requirement on how the circuit crosses the boundary of the strip. 
We prefer to avoid this, 
and instead alter the definitions. In essence, we now define satisfactory input to be those $\omega \in \zoz$ realizing  $\acon \cap \big\{ \intg \geq n^2 + nt \big\}$.
Exploiting the presence of regeneration sites in the input circuit $\cir$ known by Theorem \ref{thmmaxrg}, we may scan $\cir$ in the upper half-plane,
pushing a vertical strip of width $t$ along it horizontally, until we find a location close to the $y$-axis at which each side of the strip meets the circuit at a unique point, with the $y$-coordinates of this pair of points differing by a small constant multiple of $t$, ($ct$, say). We do the same in the lower-half plane. 
We then rip out from the plane a modified region which follows the shape of
each located strip as the strip in question crosses $\cir$ (as depicted in the upcoming Figure $7$). 
The two sides of the circuit will not necessarily line up perfectly in the full-plane output along the boundary left by the removal of the region. 
Thus, after this removal, we no longer shift the right-hand portion all the way over so that its left boundary touches the right boundary of the left-hand portion, but instead leave a narrow corridor, (which we fill in with independent percolation, in order that the percolation measure is mapped to itself by the operation). 
We define the operation to act successfully if this corridor contains two open paths, one in the upper- and one in the lower-half plane,   
that serve to glue across the corridor the fragments of $\cir$ in the left-hand and right-hand portions, 
so completing a circuit in the output. Satisfactory input and successful action thus give rise to the peculiarity of length-$t$ subcritical connections in the strip output. These are bought by paying for area-excess of $nt$ (to ensure satisfactory input) and the presence of the open gluing connections in the corridor (to ensure successful action). However, the open gluing connections, being of length only $ct$, are much more probable than the length-$t$ strip output open paths. 
Thus, area-excess of order $t$ is as improbable given $\acon$ as the unconditioned probability $\exp\big\{ - c t \big\}$ of a subcritical path of  length $t$. This is the conclusion that we seek.
\subsection{The storage-shift-replacement operation}
Prepared with this motivation, we introduce the operation formally.
\begin{definition}\label{defsro}
Let $P$ be a given measure on configurations $\zoz$. 
Let $\aaa,\bb \subseteq \R^2$ and $\bo{x} \in \Z^2$ be such that $E(\aaa) \cap E(\bb) = \emptyset$
and $E(\aaa) \cap \big( E(\bb) + \bo{x} \big) = \emptyset$. The storage-shift-replacement operation \hfff{stsh}
$$
\phi = \phi_{\aaa,\bb,\bo{x}}: \zoz \to \zoz \times \big\{ 0,1 \big\}^{E(\Z^2) \setminus \big( E(\aaa) \cup E(\bb) \big)}
$$
is the following random map.
Let $\omega \in \zoz$ denote the input of $\phi$. We will define $\phi(\omega) = \big(\omega_1,\omega_2 \big)$.
Firstly, the contents of the input in $E(\Z^2) \setminus \big( E(\aaa) \cup E(\bb) \big)$ are stored as the second component of the output:
$\omega_2 = \omega \big\vert_{E(\Z^2) \setminus \big( E(\aaa) \cup E(\bb) \big)}$.

The full-plane configuration $\omega_1$ in $E(\aaa)$ is set equal to $\omega$: $\omega_1 \big\vert_{E(\aaa)} = \omega \big\vert_{E(\aaa)}$.
The input in $\bb$ is displaced by $\bo{x}$ 
and recorded as $\omega_1$: 
for each $y \in E(\bb)$, we set $\omega_1 ( y  ) =
\omega \big( y - \bo{x} \big)$. The configuration $\omega_1$ is 
completed by assigning 
$\omega_1 \big\vert_{E(\Z^2) \setminus \big( E(\aaa) \cup ( E(\bb) + \bo{x}) \big)}$
to be random, its law being the marginal on
$E(\Z^2) \setminus \big( E(\aaa) \cup ( E(\bb) + \bo{x}) \big)$ of the
  conditional distribution of $P$
given the already assigned values 
$\omega_1 \big\vert_{E(\aaa) \cup \big( E(\bb) + \bo{x} \big)}$

We also define the shift-replacement operation 
$\xi_{\aaa,\bb,\bo{x}}: \zoz \to \zoz$ to be the first component of the map $\phi_{\aaa,\bb,\bo{x}}$.
\end{definition}
We will take the operation to act regularly, in a similar sense as that of Definition \ref{defregact}:
\begin{definition}\label{defregactn}
Let $P$, $\aaa$, $\bb$ and $\bo{x}$ be specified as in Definition \ref{defsro}.
The storage-shift-replacement operation 
$\phi_{\aaa,\bb,\bo{x}}$ will be said to act regularly if 
\begin{itemize}
 \item the input configuration has the distribution $P$, and
 \item 
the randomness of the action is chosen such that, 
given  
$\omega_1 \big\vert_{E(\aaa) \cup \big( E(\bb) + \bo{x} \big)}$, the configuration
$\omega_1 \big\vert_{E(\Z^2) \setminus \big( E(\aaa) \cup ( E(\bb) + \bo{x}) \big)}$
is conditionally independent of the stored configuration 
$\omega_2 = \omega \big\vert_{E(\Z^2) \setminus \big( E(\aaa) \cup E(\bb) \big)}$.
\end{itemize}
\end{definition}
The near-invariance of $P$ under 
$\xi_{\aaa,\bb,\bo{x}}$ follows straightforwardly from Lemma \ref{lemkapab}:
\begin{lemma}\label{lemabx}
Suppose that $P$ satisfies the ratio-weak-mixing property (\ref{defrwm}). There exists $\cgeno > 0$ such that the following holds. Let $\aaa \subseteq \R^2$, $\bb \subseteq B_{\cctwo n}$ and $\bo{x} \in \Z^2$
be such that  $\min \big\{  d(\aaa,\bb), d(\aaa,\bb + \bo{x} ) \big\} \geq \cgeno \log n$. 
Write $\tilde{P}$ for the law on $\zoz$ of $\xi_{\aaa,\bb,\bo{x}}$
applied to an input configuration having the law $P$.
Then we have that
$$
\frac{1}{2} \leq \frac{d \tilde{P}}{d P} \big( \omega \big) \leq 2,
$$
for all $\omega \in \zoz$. 
\end{lemma}
\subsection{Proof of Proposition \ref{propexc}}
By Proposition \ref{propglobdis}, 
Lemma \ref{lemmac} and 
Theorem \ref{thmmaxrg}, 
it suffices to show that there exist constants $c,C,\cprime,  \ctilde > 0$ 
such that, for $C \log n \leq t \leq cn$, 
\begin{eqnarray}
& & P \Big( \exc \geq t n  \Big\vert 
 \acon, \cir \subseteq B_{\cctwo n} \setminus B_{\ccone n},  \nonumber \\
 & & \qquad \qquad \globdis \leq n/\ctilde, \mar \leq
 \frac{t}{\cprime n} \Big) \leq \exp \big\{ - c t \big\}. \nonumber
\end{eqnarray}
We must locate the north-side and south-side strips described in the sketch. In fact, we will not insist that their width be precisely $t$. We will instead locate pairs $\bo{x_1}$ and $\bo{x_2}$ of points on the north-side whose horizontal displacement exceeds $t$, but only by a small constant multiple of $t$, and that have a vertical displacement that is a small constant multiple of $t$. In addition, we will find such points that are $\cir$-regeneration points, close to the $y$-axis. This will force each point to be the unique point of contact of $\cir$ with a vertical line segment reaching to the $x$-axis. A similar pair of points $\bo{y_1}$ and $\bo{y_2}$ will be found on the south-side. Later, we will tear out a strip whose boundary runs through each pair of points, and follow the program proposed in the sketch. The proof of the following lemma appears after the end of the argument. 
\begin{lemma}\label{lemlocate}
Let $\ctwo > 0$ satisfy $\ctwo^{-1} < \ccone \sin \big( \min\{ c_0, \qzero/2  \} \big)$. Set $\cthr = \csc\big( \qzero/2 \big)$ and $\cfour = \cthr \cot \big( 3 \qzero/4 \big)$. 
Let $\cprime > 0$ satisfy $\cprime \geq 4 \cctwo^2 \cthr$, as well as a further condition to be specified during the proof. On the event
\begin{equation}\label{inpev}
\Big\{  \cir \subseteq B_{\cctwo n} \setminus B_{\ccone n},
\globdis \leq n/\ctilde, \mar \leq
 \frac{t}{\cprime n} \Big\}, 
\end{equation}
there exist $\bo{x_1},\bo{x_2} \in \reg$ such that, writing
$\bo{x_i} = \big( x_i(1),x_i(2) \big)$, we have that
\begin{eqnarray}
& & x_1(2),x_2(2) > 0, \qquad  
  - \frac{n}{\ctwo} \leq x_1(1) \leq x_2(1) \leq \frac{n}{\ctwo} \label{xiprop} \\
 & & x_1(1) + t/(4\cctwo) \leq x_2(1) \leq x_1(1) + \Big( \frac{1}{4\cctwo} + \frac{\cthr \cctwo}{\cprime}
 \Big) t, \, \, \, \, 
\big\vert x_2(2) - x_1(2) \big\vert \leq  \frac{3 \cfour \cctwo}{\cprime}  t. \nonumber
\end{eqnarray}
(In recording coordinates such as $x_i(1)$, we dispense with the boldface notation used for points in $\R^2$.)
Under the same event, we also claim that such $\bo{y_1}$ and $\bo{y_2}$
exist that satisfy the same set of conditions, with $y_1(2),y_2(2) < 0$ in
place of the corresponding property in (\ref{xiprop}). 
Moreover, we have that 
\begin{equation}\label{xipropf}
\Big\{ y > 0: \big( x_i(1), y \big) \in \cir \Big\}
 =  \big\{ x_i(2)  \big\},
\end{equation}
for $i=1,2$, 
as well as the similar claim for $\bo{y_1}$ and $\bo{y_2}$ (with $y > 0$
replaced by $y < 0$).
\end{lemma}
Regarding the constants in Lemma \ref{lemlocate}, the important point is that we choose $\cprime$ high enough relative to the others, to ensure the approximate vertical alignment of the ends of circuit fragments that we seek to glue together.

We are ready to specify the parameters and to analyse the action of the storage-shift-replacement operation.
We present the argument 
with a similar structure to that of the proof of Proposition \ref{propmscbexc}.\\
\noindent{\bf Definition of satisfactory input and of operation parameters.}
Let
$$
H_1 = \Big\{ \big(x,y \big) \in \Z^2: y > 0 \, \, \textrm{and} \, \, x \leq
x_1(1)\Big\} \cup
 \Big\{ \big(x,y \big) \in \Z^2: y \leq 0 \, \, \textrm{and} \, \, x \leq
y_1(1)\Big\}.  
$$
Set $\aaa_{\bo{x_1},\bo{y_1}} = H_1 \cap \Big( B_{\cctwo n}  \setminus
B_{\ccone n}  \Big)$. Let 
$$
H_2 = \Big\{ \big(x,y \big) \in \Z^2: y > 0 \, \, \textrm{and} \, \, x \geq
x_2(1)\Big\} \cup
 \Big\{ \big(x,y \big) \in \Z^2: y \leq 0 \, \, \textrm{and} \, \, x \geq
y_2(1)\Big\}.  
$$
Set  $\bb_{\bo{x_2},\bo{y_2}} = H_2 \cap \Big( B_{\cctwo n}  \setminus
B_{\ccone n} \Big)$.
\begin{figure}\label{figareaexc}
\begin{center}
\includegraphics[width=0.8\textwidth]{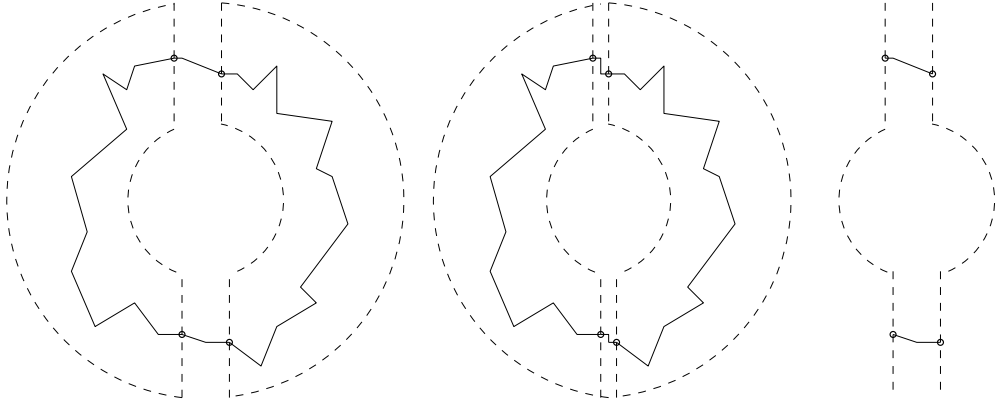} \\
\end{center}
\caption{A satisfactory input for $\phi_{\aaa,\bb,\bo{x}}$ and the full-plane and strip outputs after successful action.}
\end{figure}

We will re-use notation from the proof of Proposition \ref{propmscbexc} for analogous but distinct objects in the present setting.
We set ${\rm SAT}_1$
equal to the intersection of the event in (\ref{inpev}) with
$\big\{ \intg \geq n^2  + nt  \big\}$. 
The fact that $\big\{ \cir \subseteq B_{\cctwo n} \big\} \subseteq {\rm SAT}_1$
ensures the existence of deterministic points
$\bo{x_1^*},\bo{x_2^*},\bo{y_1^*},\bo{y_2^*} \in B_{\cctwo n}$ such that,
writing ${\rm SAT}_2 =  \big\{ \xon = \bo{x_1^*},  \xtw = \bo{x_2^*},     
  \bo{y_1} = \bo{y_1^*},  \bo{y_2} = \bo{y_2^*} \big\}$,
\begin{equation}\label{xstbd}
 P \Big( {\rm SAT}_2   \Big\vert {\rm SAT}_1  \Big)
 \geq  \frac{1}{2} \Big( \frac{1}{\pi \cctwo^2 n^2} \Big)^4.
\end{equation}
A configuration $\omega \in \zoz$ will be called satisfactory input if it realizes the event
${\rm SAT} : = {\rm SAT}_1 \cap {\rm SAT}_2$.

To specify the parameters of $\phi_{\aaa,\bb,\bo{x}}$, set
$\aaa = \aaa_{ \bo{x_1^*},\bo{y_1^*}}$ and  $\bb = \bb_{ \bo{x_2^*},\bo{y_2^*}}$.
With $\cgeno$ the constant from Lemma \ref{lemabx}, we set 
$$
\bo{x} = \Big(  \lfloor -t/(4\cctwo) +  \cgeno \log n  \rfloor, 0 \Big).
$$
Note that, as the sketch proposed, 
$\phi_{\aaa,\bb,\bo{x}}$  pushes the region $\bb$ to the left by almost $t/(4\cctwo)$ units, leaving a small corridor of width $\cgeno \log n$. \\
\noindent{\bf Definition of operation randomness and successful action.}
The operation $\phi_{\aaa,\bb,\bo{x}}$ will be taken to act regularly. 
As a temporary convenience, write $\kappa = E(\Z^2) \setminus \big( E(\aaa) \cup \big(  E(\bb) + \bo{x} \big) \big)$. Let $P_1$ be a path from $\bo{x_1^*}$ to $\bo{x_2^*} + \bo{x}$
whose edges lie in $\kappa$, with the minimal possible number of edges. Let $P_2$ denote a similar path, from $\bo{y_1^*}$ to $\bo{y_2^*} + \bo{x}$. 
The operation will be said to act successfully if the
updated configuration 
$\omega_1 \big\vert_\kappa$ 
satisfies  that $E(P_1) \cup E(P_2)$ is open. (Note that $P_1$ and $P_2$ are the open gluing connections whose use was proposed in the sketch.)\\
\noindent{\bf Properties enjoyed by the output.}
We claim that, in the event of successful action on satisfactory input, the
output of $\phi_{\aaa,\bb,\bo{x}}$ enjoys the
following two properties:
\begin{itemize}
\item \noindent{\bf Full-plane circuit property:} the full-plane configuration $\omega_1 \in \zoz$ 
contains an open circuit $\Gamma$ for which $\Gamma \subseteq B_{(\cctwo+1)n}$ and
$\big\vert {\rm INT} \big( \Gamma \big) \big\vert \geq n^2$,
\item {\bf Strip open-path property:} and the strip configuration 
$\omega_2 \in \{ 0,1 \}^{E(\Z^2) \setminus \big( E(\aaa) \cup  E(\bb) \big)}$ 
realizes the event, to be denoted by ${\rm SOPP}$, that there exists an open path $\gamma \subseteq B_{\cctwo n} \setminus B_{\ccone n}$
connecting  
$\bo{x_1^*}$ to $\bo{x_2^*}$  
and such that 
$\gamma \subseteq  \clum{\bo{x_1^*}}\cap \clu{\bo{x_2^*}}$.
\end{itemize}
\noindent{\it Remark 1.} Regarding the strip open-path property, a similar path exists between 
$\bo{y_1^*}$ and $\bo{y_2^*}$  
and such that 
$\gamma \subseteq  \clus{\bo{y_1^*}}{\bo{y_2^*}}$. However, we will not make use of this second path. \\
\noindent{\it Remark 2.} We have that 
${\rm SOPP}\in \sigma \big\{ E \big( \clum{\bo{x_1^*}}\cap \clu{\bo{x_2^*}} \big)
 \cap E \big( B_{\cctwo n} \setminus B_{\ccone n} \big) \big\}$. 
Note that there exist constants $\clemkac > 0$ and $\clemkam > 0$ such that the pair of
edge-sets $E (\aaa) \cup E(\bb)$ and  
$E \big( \clum{\bo{x_1^*}}\cap \clu{\bo{x_2^*}} \big)
 \cap E \big( B_{\cctwo n} \setminus B_{\ccone n} \big)$ are
$(\clemkam,\clemkac)$-well separated, in the sense of Lemma \ref{lemkapab}. \\
\noindent{\bf Proof of the strip open-path property.} Note that $\omega_2$ is
equal to the input configuration $\omega\big\vert_{E(\Z^2) \setminus \big( E(\aaa) \cup E(\bb) \big)}$ in $E(\Z^2) \setminus \big( E(\aaa) \cup E(\bb) \big)$.
Note that $\bo{x_1^*},\bo{x_2^*} \in \reg$ each satisfy (\ref{xipropf}). By this and 
$\cir \subseteq B_{\cctwo n} \setminus B_{\ccone n}$, $\cir$ contains a subpath from $\bo{x_1^*}$ to $\bo{x_2^*}$ lying in $\R^2 \setminus \big( \aaa \cup \bb \big)$. We choose $\gamma$ to be this subpath. From $\bo{x_1^*},\bo{x_2^*} \not\in B_{\ccone n}$ and  the fact that each of these points satisfies (\ref{xiprop}), we obtain, via $\ctwo^{-1} \leq \ccone \sin (c_0/2)$, that 
$\big\vert \argu(\bo{z}) - \pi/2 \big\vert \leq c_0/2$  for $\bo{z} = \bo{x_1^*},\bo{x_2^*}$, whence $\ang(\bo{x_1^*},\bo{x_2^*}) \leq c_0$. Hence, $\bo{x_1^*},\bo{x_2^*} \in \reg$ implies that  
$\gamma \subseteq  \clum{\bo{x_1^*}}\cap \clu{\bo{x_2^*}}$, as claimed. \\
\noindent{\bf Proof of the full-plane circuit property.}
Let $\gamma_\aaa = \cir \cap \aaa$ and $\gamma_\bb = \cir \cap \bb$. Note that each of $\gamma_\aaa$ and $\gamma_\bb$
is an $\omega$-open path, and that $\gamma_\bb$ connects $\bo{x_1^*}$ to $\bo{y_1^*}$
and  $\gamma_\aaa$ connects $\bo{y_2^*}$ to $\bo{x_2^*}$. Set $\Gamma$
equal to the union of $\gamma_\bb + \bo{x}$, $P_1$, $\gamma_\aaa$ and $P_2$. Then $\Gamma$ is an $\omega_1$-open circuit.  We claim that 
\begin{equation}\label{eqdoubsup}
\sup \big\{ \vert\vert \bo{z} \vert\vert: \bo{z} \in \Gamma \big\} \leq
\sup \big\{ \vert\vert \bo{z} \vert\vert: \bo{z} \in \cir \big\} + \vert\vert \bo{x} \vert \vert.
\end{equation}
Indeed, by the definition of the path $P_1$, the maximal norm among points in $P_1$ is attained at one of its endpoints
$\bo{x_1^*} \in \gamma_\aaa$ or 
$\bo{x_1^*} \in \gamma_\bb + \bo{x}$. A similar property holds for $P_2$. Hence, (\ref{eqdoubsup}) follows, since $\gamma_\aaa \cup \gamma_\bb \subseteq \cir$. From $\vert\vert \bo{x} \vert\vert \leq t/(4\cctwo) + \cgeno \log n \leq \csmall n/(3\cctwo) \leq n$, we see that, indeed, $\Gamma \subseteq B_{(\cctwo+1)n}$.

  Let $E$ denote the region $E = (E^+ \cup E^- ) \cap \big\{ \bo{z} \in \R^2: - \cctwo n \leq z_2 \leq \cctwo n  \big\}$, with
  $E^+ = \big\{ \bo{z} \in \R^2:   x_1^*(1) \leq z_1 \leq x_2^*(1),  z_2 \geq 0  \big\}$
  and 
  $E^- = \big\{ \bo{z} \in \R^2:   y_1^*(1) \leq z_1 \leq y_2^*(1),  z_2  \leq  0  \big\}$.
  From $\cir \subseteq B_{\cctwo n}$, each point in $\intg$ shares a $y$-coordinate with some point in $E$. Let ${\rm LEFT},{\rm RIGHT} \subseteq \intg$ denote the set of $\bo{u} \in \intg$ 
  such that $E$ is encountered to the right, respectively left, of $\bo{u}$ on a horizontal line.
  
Note the following properties:
$$
\intg \subseteq {\rm LEFT} \cup E \cup {\rm RIGHT}, \qquad
{\rm LEFT} \cup \big( {\rm RIGHT} + \bo{x} \big) \subseteq {\rm INT} \big( \Gamma \big),
$$
$$
{\rm LEFT} \cap \big( {\rm RIGHT} + \bo{x} \big) = \emptyset, \, \, \, \textrm{and} \qquad  \vert E \vert \leq tn.
$$
The second of these follows directly from $\bo{x_1^*}$, $\bo{x_2^*}$, $\bo{y_1^*}$ and $\bo{y_2^*}$
satisfying (\ref{xipropf}).  The third is implied by $\vert\vert \bo{x} \vert\vert \leq t/(4\cctwo) \leq \min \big\{ y_2^*(1) - y_1^*(1), x_2^*(1) - x_1^*(1) \big\}$. For the fourth, note that the displacement in $y$-coordinate of any pair of elements in $E$ is at most $2\cctwo n$, while this quantity for the $x$-coordinate is given by 
$\max \big\{ y_2^*(1) - y_1^*(1), x_2^*(1) - x_1^*(1) \big\} \leq \big( (4\cctwo)^{-1} + \cthr \cctwo/\cprime \big)t \leq t/(2\cctwo)$, since $\cprime \geq 4 \cctwo^2 \cthr$.

From these properties and $\big\vert \intg \big\vert \geq n^2 + nt$, it follows that
$$
  \big\vert \intg \big\vert \geq \big\vert {\rm LEFT} \big\vert + \big\vert {\rm RIGHT} \big\vert
   \geq \big\vert \intg \big\vert - \vert E \vert \geq n^2,
$$  
as required. \\
\noindent{\bf The lower bound on the probability of satisfactory input and successful action.}
Note that
$$
 P \big( {\rm SAT}_1 \big) 
 \geq \frac{1}{2} P \Big( \big\vert \intg \big\vert \geq n^2 + nt  \Big)
  $$
by Proposition \ref{propglobdis}, Lemma \ref{lemmac},
and Theorem \ref{thmmaxrg}.
We find then 
that the input is satisfactory with probability at least $2^{-1} \pi^{-4} 
\cctwo^{-8} n^{-8}
P \big( \intg \geq n^2 + nt \big)$. 

To bound the probability of successful action of the operation, note that
\begin{eqnarray}
& &  \big\vert E(P_1) \big\vert = d_{\ell_1} \Big(  \bo{x_1^*},\bo{x_2^*} + \big(  \lfloor - t/2 + 
\cgeno \log
n \rfloor , 0 \big) \Big) \nonumber \\
 & \leq & \big\vert x_1^*(1) -  x_2^*(1) + t/(4C) - \cgeno \log n
   \big\vert + 1 +  \big\vert x_1^*(2) - x_2^*(2) \big\vert \nonumber \\
    &
  \leq & \frac{\cthr \cctwo}{C'}t + \cgeno \log n + 1 + \frac{3 \cfour \cctwo}{C'} t, \nonumber
\end{eqnarray}
since $\big(  \bo{x_1^*} , \bo{x_2^*} \big)$ is in the range of the random variable 
$\big( \xon,\xtw \big)$ satisfying the properties listed in
(\ref{xiprop}). The same bound holds for $\big\vert E(P_2) \big\vert$.

The operation $\phi_{\aaa,\bb,\bo{x}}$ acting regularly, 
 the conditional probability of successful action, given satisfactory input, is at least
\begin{eqnarray}
& & \inf  \Big\{ P_{\tilde{w}}
 \Big( E(P_1) \cup E(P_2) \, \, \textrm{is open}  \Big): \tilde\omega \in \{0,1\}^{E(\aaa) \cup (E(\bb) + \bo{x})}  \Big\}
 \nonumber \\
 & \geq & \cposen^{\vert E(P_1) \vert + \vert E(P_2) \vert} \geq
  \cposen^{\frac{2 \cthr \cctwo}{C'}t  +   \frac{6 \cfour \cctwo}{C'}  t     + 2\cgeno \log n + 2}, \nonumber
\end{eqnarray}
the displayed inequalities by the bounded energy property of $P$ and the preceding bounds on $\big\vert E(P_1) \big\vert$ and $\big\vert E(P_2) \big\vert$.

In summary, the probability of satisfactory input and successful action is at least
\begin{equation}\label{satsucprob}
P \Big( \big\vert \intg \big\vert  \geq n^2 + nt \Big)
 2^{-1} \pi^{-4} 
\cctwo^{-8} n^{-8} 
  \cposen^{\frac{2 \cthr \cctwo}{C'}t  +   \frac{6 \cfour \cctwo}{C'}  t     + 2C \log n + 2}.
 \end{equation}
\noindent{\bf The upper bound on the probability of the two output properties.}
 The hypotheses of Lemma \ref{lemabx} being satisfied for the present choice of $\aaa,\bb$ and $\bo{x}$, the law 
 $\tilde{P}$ of the first component of the output
satisfies $1/2 \leq \frac{d \tilde{P}}{d P} \big( \omega \big) \leq 2$
for all $\omega \in \zoz$.  
Hence, the output satisfies the full-plane circuit property
with probability at most
$$
2 P \Big( \exists \, \textrm{an open circuit $\Gamma$: $\Gamma \subseteq B_{(\cctwo+1)n}$, $\aconv$}
\Big) \leq \cpi (2\cctwo)^2 n^2 P \Big( \acon \Big).
$$
The inequality is due to Lemma \ref{lemcirprop}. 

By the analogous argument to that presented in the paragraph following (\ref{proponebd}), the conditional probability of the strip open-path property, given the full-plane circuit property, is at most
$\sup_{\tilde\omega \in \{ 0,1\}^{E(\aaa)\cup E(\bb)}} P_{\tilde\omega}\big( {\rm SOPP}\big)$.
By the second remark after the statement of the strip open-path property, and by Lemma \ref{lemkapab}, this quantity is at most 
$\conka P \big( {\rm SOPP} \big) 
\leq  
\conka P \big( \bo{x_1^*} \leftrightarrow \bo{x_2^*} \big)
 \leq 
 \conka \exp \big\{ - ct \big\}$, the final inequality due to  
$\vert\vert   \bo{x_1^*} -  \bo{x_2^*} \vert\vert \geq t/(4\cctwo)$ and the exponential decay of connectivity satisfied by $P$. 

The probability, in acting on an input having the law $P$, that the output
satisfies these two properties, is therefore at most
\begin{equation}\label{nsqbd}
  \cpi (\cctwo + 1)^2 n^2 P \Big( \acon \Big) \conka \exp \big\{ -ct \big\}.
\end{equation}
\noindent{\bf Conclusion by comparing the obtained bounds.} This circumstance assuredly arising with probability at least
(\ref{satsucprob}), we find that
(\ref{satsucprob}) is at most 
(\ref{nsqbd}). Rearranging,
\begin{eqnarray}
 & & P \Big( \exc \geq nt  \Big\vert \acon \Big) \nonumber \\
 & \leq &  40 \pi^5 \cctwo^8  n^{10} (\cctwo + 1)^2 \conka 
   \cposen^{-\frac{2 \cthr \cctwo}{\cprime}t  -   \frac{6 \cfour \cctwo}{\cprime}  t     - 2\cgeno \log n - 2}
  \exp \big\{ - ct \big\}. \label{probestf}
\end{eqnarray}
Choosing $C' > 0$ sufficiently high gives that
$$
P \Big( \exc \geq nt  \Big\vert \acon \Big)  \leq n^{3C\log(c^{-1})} \exp \big\{ - c t \big\},
$$
from which the statement of the proposition follows by means of Lemma \ref{lemcendisp}. \qed
\noindent{\it Remark.}
Note that, in (\ref{probestf}), the competition between the decaying  $\exp \big\{ -ct \big\}$  and  the exploding $\cposen^{-\frac{2 \cthr \cctwo}{\cprime}t  -   \frac{6 \cfour \cctwo}{\cprime} t}$  terms is resolved in favour of the decaying term by choosing $\cprime > 0$ to be high enough. This is just as we proposed in the sketch: the diagnozable peculiarity of long subcritical connections in the stored ``strip'' part of the output is a rarer thing than the price of buying gluing path connections, if we insist that our pairs of points $\bo{x_1}$ and $\bo{x_2}$ (and $\bo{y_1}$ and $\bo{y_2}$) 
have a much smaller vertical than horizontal displacement. \\ 
\noindent{\bf Proof of Lemma \ref{lemlocate}.}
We begin by showing (\ref{xipropf}). To this end, note that, from the bound satisfied by $\ctwo$, we have that
\begin{equation}\label{argminbd}
 \Big\{ \big( x,y \big) \in \R^2: \vert x \vert \leq n/\ctwo, y > 0 \Big\} \cap B_{\ccone n}^c
 \subseteq \Big\{ \bo{z} \in \R^2: \big\vert \argu(\bo{z}) - \pi/2 \big\vert \leq \min \big\{ c_0,\qzero/2 \big\}   \Big\}.
\end{equation}
Hence, it follows from (\ref{xiprop}) that
\begin{equation}\label{xyincwb}
 \Big\{ \big( x_1(1) , y \big) : y > 0 \Big\}
 \subseteq    \cco{\xon}   \cup   B_{\ccone n}.
\end{equation}
Intersecting (\ref{xyincwb}) with $\cir$, the second set in the union on the right-hand-side 
is empty by
assumption. By $\xon \in \reg$,
$$
\cir \cap \cco{\xon} \subseteq \clus{\xon}{\xon}.
$$
Note that 
$$
\clum{\xon} \cap \Big\{ \big( x_1(1),y \big): y > 0 \Big\}
 = \big\{ \xon \big\}
$$
(and similarly for $\clu{\xon}$), since $\big\vert \argu(\xon) - \pi/2
\big\vert < \qzero$ by (\ref{argminbd}). Thus, $\cir \cap \big\{ (x_1(1),y): y > 0 \big\} =
\big\{ \xon \big\}$, which is (\ref{xipropf}) for $\xon$. The other three assertions
follow identically.

We now establish the existence of $\xon$ and $\xtw$. To locate this
pair of points, we begin by using a procedure ${\rm LOCATE}$. 
Procedure ${\rm LOCATE}$ will construct a finite sequence 
\begin{equation}\label{proctoutp}
\Big\{ \big( \bo{u_i}, \bo{v_i} \big) : i \in \big\{ 0,\ldots,M \big\} \Big\}
\end{equation}
of pairs $\big( \bo{u_i}, \bo{v_i} \big) \in \reg^2$ of $\cir$-regeneration
sites. Acting on an input $\omega \in \zoz$ that realizes the event 
(\ref{inpev}), there will exist an $m \in \{ 0,\ldots,M \}$
such that we may take 
$\big( \xon,\xtw \big) = \big( \bo{u_m},\bo{v_m} \big)$, as we will show.
To define procedure ${\rm LOCATE}$, we firstly record as $\big\{
\bo{u_0},\ldots,\bo{u_{M'}} \big\}$ the elements $\bo{u}$ of $\reg$ for
which $u_2 > 0$ and $\vert u_1 \vert \leq n/\ctwo$. 
The list is made in order of increasing $x$-coordinate. 
(There is no ambiguity here, because distinct elements in the list have distinct $x$-coordinates. Indeed, (\ref{xipropf}) applies with $\bo{x_i}$ replaced by any such $\bo{u}$.)
For each $i \in \big\{ 0,\ldots,M' \big\}$, we
define $\bo{v_i}$ to be the element of $\reg$ with $v_i(2) > 0$
and of minimal $x$-coordinate subject to $v_i(1) \geq u_i(1) + t/(4\cctwo)$, (provided, of course, that such a
vertex exists). Let $M \leq M'$ be maximal such that $\bo{v_i}$ exists and
satisfies $v_i(1) \leq n/\ctwo$. In this way, we define the output
(\ref{proctoutp}) of procedure ${\rm LOCATE}$. 

We will establish the following properties of the sequence constructed by procedure
${\rm LOCATE}$, which hold in the case that the input realizes the event
(\ref{inpev}). (The bold-face labels will be used to refer to the properties during their proof.)
\begin{itemize}
\item for $i \in \{0,\ldots,M-1 \}$, {\bf [A1]:} 
$0 \leq u_{i+1}(1) - u_i(1) \leq \frac{\cthr \cctwo}{\cprime}t$,
{\bf [A2]:} $0 \leq v_{i+1}(1) - v_i(1) \leq \frac{2 \cthr \cctwo}{\cprime}t$,
\item  for $i \in \{0,\ldots,M-1 \}$, {\bf [B1]:} 
$\big\vert u_{i+1}(2) - u_i(2) \big\vert \leq \cfour t$,
{\bf [B2]:} $\big\vert v_{i+1}(2) - v_i(2) \big\vert \leq 2 \cfour t$,
\item {\bf [C]:}  for $i \in \{0,\ldots,M-1 \}$, 
$tn/(4C)  \leq v_i(1) - u_i(1) \leq  \big( \frac{1}{4\cctwo} + \frac{\cthr \cctwo}{\cprime} \big) t$,
\item {\bf [D]:}  there exist indices $I_1,I_2 \in \big\{ 0, \ldots, M - 1 \big\}$,
 $I_1 < I_2$, such that $v_{I_1}(2) - u_{I_1}(2) > 0$ and 
$v_{I_2}(2) - u_{I_2}(2) < 0$.
\end{itemize}
\noindent{\bf Derivation of properties A-D.} We begin with a useful fact.
Let $i \in \big\{ 0,\ldots, M - 1 \big\}$.
We claim that
\begin{equation}\label{rgainc}
 \reg \cap A_{\bo{u_{i+1}},\bo{u_i}} =   \big\{ \bo{u_i}, \bo{u_{i+1}} \big\}.
\end{equation}
Indeed, by virtue of $\bo{u_i},\bo{u_{i+1}} \in \reg$ and 
${\ang} \big( \bo{u_i},\bo{u_{i+1}} \big) \leq c_0$ (which we justify
shortly),
$$
\reg \cap  A_{\bo{u_{i+1}},\bo{u_i}} \subseteq 
 \clum{\bo{u_i}} \cap
 \clu{\bo{u_{i+1}}}.
$$
Note that
$\cir \cap B_{\ccone n} = \emptyset$ implies that $\vert\vert \bo{u_i} \vert\vert,\vert\vert \bo{u_i} \vert\vert \geq \ccone n$. From $\vert u_i(1) \vert \leq n/\ctwo$ and (\ref{argminbd}), 
we find that $\big\vert  \argu(\bo{z}) - \pi/2 \big\vert < \qzero/2 \leq \qzero$ for
$\bo{z} = \bo{u_i},\bo{u_{i+1}}$. Hence,
\begin{equation}\label{cbcfinc}
 \clum{\bo{u_i}} \cap \clu{\bo{u_{i+1}}} 
\subseteq \big\{ \bo{u_i} \big\} \cup \big\{ \bo{u_{i+1}} \big\}
 \cup \Big\{ \bo{x} \in \R^2: u_i(1) < x_1 < u_{i+1}(1) \Big\}.
\end{equation}
By the construction of $\bo{u_{i+1}}$, there is no element $\bo{u}$ of
 $\reg$ that satisfies $u(2) > 0$ and $u_i(1) < u(1) < u_{i+1}(1)$.
This establishes (\ref{rgainc}). From (\ref{rgainc}) follows
\begin{equation}\label{marangk}
\mar \geq \ang \big( \bo{u_i}, \bo{u_{i+1}} \big).
\end{equation}   
Hence, the occurrence of (\ref{inpev}) and the assumption that $c \leq c_0 \cprime$ ensures the
promised $\ang\big(\bo{u_i},\bo{u_{i+1}} \big) \leq c_0$. From this, $\bo{u_i} \in \reg$ and $\bo{u_{i+1}} \in \cir$, we obtain $\bo{u_{i+1}} \in \clum{\bo{u_i}} \cup \clu{\bo{u_i}}$. By Lemma \ref{lemdistang}, $\bo{u_i} \in B_{C n}$, (\ref{marangk}) and $(\ref{inpev}) \subseteq \big\{ \mar \leq t/(\cprime n) \big\}$, we obtain
\begin{equation}\label{uionei} 
\vert\vert \bo{u_{i+1}} - \bo{u_i} \vert\vert
\leq
\frac{\cthr \cctwo}{\cprime}t,
\end{equation}
giving properties $A1$ and $B1$.

Let $\bo{w_j}$ denote the element $\bo{w}$ of $\reg$ such that $w(2) > 0$
with $w(1)$ maximal subject to $w(1) \leq u_j(1) + t/(4\cctwo)$. Then
\begin{equation}\label{rgawv}
 \reg \cap A_{\bo{w_j},\bo{v_j}} \subseteq \big\{ \bo{w_i} \big\}
 \cup \big\{ \bo{v_i} \big\}.
\end{equation}
Indeed, note that there is no element $\bo{u}$ of $\reg$ that satisfies $u_2 > 0$ 
and $w_j(1) < u_1 < v_j(1)$. Noting also that 
$w_j(1) \geq u_j(1) \geq -n/\ctwo$, we may obtain (\ref{rgawv}) by
reprising the argument that gives (\ref{rgainc}).
By substituting $\bo{w_{i+1}}$ for $\bo{u_i}$, and $\bo{v_{i+1}}$ for
$\bo{u_{i+1}}$, in the argument leading to (\ref{uionei}), we obtain 
$\vert\vert v_{i+1} - w_{i+1} \vert\vert \leq \frac{\cthr \cctwo}{\cprime}t$ for $n$ high. Hence,
\begin{equation}\label{vutineq}
 v_{i+1}(1) - u_{i+1}(1) - t/(4\cctwo) \leq \frac{\cthr \cctwo}{\cprime}t.
\end{equation}
Note then that $v_{i+1}(1) - v_i(1) \leq v_{i+1}(1) - u_i(1) - t/(4\cctwo) \leq u_{i+1}(1) - u_i(1) + \frac{\cthr \cctwo}{\cprime} t \leq \frac{2 \cthr \cctwo}{\cprime} t$, the last inequality by property $A1$.
The bound $v_{i+1}(1) \geq v_i(1)$ being trivial, we obtain property $A2$.


The vertices $\bo{v_i}$  and $\bo{v_{i+1}}$ belonging to the left-hand-side of (\ref{argminbd}), we have that
$\ang\big( \bo{v_i},\bo{v_{i+1}} \big) \leq c_0$. From 
$\bo{v_i} \in \reg$ and $\bo{v_{i+1}} \in \cir$, we find that $\bo{v_{i+1}} \in \clum{\bo{v_i}} \cup \clu{\bo{v_i}}$. Using $\big\vert \arg\big( \bo{v_i }\big) - \pi/2 \big\vert \leq c_0/2$ (from (\ref{argminbd})) and $c_0 \leq \qzero/2$, we see that 
$\big\vert \arg\big( \bo{v_{i+1}} -  \bo{v_i} \big) - \pi/2 \big\vert \geq \qzero - c_0/2 \geq 3 \qzero/4$.
Hence,  $\big\vert v_{i+1}(2) -  v_i(2) \big\vert \leq \big\vert v_{i+1}(1) -  v_i(1) \big\vert \cot \big( 3\qzero/4 \big)$.

Property $C$ is a reindexed (\ref{vutineq}).

We now establish property $D$. 
We begin by showing that, if (\ref{inpev}) occurs, then
\begin{equation}\label{eqxygam}
\bo{x},\bo{y} \in \cir,
x(2),y(2) > 0, 
\vert x(1) \vert \leq \frac{n}{8\ctwo},
\vert y(1) \vert \in   \Big[\frac{n}{2\ctwo},\frac{n}{\ctwo} \Big] \implies
x(2) > y(2).
\end{equation} 
To this end, we write $\tilde\Gamma : = n \partial
\wulff$, and note that this dilate of the Wulff shape boundary satisfies
$\cir \subseteq \tilde\Gamma + B_{\globdis}$. Let $\bo{x'},\bo{y'} \in
\tilde\Gamma$ satisfy $\vert\vert \bo{x'} - \bo{x}  \vert\vert \leq
\globdis$ and $\vert\vert \bo{y'} - \bo{y}  \vert\vert \leq
\globdis$. By $\globdis \leq n/\ctilde$, $\vert x'(1)  \vert \leq n/(4\ctwo)$ (since $\ctilde \geq 8 \ctwo$) and  
$ \vert y'(1) \vert \in   \big[3n/(8\ctwo),9n/(8\ctwo) \big]$.
By convexity of $\partial \wulff$ and its symmetry in the $y$-axis, we have that  
$$
x'(2) - y'(2)  \geq n \bigg( \wulff^+ \Big( \frac{1}{4\ctwo}\Big) - \wulff  \Big( \frac{3}{8\ctwo}  \Big) \bigg),
$$
where here $\wulff^+:(-c,c) \to [0,\infty)$, $c = \sup \big\{ w_1: \bo{w} \in \wulff \big\}$, denotes the curve $\partial \wulff$ in the
upper-half-plane as a function of the $x$-coordinate. Note that 
$\big\vert x(2) - x'(2) \big\vert, \big\vert  y(2) - y'(2) \big\vert \leq \globdis \leq n/\ctilde$,
so that the desired $x(2) > y(2)$ is ensured by
$$
 \frac{1}{\ctilde} < \frac{1}{2} \bigg( \wulff^+ \Big( \frac{1}{4\ctwo}\Big) - \wulff  \Big( \frac{3}{8\ctwo}  \Big) \bigg).
$$
We now define a partition $\mathcal{P} = \big\{ A_0,A_1,\ldots,A_{Q-1}
\big\}$ of a certain interval whose left-hand endpoint is $u_0(1)$. 
The first
interval $A_0$ in $\mathcal{P}$ is given by $A_0 = \big[ u_0(1), v_0(1) \big]$. 
There exists $j_1 \in \big\{ 0,\ldots,M\big\}$ such that $\bo{u_{j_1}} =
\bo{v_0}$, because the $u$-sequence enumerates 
$\big\{ \bo{u} \in \reg: u_2 > 0, -n/\ctwo \leq u_1 \leq n/\ctwo \big\}$.
Let $A_1 = \big[ u_{j_1}(1),v_{j_1}(1) \big]$. Iteratively construct
$A_n = \big[ u_{j_n}(1),v_{j_n}(1) \big]$. Set $j_{n+1}$ so that
$\bo{u_{j_{n+1}}} = \bo{v_{j_n}}$, and set 
$A_{n+1} = \big[ u_{j_{n+1}}(1), v_{j_{n+1}}(1) \big]$.
Let $Q \in \N$ be such that 
$\big\vert v_{j_{Q-1}}(1) \big\vert = \big\vert u_{j_Q}(1) \big\vert \in \big[ n/(2\ctwo), n/\ctwo 
 \big]$. (Such a $Q$  exists because each constructed interval $A_i$ 
satisfies $\vert A_n \vert \leq n/(16 \ctwo)$, by property $C$, $t \leq \csmall n$ and
the bound $\csmall \leq \big( 1/(4\cctwo) + \frac{\cthr \cctwo}{\cprime} \big)^{-1}  (16 \ctwo)^{-1}$ that we may impose on $\csmall$.)
We stop constructing the sequence $A_i$ at $i = Q -1$. Let 
$R \in \big\{ 0,\ldots, Q \big\}$ be such that 
$\big\vert u_{j_{R + 1}}(1) \big\vert = \big\vert v_{j_R}(1) \big\vert  \leq n/(8\ctwo)$.
By (\ref{eqxygam}),
$$
u_{j_{R + 1}}(2) > \max \Big\{ u_{j_1}(2), u_{j_Q}(2) \Big\},
$$
since $\bo{u_{j_1}} = \bo{v_0}$ satisfies $-n/\ctwo \leq v_0(1) \leq
-n/\ctwo + \big( 1/2 + \frac{2\cctwo}{\cprime}\big)t \leq -n/(2\ctwo)$,
as well as the stated bound on $u_{j_Q}(1)$.

Note that each of the quantities $u_{j_{R + 1}}(2) - u_{j_1}(2) > 0$
and  $u_{j_Q}(2) - u_{j_{R + 1}}(2) < 0$
may be expressed as a telescoping sum of terms of the form 
$v_{j_i}(2) - u_{j_i}(2)$, with summand index sets respectively 
$i = 1,\ldots,R$ and $i = R+1,\ldots,Q-1$. Thus, the indices 
$I_1$ and $I_2$ whose existence is claimed by  property D 
may be chosen so that $I_1 \in \big\{ j_1,j_2,\ldots,j_R \big\}$
and  $I_2 \in \big\{ j_{R+1},j_{R+2},\ldots,j_{Q-1} \big\}$.  

These properties
established, we may now find the pair $\big( \xon,\xtw \big) \in \reg^2$
satisfying (\ref{xiprop}). Indeed, note that the sequence of terms 
$\big\{ v_i(2) - u_i(2): i \in \{ I_1,\ldots,I_2 \} \big\}$
begins at $v_{I_1}(2) - u_{I_1}(2) > 0$,
ends at  $v_{I_2}(2) - u_{I_2}(2) < 0$, and has a difference 
between successive terms (indexed by $i$ and $i+1$)
being in absolute value at most 
$$
\big\vert v_{i+1}(2) - v_i(2) \big\vert +
\big\vert u_{i+1}(2) - u_i(2) \big\vert \leq 3 \frac{\cfour \cctwo}{\cprime} t, 
$$
by property $B$. Setting $I_0$ to be the maximal $i \in \{ I_1,\ldots,I_2 \}$ for
which $v_i(2) - u_i(2) > 0$, we see that 
$\big\vert v_{I_0}(2) - u_{I_0}(2) \big\vert \leq  3 \cfour (\cctwo/C') t$.
We set $\xon = \bo{u_{I_0}}$
and $\xtw = \bo{v_{I_0}}$ to obtain the pair $\big( \xon,\xtw \big)$ as we
sought. \qed
\end{section}
\begin{section}{The upper bound on maximum local roughness}\label{secubdmlr}
In this section, we prove Theorem \ref{thmmlrbd} by making rigorous the argument presented in Subsection \ref{secmlrlog}.
The following definition is convenient.
\begin{definition}
For $\bo{x} \in \vcir$, write $\bo{x'}$ and $\bo{x''}$
for the elements of $\reg$ first encountered in a counterclockwise,
respectively clockwise, search that is centred at $\bo{0}$ and begins at
$\bo{x}$. 
Set the {\it maximum point-to-regeneration-site distance} \hfff{mprg} $\maxseprg$ equal to  
$$
\maxseprg = \max \Big\{ d \big( \bo{x}, \bo{x'} \big),  d \big(
\bo{x}, \bo{x''} \big) : \bo{x} \in \vcir \Big\}.
$$
\end{definition}
\begin{lemma}\label{lemmarbd}
Let $R = \sup \big\{ \vert\vert \bo{x} \vert\vert: \bo{x} \in \cir \big\}$. 
If $\mar \leq 2 c_0$, then 
$$
\maxseprg \leq R \csc(\qzero/2)  \mar.
$$
\end{lemma}
\noindent{\bf Proof.} Let $\bo{x} \in \vcir$.
For $\zeta = \mar$, let $\bo{u} \in \reg \cap W_{\bo{x},\zeta/2}\big( \bo{0} \big)$ be such that no element of $\reg$ has an argument lying strictly between $\argu(\bo{u})$ and $\argu(\bo{x})$. 
Without loss of generality, $\argu(\bo{u}) \leq \argu(\bo{x})$.
From $\argu(\bo{u}) \leq \argu(\bo{x}) \leq \argu(\bo{u}) + \zeta/2 \leq
\argu(\bo{u}) + c_0$, we find that $\bo{x} \in W_{\bo{u},c_0}^+$.
By $\bo{u} \in \reg$, we obtain $\bo{x} \in \clu{\bo{u}}$. By Lemma \ref{lemdistang} and 
$\ang\big( \bo{u},\bo{x}\big) \leq \mar$, we obtain the result. \qed
\noindent{\bf Proof of Theorem \ref{thmmlrbd}.}
Given $\omega \in \zoz$ realizing $\acon$, let $\xmlr$
denote a point $\bo{x} \in \vcir$ such that $\mlr = d \big( \bo{x},
\delconv \big)$. (An arbitrary deterministic rule should be used to select
$\xmlr$ if there is a choice to be made.)
Let $\xmlrm$ and $\xmlrp$ be the extreme points of $\delconv$
that are first encountered from $\xmlr$ in a clockwise, or
counterclockwise, search whose centre is the origin. Let the length 
$\vert\vert  \xmlrp - \xmlrm \vert\vert$ be denoted by \hfff{mlrf} $\mlrs$ 
(the length of the {\it maximum local roughness facet}). 
Note that $\mlrs \leq \mfl$, since the interval 
$\big[ \xmlrm,\xmlrp \big]$ is one of the line segments of which $\delconv$
is comprised. 

The statement of the proposition will be obtained by considering
several cases, as indicated by the following
inclusion, in which $\delta > 0$ is a parameter to be specified later, and where integer-rounding has been omitted for ease of notation:
\begin{eqnarray}
 & & \Big\{ \mlr \geq n^{1/3} \big( \log n \big)^{2/3} t , \acon \Big\}
 \label{mlrinc} \\
 & \subseteq &  \Big\{ \mlrs \geq n^{2/3}\big( \log n \big)^{1/3} t^{2 - \delta} , \acon \Big\}
 \nonumber \\
 & & \quad  
  \cup \, \, \Big(  \bigcup_{j=n^{1/3}}^{n^{2/3} \big( \log n \big)^{1/3} t^{2-\delta}} A_j \Big) 
   \, \, \cup \, \, A  \cup 
   \Big\{ \maxseprg >  n^{1/6}, \cir \subseteq B_{\cctwo n},  \acon\Big\} 
 \nonumber \\
 & &  \qquad 
\cup \, \, \Big\{ \cir
   \not\subseteq B_{\cctwo n} \setminus  B_{\ccone n} , \acon \Big\}   , \nonumber
\end{eqnarray}
where 
\begin{eqnarray}
A_j : & = & \Big\{ \mlr \geq n^{1/3}\big( \log n \big)^{2/3} t ,
 \mlrs = j, \maxseprg \leq n^{1/6}, \nonumber \\
  & & \qquad \acon,  \cir \subseteq B_{\cctwo n}  \setminus  B_{\ccone n} =
 \emptyset \Big\}
\end{eqnarray}
and
$$
A : =  \Big\{ \mlr \geq n^{1/3} \big( \log n \big)^{2/3} t ,
\mlrs < n^{1/3} , 
\maxseprg \leq n^{1/6},
\cir \cap B_{\ccone n} = \emptyset \Big\}.
$$
The case of principal interest is now treated. Integer-rounding is not indicated for ease of notation.
\begin{lemma}\label{lemprinint}
Let $j \in \big\{ n^{1/3},\ldots, n^{2/3} \big( \log n \big)^{1/3}  t^{2 - \delta} \big\}$.
There exists $t_0 > 0$ such that, for $t \geq t_0$,
\begin{equation}\label{probajbd}
P \big( A_j \big) \leq \exp \Big\{ - c \big\{ t^\delta \log n , 
  n^{1/3} \big\} \Big\}
 P \Big( \acon \Big),
\end{equation}
\end{lemma}
\noindent{\bf Proof.}
The proof is an application of the sector storage-replacement operation. We use the template provided by the proof of Proposition \ref{propmscbexc}, re-using notation with altered definitions when this is convenient. \\
\noindent{\bf Definition of satisfactory input, operation parameters and successful action.} 
Set ${\rm SAT}_1 = A_j$. Let $\bo{x^-},\bo{x^+} \in \reg$
be the first elements of $\reg$ encountered in a clockwise, respectively
counterclockwise, search from $\xmlrm$ or from $\xmlrp$.

From ${\rm SAT}_1  \subseteq \big\{ \cir \subseteq B_{\cctwo n} \big\}$, we may
find deterministic  $\bo{x_0},\bo{y_0} \in B_{\cctwo n}$ such that
$$
 P \Big(   \bo{x^-} = \bo{x_0}, \bo{x^+} = \bo{y_0}  \Big\vert {\rm SAT}_1 \Big) \geq \frac{1}{2 \pi^2 \cctwo^4 n^4}.
$$
The operation 
$\sigma_{\xo,\yo}$ will be taken to act regularly, in the sense of Definition \ref{defregact}.
Set ${\rm SAT}_2 = \big\{ \bo{x^-} = \bo{x_0}, \bo{x^+} = \bo{y_0}
\big\}$. 
The input 
$\omega \in \zoz$ will be said to be  
{\it satisfactory} if it realizes the event ${\rm SAT} : = {\rm
  SAT}_1 \cap {\rm SAT}_2$.
The operation will be said to {\it act successfully} 
if the updated configuration $\omega_1 \big\vert_{\axyoe}$ produced by
$\sigma_{\bo{x_0},\bo{y_0}}$ realizes the event ${\rm GAC}\big( \xo,\yo
\big)$ specified in Definition \ref{defgac}.\\
\noindent{\bf Output properties.}
We claim that, if the input configuration $\omega$ is satisfactory  and the operation $\sigma_{\bo{x_0},\bo{y_0}}$ acts successfully, then the output 
$\sigma_{\bo{x_0},\bo{y_0}}(\omega) = \big( \omega_1,\omega_2 \big)$ has the properties that:
\begin{enumerate}
\item {\bf Full-plane circuit property:} the full-plane configuration $\omega_1 \in \zoz$ contains an open
  circuit $\Gamma$ for which $\Gamma \subseteq B_{5 \cctwo n}$, 
and  
 $\Big\vert  {\rm INT}  \big( \Gamma \big) \Big\vert \geq n^2$;
\item {\bf Sector open-path property:}
 and the sector configuration $\omega_2 \in \big\{ 0,1 \big\}^{\axyoe}$
  realizes the event, to be denoted by ${\rm SOPP}^+$,
that there exists an open path $\gamma$ 
connecting $\xo$ to 
$\yo$ and such that $\{ \xo \} \cup \{ \yo \} \in V(\gamma)$,  
$\gamma 
\subseteq \axyo \cap \big( B_{\cctwo n} \setminus B_{\ccone n} \big)$,
$\gamma  \cap W_{\xo,Cn^{-1}\log n}^+ \subseteq \clu{\xo}$ and  
$\gamma \cap W_{\yo,Cn^{-1}\log n}^- \subseteq \clum{\yo}$, 
along with
\begin{equation}\label{dflucbd}
\sup \Big\{ d \big( \bo{v}, \big[ \xo,\yo \big] \big): \bo{v} \in
\gamma  \Big\} \geq n^{1/3} \big( \log n \big)^{2/3}  t - n^{1/6}.
\end{equation}
\end{enumerate}
\noindent{\bf Proof of the sector-open path property.} 
Note that ${\rm SOPP}^+ = {\rm SOPP} \cap (\ref{dflucbd})$, where ${\rm SOPP}$ is defined in sector-open path property that appears in the proof of Proposition \ref{propmscbexc}.
We define $\gamma = \cir \cap \axyo$ as in the earlier proof.
The previous argument shows that ${\rm SOPP}$ is satisfied in this case also.
To obtain (\ref{dflucbd}), note that from the form of $\gamma$
 that $\xmlr \in \gamma$, since
\begin{equation}\label{xoargyo}
\argu(\yo) =
\argu(\bo{x^+}) \geq
\argu(\xmlrp) \geq
\argu(\xmlr) \geq
\argu(\xmlrm) \geq
\argu(\bo{x^-}) =
\argu(\xo).
\end{equation}
We have that
\begin{equation}\label{xmlrnt}
d \Big( \xmlr, \big[ \xmlrm,\xmlrp \big] \Big) \geq
d \Big( \xmlr, \delconv \Big) = \mlr \geq n^{1/3} \big( \log n \big)^{2/3} t. 
\end{equation}
Note that
\begin{equation}\label{dnonesi}
\sup \Big\{ d \Big(  \bo{x}, \big[ \xmlrm,\xmlrp \big] \Big): \bo{x} \in \big[
\xo,\yo \big]  \Big\}
\leq \max \Big\{ d \big( \xo, \xmlrm \big), d \big( \yo, \xmlrp \big)
\Big\}
 \leq n^{1/6},
\end{equation}
the first inequality because any two planar line segments $\ell_1$ and
$\ell_2$ are such that $d\big(\bo{x},\ell_2\big)$ attains its supremum
among $\bo{x} \in \ell_1$ at one of its endpoints. The second
is due to the definition of  $\maxseprg$ and to $A_j \subseteq \big\{ \maxseprg \leq n^{1/6} \big\}$.

From (\ref{xmlrnt}) and (\ref{dnonesi}), we obtain 
\begin{equation}\label{xmlrxoyo}
d \big( \xmlr , [\xo,\yo] \big) \geq n^{1/3}\big( \log n \big)^{2/3} t -
n^{1/6}.
\end{equation}
By $\xmlr \in V\big( \gamma \big)$, we obtain (\ref{dflucbd}). \\
\noindent{\bf Proof of the full-plane circuit property.}
This proof follows that of its counterpart for Proposition \ref{propmscbexc}, except for a few changes. 
We replace $\bo{x}$ by $\xmlrm$, and $\bo{y}$ by $\xmlrp$. 
To verify (\ref{propnpl}) in the present case, note that 
$\max \big\{ \vert\vert \xo - \xmlrm \vert\vert, \vert\vert \yo - \xmlrp
\vert\vert  \big\} \leq \maxseprg \leq n^{1/6}$ by the occurrence of ${\rm
  SAT}_1 = A_j$, while  
$\vert\vert \xmlrm - \xmlrp \vert\vert = \mlrs = j \geq n^{1/3}$. Thus, 
$$
\max \Big\{ \vert\vert \xo - \xmlrm \vert\vert, \vert\vert \yo - \xmlrp
\vert\vert  \Big\} \leq \maxseprg \leq \vert\vert \xmlrm - \xmlrp \vert\vert^{1/2}, 
$$
which is (\ref{propnpl}) in the present setting. 
We obtain the counterpart of the third inequality in (\ref{intlbd}),
$$
\big\vert {\rm INT} \big( \Gamma \big)
\big\vert \geq 
\big\vert {\rm INT} \big( \Gamma_0 \big)
\big\vert +  \frac{1}{40} \vert\vert \xmlrm - \xmlrp \vert\vert^{3/2}
  \Big( \log \vert\vert \xmlrm - \xmlrp \vert\vert \Big)^{1/2}, 
$$
That $A_j \subseteq \acon$ yields 
$\big\vert {\rm INT} \big( \Gamma \big) \big\vert \geq n^2$. \\
\noindent{\bf The lower bound on the probability of satisfactory input and successful action.}
We claim that, as previously, this probability at least (\ref{prsatbd}).
Pursuing the earlier argument, we must in addition check that the
present choice of $\xo$ and $\yo$ satisfy the hypotheses of Lemma
\ref{lemgac}.
We obtain $\vert\vert \xo - \yo \vert\vert \geq C \log n$, (indeed,
$\vert\vert \xo - \yo \vert\vert \geq c n^{1/3}$),  from 
$\max \big\{ \vert\vert \xo - \xmlrm \vert\vert, \vert\vert \yo - \xmlrp
\vert\vert  \big\} \leq \vert\vert \xmlrm - \xmlrp \vert\vert^{1/2}$ 
and
$\vert\vert \xmlrm - \xmlrp \vert\vert = \mlrs = j \geq n^{1/3}$.
That $\yo \in \clu{\xo}$ and $\xo \in \clum{\yo}$ follows from $\xo,\yo \in
\reg$, since, as we now argue, $\ang \big( \xo,\yo \big) \leq c_0$. Indeed, 
\begin{equation}\label{xoyombd}
 \vert\vert \xo - \yo \vert\vert \leq \mlrs + 2 \maxseprg \leq j + 2 n^{1/6},
\end{equation}
the latter inequality by the occurrence of $A_j$. From $j \leq n^{2/3}
\big( \log n  \big)^{1/2} t^{2 - \delta}$ and $t = o(n^{1/6})$, (a
weaker condition than our assumption),
we see that $\vert\vert \xo - \yo \vert\vert = o(n)$. Thus,
$\cir \cap B_{\ccone n} = \emptyset$ and $\xo,\yo \in \vcir$ yield $\ang
\big( \xo, \yo \big) = o(1)$, as we sought (to verify the hypotheses of Lemma \ref{lemgac}). \\
\noindent{\bf The upper bound on the probability of the two output properties.} The upper bound (\ref{proponebd}) on the probability of the full-plane circuit property is derived without change. Likewise, the conditional probability of the sector open-path property, given the full-plane circuit property, is shown to be bounded above by $\conka P \big( {\rm SOPP}^+ \big)$. 
To bound $P \big( {\rm SOPP}^+ \big)$,
note that
the final
requirement of the event ${\rm SOPP}^+$ 
entails that the set $\gamma$ satisfies
$$
 \{ \xo \} \cup \{ \yo \}  \subseteq \vcir,  \, \,
 \sup \Big\{ d \Big( \bo{v}, \big[ \xo,\yo \big] \Big): \bo{v} \in
 V(\gamma) \Big\} \geq   \vert\vert \xo - \yo
\vert\vert^{1/2}  \frac{n^{1/3} \big( \log n \big)^{2/3} t}{2 \sqrt{2j}},
$$
because $n^{1/6} \leq n^{1/3}\big( \log n \big)^{2/3} t/2$ (since $t \geq
t_0 \geq 1$), and $j \geq
\vert\vert \xo - \yo \vert\vert/2$, (which is due to (\ref{xoyombd}) and 
 $j \geq n^{1/3} \geq 2 n^{1/6}$).

Applying Lemma \ref{lemmdf}, we find that
$$
 P \big( {\rm SOPP}^+ \big) \leq  \exp \bigg\{ - c  \min \Big\{ \frac{n^{2/3}\big(
   \log n \big)^{4/3} t^2}{8j}, c^2 \vert\vert \xo - \yo \vert\vert \Big\} \bigg\}
 P \Big( \xo \leftrightarrow \yo \Big).
$$
As noted before (\ref{xoyombd}), 
$\vert\vert \xo - \yo \vert\vert \geq c n^{1/3}$. Thus,
$$
 P \big( {\rm SOPP}^+ \big) \leq  \exp \bigg\{ - c  \min \Big\{ \frac{n^{2/3}\big(
   \log n \big)^{4/3} t^2}{8j}, c^2 c n^{1/3} \Big\} \bigg\}
 P \Big( \xo \leftrightarrow \yo \Big).
$$ 
Hence, in acting on an input having the law $P$, $\sigma_{\xo,\yo}$
will produce an output having the full-plane circuit and sector open-path 
properties with probability at most
\begin{equation}\label{upbda}
 \cpi \big( 5 \cctwo \big)^2  \conka n^2 P \Big( \acon \Big)  \exp \bigg\{ - c  \min \Big\{ \frac{n^{2/3}\big(
   \log n \big)^{4/3} t^2}{8j}, c^2 c n^{1/3} \Big\} \bigg\}
 P \Big( \xo \leftrightarrow \yo \Big).
\end{equation}
\noindent{\bf Conclusion by comparing the obtained bounds.}
Analogously to the proof of Proposition \ref{propmscbexc}, 
(\ref{upbda}) is at least (\ref{prsatbd}). 
Recalling that ${\rm SAT}_1 = A_j$, we obtain
$$
P \big( A_j \big) \leq n^{10^{-2} \clemgac + 7} 
 \frac{P \big( \xo \leftrightarrow \yo \big)}{P\big( \xo \build\leftrightarrow_{}^{\axyo} \yo \big)}
 \exp \bigg\{ - c  \min \Big\{ \frac{n^{2/3}\big(
   \log n \big)^{4/3} t^2}{8j}, c c^2 n^{1/3} \Big\} \bigg\}
   P \Big( \acon \Big),
$$
from which
$$
P \big( A_j \big) \leq c^{-1} n^{10^{-2} \clemgac + 7}  
 \exp \bigg\{ - c  \min \Big\{ \frac{n^{2/3}\big(
   \log n \big)^{4/3} t^2}{j}, n^{1/3} \Big\}
 \bigg\}
  P \Big( \acon \Big),
$$
because, due to an application of Lemma \ref{lemoznor}, we have that 
\begin{equation}\label{xyarbd}
 \frac{P\big( \xo \build\leftrightarrow_{}^{\axyo} \yo \big)}{P \big( \xo
   \leftrightarrow \yo \big)} \geq c.
\end{equation}
In applying Lemma \ref{lemoznor}, we take $\delta$ to be any positive value
less than $\qzero/2$. Indeed, as we have seen, 
$\yo \in \clu{\xo}$ and $\xo \in \clum{\yo}$. These imply that 
$W_{\xo - \yo,\qzero/2}\big( \yo \big) \cap W_{\yo - \xo,\qzero/2}\big( \xo \big)
\subseteq \axyo$. 
Hence, by Lemma \ref{lemoznor}, 
we reach the conclusion (\ref{xyarbd}), with $\axyo \cup
B_K(\xo) \cup B_K(\yo)$ in place of $\axyo$ in the numerator. Finally, the
desired form may be obtained, by altering the configuration inside
$B_K(\xo) \cup B_K(\yo)$, since two configurations differing on $N$ edges have conditional probabilities differing by a factor of at most $\cposen^{-N}$, by the bounded energy property of $P$.

From $n^{1/3} \leq j \leq n^{2/3} \big( \log n \big)^{1/3} t^{2 - \delta}$, we conclude that the statement of the lemma holds. \qed
Note that
\begin{eqnarray}
 &  &     P \Big(  \mlrs \geq n^{2/3} \big( \log n \big)^{1/3} t^{2 - \delta} , \acon \Big)
 \label{mlrsbd} \\
 &  \leq &     P \Big(  \mfl \geq n^{2/3} \big( \log n \big)^{1/3}  t^{2 - \delta} , \acon \Big)
 \nonumber \\
 & \leq &    n^C
  P \Big(  \big\vert \intg \big\vert \geq  n^2 +  c n t^{\frac{3}{2}(2 -
    \delta)} \log n 
  \Big) \nonumber \\
  &  & \qquad + \, \, \,  
 \exp \Big\{  - c n^{1/3} \big( \log n \big)^{1/6} 
   t^{1 - \delta/2}  \Big\} P \Big( \acon \Big)
 \nonumber \\
  & \leq & \bigg(    
  n^C \exp \Big\{  - c t^{\frac{3}{2}(2 - \delta)} \log n \Big\}
 + \exp \Big\{ - c n^{1/3}  \big( \log n \big)^{1/6}   t^{1 - \delta/2}  \Big\} 
 \bigg)  P \Big( \acon \Big).
 \nonumber 
\end{eqnarray}
The second inequality here is due to Proposition \ref{propmscbexc} and
requires that 
\begin{equation}\label{ttwodo}
t^{2 - \delta} = o \big( n^{4/3}  \big)
\end{equation}
The third is due to Proposition \ref{propexc} and
requires that 
\begin{equation}\label{tthrtwc}
t^{(3/2)(2 - \delta)} \leq cn,
\end{equation}
as well as $t \geq t_0 \geq C$.

We will show that $A = \emptyset$.
Noting that (\ref{xmlrxoyo}) relied on $\maxseprg \leq n^{1/6}$, which holds in
the present case also, we find that
\begin{equation}\label{xmxmineq}
 d \Big( \xmlr, \big[ \bo{x^-},\bo{x^+} \big] \Big)
 \geq   n^{1/3} \big( \log n \big)^{2/3}  t - n^{1/6}.
\end{equation}
The conditions that define the event $A$ ensure that 
$\vert\vert \bo{x^-} \vert\vert , \vert\vert \bo{x^+} \vert\vert \geq
\ccone n$
and $\vert\vert \bo{x^-} - \bo{x^+} \vert\vert \leq n^{1/3} + 2n^{1/6}$,
and these two lead to $\ang \big( \bo{x^-}, \bo{x^+} \big) = o(1)$. From
this, and $\bo{x^-} \in \reg$, we learn that
$W_{\bo{x^-},c_0}^+ \cap \cir   
\subseteq \clu{\bo{x^-}}$.
We use this as well as $\ang \big( \bo{x^-}, \bo{x^+} \big) \leq c_0$ and
$\xmlr \in A_{\bo{x^-},\bo{x^+}}$ to obtain 
$\xmlr \in \clu{\bo{x^-}}$. 
It follows from $\ang\big( \bo{x^-},\bo{x^+} \big) = o(1)$ and Lemma \ref{lemdistang}
that 
$$
  {\rm diam}  \Big(  \clu{\bo{x^-}}  \cap
   A_{\bo{x^-},\bo{x^+}} \Big) \leq \csc(\qzero/2) \vert\vert \bo{x^-} - \bo{x^+} \vert\vert.
$$
Hence,
\begin{eqnarray}
 & & d \Big( \xmlr, \big[ \bo{x^-},\bo{x^+} \big] \Big)
 \leq  d \big( \xmlr , \bo{x^+} \big) \nonumber \\
 & \leq &   {\rm diam}  \Big(  \clu{\bo{x^-}} \cap
   A_{\bo{x^-},\bo{x^+}} 
 \Big)  \leq  \csc(\qzero/2) \vert\vert \bo{x^-} - \bo{x^+} \vert\vert \leq 2 \csc(\qzero/2) n^{1/3}, \nonumber 
\end{eqnarray}
contradicting (\ref{xmxmineq}) and proving that $A = \emptyset$.

By Theorem \ref{thmmaxrg} and Lemma \ref{lemmarbd},
\begin{equation}\label{maxseprgbd}
 P \Big( \maxseprg >  n^{1/6}, \cir \subseteq B_{\cctwo n}, \acon \Big) \leq \exp \Big\{ - c n^{1/6} \Big\}
  P \Big( \acon \Big).
\end{equation}
By means of the inclusion (\ref{mlrinc}), and the deductions (\ref{mlrsbd}), (\ref{probajbd}), $A =
\emptyset$ and (\ref{maxseprgbd}), as well as Lemma \ref{lemmac},
\begin{eqnarray}
 & &  P \Big(  \mlr \geq n^{1/3} \big( \log n \big)^{2/3} t , \acon \Big) \nonumber \\
 & \leq & 
 \bigg(    
  n^C \exp \Big\{  - c t^{\frac{3}{2}(2 - \delta)} \log n \Big\}
 + \exp \Big\{ - c n^{1/3}  \big( \log n \big)^{1/6}   t^{1 - \delta/2}  \Big\} 
 \bigg)  P \Big( \acon \Big) \nonumber \\
 & &  \qquad + \, \, \, \sum_{j = n^{1/3}}^{ n^{2/3} (\log n)^{1/3} t^{2 - \delta}}
  \exp \Big\{  - c \min \big\{  t^{\delta} \log n, n^{1/3} \big\} \Big\}
   P \Big( \acon \Big)  \nonumber \\ 
& &  \qquad + \, \, \,
 \bigg(
\exp \Big\{  - c n^{1/6} \Big\} + \exp \big\{ - c  n \big\}  \bigg)   P \Big( \acon \Big)  \nonumber \\ 
& \leq &   \bigg(  \exp \Big\{  -  c t^{3 - \frac{3}{2}\delta} \log n
\Big\} +  n^{2/3} \big( \log n \big)^{1/3} t^{2 - \delta} \exp \Big\{ -
  c  \min \big\{  t^{\delta} \log n, n^{1/3} \big\} \Big\} \nonumber \\
  & & \qquad \qquad \qquad
 +  
\exp \Big\{  - c n^{1/6} \Big\}  \bigg)   P \Big( \acon \Big),  \nonumber
\end{eqnarray}
the first inequality requiring (\ref{ttwodo}) and (\ref{tthrtwc}), 
and the second, $t \geq t_0$.
Choosing $\delta = 6/5$, we obtain, for $t \geq t_0$,
$t = O \big( n^{5/18} (\log n)^{-C} \big)$, 
\begin{eqnarray}
 & &  P \Big(  \mlr \geq n^{1/3} \big( \log n \big)^{2/3} t  \Big\vert \acon \Big) \nonumber \\ 
 & \leq & 2 n^{2/3} \big( \log n \big)^{1/3} t^{2 - \delta} \exp \Big\{ - c t^{6/5}
 \log n \Big\} +    
\exp \Big\{  - c n^{1/6} \Big\}. \nonumber 
\end{eqnarray}
The condition $t = O \big( n^{5/18} (\log n)^{-C} \big)$ ensures that 
$t^{6/5} \log n \leq n^{1/3}$. We also require that  (\ref{ttwodo}) and
(\ref{tthrtwc}), and this condition ensures them. 
We obtain, for  $t \geq t_0$,
$t = O \big( n^{5/36} (\log n)^{-C} \big)$,
$$
 P \Big(  \mlr \geq n^{1/3} \big( \log n \big)^{2/3} t  \Big\vert \acon \Big)
 \leq  
 \exp \Big\{ - c t^{6/5} \log n    \Big\}.
$$
An application of Lemma \ref{lemcendisp} completes the proof. \qed
\begin{subsection}{The proof of Theorem \ref{thmfixedarea}}\label{secfixedarea}
This is a simple consequence of the regeneration structure Theorem \ref{thmmaxrg}. We only sketch the details. The reader who wishes to construct a complete proof is invited to consult Proposition 
3
in \cite{hammondthr}, where the surgery used here is undertaken to prove a lower bound on the area-excess $\exc$ (with a slight change in the choice of parameters).

It suffices to show that $P \big( \big\vert \intg \big\vert = n^2 \big) \geq n^{-C} P \big( \areacon \big)$. Note that Proposition \ref{propglobdis} implies  
$P \big( \big\vert \intg \big\vert \geq n^2, \centre(\cir) = \bo{0} \big) \geq c n^{-2} P \big( \areacon \big)$.
Defining $p_k =  P \big(  \big\vert \intg \big\vert = k, \centre(\cir) = \bo{0} \big)$,
it thus suffices to show that $p_{n^2} \geq n^{-C} P \big( \acon \big)$.
Fixing $n \in \N$, call $m \in \N$ good if   $m \in \{0,\ldots,C n \log n \}$ and $p_{n^2 - m} \geq n^{-C} P \big( \acon \big)$. 
By the area-excess Proposition \ref{propexc}, $\sum_{m \, \textrm{good}}{p_{n^2 - m}} \geq 2^{-1} P \big(  \acon \big)$, so that it is enough to show that
\begin{equation}\label{pnmpn}
p_{n^2} \geq n^{-C} p_{n^2 - m}
\end{equation}
for all good $m$. Given the regeneration structure Theorem \ref{thmmaxrg}, this is 
easily accomplished by a fairly local surgery. Fix $c'$ small and write $m = c' n k + r$ with $0 \leq k \leq C \log n$ and $r = \Theta(\log n)$. Note that the statement of Theorem \ref{thmmaxrg} holds for the measure 
$P \big( \cdot \big\vert \intg \big\vert = n^2 - m, \centre(\cir) = \bo{0}   \big)$ because of the lower bound on $p_{n^2 - m}$. 
Under a typical realization of $P \big( \cdot \big\vert \intg \big\vert = n^2 - m, \centre(\cir) = \bo{0}   \big)$, fix $\bo{x}$ to be the element of $\reg$ in the right-hand half-plane closest to the $x$-axis. Moving counterclockwise along $\cir$, the $y$-coordinate of the encountered elements of $\reg$ increases in steps of size $\Theta(\log n)$. Let $\bo{y}$ be the first such element having a height difference with $\bo{x}$ of at least $c' n$. 
Now, the configuration between $\bo{x}$ and $\bo{y}$  may be shifted to the right by $k$ places, with the configuration in $A_{\bo{x},\bo{y}}^c$ held fixed and the remaining region updated, with this random map sending $P$ to an output
that is absolutely continuous with respect to $P$. (Formally, we set $\xo,\yo$ to be a deterministic pair of points with a polynomially decaying probability of coinciding with $\bo{x}$ and $\bo{y}$, then
apply the shift-replacement operation $\xi_{A,B,\bo{z}}$ specified at the end of Definition \ref{defsro} with $F = A_{\xo,\yo}^c$, $G = \clu{\bo{x}} \cap \clum{\bo{y}}$ and $\bo{z} = k \bo{e_1}$.) Action is defined to be  successful if the two horizontal paths of length $k$ along which the points $\xo$ and $\yo$ were moved are open under the updated configuration.  
The  operation is successful (due to $\xo$ and $\yo$ hitting $\bo{x}$ and $\bo{y}$, and the paths being open) with a polynomially decaying probability in $n$, and it maps $\big\{   \big\vert \intg \big\vert =n^2 - m, \centre(\cir) = \bo{0}  \big\}$ into  $\big\{ \big\vert \intg \big\vert = n^2  + t \big\}$, with $t = \Theta(\log^2 n)$. We may then excise an area-$t$ region trapped inside but touching the circuit at polynomial cost in $n$, since such a region may be chosen with a boundary of length $\Theta\big(\log n\big)$.
The output circuit may be recentred at $\bo{0}$ by a shift. Thus, we obtain (\ref{pnmpn}). 
\qed
\end{subsection}
\end{section}
\bibliographystyle{plain}
\bibliography{mlrbib}
\end{document}